\documentclass[smallextended]{svjour3}       
\smartqed  

\usepackage{amsfonts,amsmath,amssymb,amstext,graphicx}

\usepackage{mathrsfs}
\usepackage[margin=1in]{geometry}
\RequirePackage[OT1]{fontenc}
\RequirePackage{natbib}
\RequirePackage
{hyperref}
\RequirePackage{hypernat}

\font\msbmx=msbm10                   
\font\msbmvii=msbm7                  
\font\msbmv=msbm5
\def\varnothing{\mathchoice{\mbox{\msbmx\char'077}}%
{\mbox{\msbmx\char'077}}{\mbox{\msbmvii\char'077}}{\mbox{\msbmv\char'077}}}%

\newcommand{\si}{\sigma}

\newcommand{\vae}{\varepsilon}

\newcommand{\wt}{\widetilde}

\def\bbn{{\mathbb N}}

\def\bbr{{\mathbb R}}

\newcommand{\mb}[1]{\mathbf{#1}} 
\newcommand{\Qb}{{\mb{Q}}}

\newcommand{\Xb}{{\mb{X}}}

\newcommand{\mc}[1]{\mathcal{#1}} 
\newcommand{\Rc}{{\mc{R}}}

\newcommand{\Nc}{{\mc{N}}}
\newcommand{\Fc}{{\mc{F}}}
\newcommand{\Gc}{{\mc{G}}}
\newcommand{\Hc}{{\mc{H}}}

\newcommand{\Mcc}{{\mc{M}}}

\newcommand{\Xc}{{\mc{X}}}
\newcommand{\Bc}{{\mc{B}}}

\def\text#1{\hbox{#1}}
\def\proof{{\noindent \bf Proof. }}
\def\mes{{\bf mes}}
\def\endproof{\mbox{\ $\qed$}}

\def\E{{\bf E}}

\def\C{{\bf C}}

\def\D{{\bf D}}
\def\B{{\bf B}}
\def\M{{\bf M}}
\def\A{{\bf A}}
\def\U{{\bf U}}
\def\H{{\bf H}}
\def\V{{\bf V}}

\def\G{{\bf G}}

\def\r{{\bf r}}
\def\x{{\bf x}}
\def\a{{\bf a}}
\def\b{{\bf b}}

\def\l{{\bf l}}
\def\v{{\bf v}}
\def\u{{\bf u}}
\def\m{{\bf m}}

\def\d{\mathrm{d}}
\def\build #1_#2{\mathrel{\mathop{\kern 0pt #1}\limits_\zs{#2}}}
\newcommand{\zs}[1]{{\mathchoice{#1}{#1}{\lower.25ex\hbox{$\scriptstyle#1$}}
{\lower0.25ex\hbox{$\scriptscriptstyle#1$}}}}

\numberwithin{equation}{section}
\def\proof{{\noindent \bf Proof. }}
\def\endproof{\mbox{\ $\qed$}}

\usepackage{mathrsfs}
\newcommand{\Pb}{{\mathsf{P}}} 
\newcommand{\EV}{{\mathsf{E}}} 
\newcommand{\PFA}{\mathsf{PFA}}
\newcommand{\LPFA}{\mathsf{LPFA}}
\newcommand{\LCPFA}{\mathsf{LCPFA}}
\newcommand{\ESADD}{{\mathsf{ESADD}}} 
\newcommand{\SADD}{{\mathsf{SADD}}} 
\newcommand{\Hyp}{{\mathsf{H}}} 

\newcommand{\mrm}[1]{\mathrm{#1}}

\newcommand{\F}{{\mrm{F}}}

\def\One{\mathchoice{\rm 1\mskip-4.2mu l}{\rm 1\mskip-4.2mu l}
{\rm 1\mskip-4.6mu l}{\rm 1\mskip-5.2mu l}}
\newcommand\Ind[1]{{\One_{\{#1\}}}} 

\newcommand{\xra}{\xrightarrow} 

\newcommand{\abs}[1]{\left\vert#1\right\vert}
\newcommand{\set}[1]{\left\{#1\right\}}

\newcommand{\brc}[1]{\left(#1\right)}
\newcommand{\brcs}[1]{\left[#1\right]}

\renewcommand{\le}{\leqslant} 
\renewcommand{\ge}{\geqslant}
\newcommand{\esssup}{\operatornamewithlimits{ess\,sup}}

\newcommand{\ignore}[1]{} 

\usepackage{color}

\begin{document}

\title{Asymptotically optimal pointwise and minimax quickest change-point detection for dependent data}  


\author{Serguei Pergamenchtchikov \and
Alexander G. Tartakovsky} 

\institute{Serguei Pergamenchtchikov \at
UMR 6085 CNRS-Universite de Rouen, France\\
\email{serge.pergamenchtchikov@univ-rouen.fr}
 \and
 Alexander G. Tartakovsky \at
 Somers, CT USA
 \email{alexg.tartakovsky@gmail.com}
 }

\date{Received: date / Accepted: date}

\maketitle


\begin{abstract}
We consider the  quickest change-point detection problem in pointwise and minimax settings for general dependent data models.
Two new classes of sequential detection procedures associated with the maximal ``local'' probability of a false alarm within a period of some fixed length  
are introduced. 
For these classes of detection procedures, we consider  two popular risks: the expected positive part of the delay to detection and the conditional delay to detection. 
Under very general conditions for the observations, we show that the popular Shiryaev--Roberts  procedure is asymptotically optimal, as the local probability of false alarm goes to zero,
with respect to both these risks pointwise (uniformly for every possible point of change) 
and in the minimax sense (with respect to maximal over point of change expected detection delays). 
The conditions are formulated in terms of the rate of convergence in the strong law of large numbers for the log-likelihood ratios between the ``change'' and ``no-change'' hypotheses, specifically as a 
uniform complete convergence of the normalized log-likelihood ratio to a positive and finite number. We also develop tools and a set of sufficient conditions for verification 
of the uniform complete convergence for a 
large class of Markov processes. These tools are based on concentration inequalities for functions of Markov processes and the Meyn--Tweedie geometric ergodic theory. 
Finally, we check these sufficient conditions for a number of challenging examples (time series) frequently arising in applications,  such as autoregression, autoregressive GARCH, etc.
\end{abstract}



%

%
\section{Introduction} \label{sec:intro}
The problem of rapid detection of abrupt changes in a state of a process or a system arises in a variety
of applications from engineering problems (e.g., navigation integrity monitoring \cite{BassevilleNikiforovbook93, TNBbook14}), military applications (e.g., target detection and tracking in heavy clutter
\cite{Tarbook91, TNBbook14}) to cyber security (e.g., quick detection of attacks in computer networks
\cite{Kent, Tartakovskyetal-SM06, Tartakovskyetal-IEEESP06,Tartakovsky-Cybersecurity13, TNBbook14}). 
In the present paper, we are interested in a sequential setting assuming that 
as long as the behavior of the observation process is consistent with a ``normal'' (initial in-control) state, we allow the process to continue. If the state changes, then we need to detect this event as rapidly as possible while 
controlling for the risk of false alarms. In other words, we are interested in designing the quickest change-point detection procedure that optimizes
the tradeoff between a measure of detection delay and a measure of the frequency of false alarms.

There are four conventional approaches to the optimum tradeoff problem: Bayesian, generalized Bayesian, multicyclic detection of changes in a stationary regime, and minimax (see \citet[Ch 6]{TNBbook14}).  
In the Bayesian context, proposed by~\citet{girshick-ams52} and \citet{Shiryaev61,ShiryaevTPA63}, the change point is assumed to be random with
a geometric prior distribution, and the optimality criterion is to minimize the weighted Bayes-type expected detection delay subject to an upper bound on the weighted probability of a false alarm. 
Until the 1990s, most of the work related to the optimality issue in change detection had been done in the iid case, assuming that observations are independent and identically distributed (iid) with 
one law before the change 
and with another distribution after the change. In particular, in the 1960s, \citet{Shiryaev61,ShiryaevTPA63} found an optimal Bayes solution showing that a detection procedure 
based on thresholding the posterior probability of the change up to the current moment is strictly optimal for any value of the weighted false alarm probability. 
Much later, in 2004--2006, a general Bayesian asymptotic theory of change-point detection (for very general non-iid models and arbitrary prior distributions of the change point) was 
developed by~\citet{TV2004} in discrete time and~\citet{BT2006} in continuous time. 

By contrast, in a minimax formulation, proposed by~\citet{LordenCPDAS71} and \citet{PollakAS85}, the change point is assumed to be an unknown non-random number and the goal 
is to minimize the worst-case delay (with respect to the point of change) subject to a lower bound on the mean time until false alarm.   Specifically,  in 1971,
\citet{LordenCPDAS71} suggested the worst-worst-case average delay to detection measure $\ESADD(\tau) = \sup_{\nu\ge 0} \esssup \EV_\nu(\tau - \nu | \tau > \nu, \Fc_\nu)$ 
that should be minimized in the class of procedures $\Hc_\gamma=\{\tau: \EV_\infty \tau \ge \gamma\}$ for which the average run length (mean time) to false alarm $\EV_\infty \tau$ is not smaller than a 
given number $\gamma>1$. Here $\tau$ is a generic change detection procedure (stopping time), $\EV_\nu$ stands for the operator of expectation when the change point is $\nu$ ($\nu=\infty$ corresponds 
to a no-change scenario) and $\Fc_\nu=\sigma(X_1,\dots,X_\nu)$ is
the sigma-algebra generated by the first $\nu$ observations $X_1,\dots,X_\nu$. \citet{LordenCPDAS71} developed an asymptotic minimax theory of change detection (in the iid case) as $\gamma \to \infty$, 
proving in particular that Page's CUSUM procedure \cite{PageBiometrika54} is asymptotically first-order minimax. Later in 1986, \citet{MoustakidesAS86} established strict optimality of CUSUM for any value of 
the average run length to false alarm $\gamma >1$. In the 1980s, \citet{PollakAS85} introduced a less pessimistic worst-case detection delay measure ---
maximal conditional average delay to detection, 
\begin{equation} \label{SADDdef}
\SADD(\tau) = \sup_{\nu\ge 0} \EV_\nu(\tau - \nu | \tau > \nu) ,
\end{equation}
and found an almost optimal procedure that minimizes $\SADD(\tau)$ subject
to the constraint on the average run length to false alarm (i.e., in the class $\Hc_\gamma$) as $\gamma$ becomes large. Pollak's idea was to modify the Shiryaev--Roberts
statistic by randomization of the initial condition in order to make it an equalizer. Pollak proved that the randomized Shiryaev--Roberts procedure that starts from a random point sampled from the 
quasi-stationary distribution of the Shiryaev--Roberts statistic is asymptotically nearly minimax within an additive vanishing term. Since the Shiryaev--Roberts--Pollak
procedure is an equalizer, it is tempting to conjecture that it may be strictly
optimal  for any value of $\gamma$, which is not true, as the articles of \citet{MoustPolTarSS11}
and \citet{PolunTartakovskyAS10} indicate.  

As we already mentioned above, in the early stages the theoretical development was focused primarily on the iid case.
However, in practice the observations may be non-identically distributed and dependent.
A general asymptotic minimax theory of change-point detection for non-iid models was developed by \citet{LaiJRSS95,LaiIEEE98} (see also \citet{Fuh03} for hidden Markov models with a finite state-space). 
In particular, for a low false alarm rate (large $\gamma$) the asymptotic minimaxity of the CUSUM procedure was established in \cite{Fuh03,LaiIEEE98}. 

In the iid case, the suitably standardized distributions of the stopping times of the CUSUM and Shiryaev--Roberts detection procedures
are asymptotically exponential for large thresholds and fit well into the geometric distribution
even for a moderate false alarm rate (see \citet{PollakTartakovskyTPA09}). In this case, the average run length 
to false alarm is an appropriate measure of false alarms. However,
for non-iid models the limiting distribution is not guaranteed to be exponential or even close to
it. In general, we cannot even guarantee that large values of the average run length to false alarm
will produce small values of the maximal local false alarm probability. Therefore, the average run length 
to false alarm is not appropriate in general, and instead it is more adequate to use the local conditional false alarm probability, as suggested in~\cite{TartakovskyIEEECDC05,TNBbook14}.  This issue is extremely
important for non-iid models, as a discussion in \cite{Mei-SQA08,Tartakovsky-SQA08a} shows.

In the present paper, we pursue two objectives. First, in Section~\ref{sec:PF}, we introduce two novel classes of change-point detection procedures, which, instead of imposing a lower bound on the average run 
length to false alarm, 
require more adequate upper bounds on the uniform probability of false alarm or uniform conditional probability of false alarm in the spirit of works by \citet{LaiIEEE98}, \citet{TartakovskyIEEECDC05} 
and \citet{TNBbook14}. However, these classes slightly differ from those proposed in \cite{LaiIEEE98, TartakovskyIEEECDC05,TNBbook14}. This modification 
allows us to substantially relax Lai's essential supremum 
conditions \cite{LaiIEEE98}, which do not hold for certain interesting practical models. In fact, our conditions are equivalent to the uniform 
version of the complete convergence for the log-likelihood ratio processes, i.e., they are related to the rate of convergence in the strong law of large numbers for the log-likelihood ratio between 
the ``change'' and ``no-change'' hypotheses.
We concentrate on a minimax problem of minimizing Pollak's maximal conditional average delay to detection defined in \eqref{SADDdef} as well as on a pointwise problem of minimizing the
conditional average delay to detection $\EV_\nu(\tau-\nu | \tau>\nu)$ for every change point $\nu\ge 0$. For the sake of completeness, 
we also consider the other popular risks $\sup_{\nu\ge 0} \EV_\nu(\tau-\nu)^+$ and $\EV_\nu(\tau-\nu)^+$, $\nu\ge 0$, while we strongly believe that the conditional versions $\EV_\nu(\tau-\nu | \tau>\nu)$ and
\eqref{SADDdef} are more appropriate for most applications. 
We consider extremely general non-iid stochastic models  for the observations, and it is our goal to find reasonable sufficient conditions for the observation models under which the Shiryaev--Roberts 
(or CUSUM) procedure is asymptotically optimal. 
To achieve the first goal we exploit the asymptotic Bayesian theory of change-point detection developed by \citet{TV2004} that offers a constructive and flexible approach for 
studying asymptotic efficiency of 
Bayesian type procedures. It turns out that a similar method can be used for the analysis of minimax risks and that the complete convergence type conditions for the log-likelihood ratio processes proposed 
in~\cite{TV2004,TNBbook14} are also sufficient in the minimax setting. These sufficient conditions as well as the main results related to asymptotic optimality of the Shiryaev--Roberts procedure in the classes of 
procedures with upper bounds 
on the weighted false alarm probability and local false alarm probabilities are given, correspondingly,  in Section~\ref{sec:Bay} and Section~\ref{sec:MaRe}.  

The second objective is to find a method for verification of the required sufficient conditions in a number of particular, still very general, challenging models. The natural question is how one may check the proposed 
sufficient conditions
and even  whether there are more or less general models, except of course the iid case, for which these conditions hold. To this end, we focus on the class of data models 
for which one can exploit the method of geometric ergodicity for homogeneous Markov processes, first proposed by \citet{MeTw93} and then further developed by \citet{GaPe13, GaPe14}
for statistical applications. These results are presented in Section~\ref{sec:Mrk} and show that our sufficient conditions for pointwise and minimax optimality hold for  homogeneous Markov ergodic processes.  
In Section~\ref{sec:Ex}, these 
conditions are further illustrated for several examples that include autoregressive, autoregressive GARCH, and other models widely used in many applications, in particular for modeling of dynamics of financial 
indices; see, e.g., \citet{ShiryaevSPR06}. All auxiliary results needed for the proofs are presented in Appendix~\ref{A}, and in Appendix~\ref{B} we give certain useful results from the geometric 
ergodic theory of Markov processes.

\section{Notation, problem formulation and detection procedures} \label{sec:PF}

Assume that we are able to observe a series of consecutive random variables $X_1,X_2, \dots$, which may change statistical properties at an unknown point in time $\nu\in\{0,1, 2, \dots\}$. 
We use the convention that $X_\nu$ is the {\em last pre-change} observation. 
Write $\Xb^{n}=(X_\zs{1},\dots,X_{n})$ for the concatenation of the first $n$ observations. Let $p_\nu(\Xb^n)= p(\Xb^n|\nu)$ be the joint probability density of the vector 
$\Xb^n$ when the change point $\nu$ is fixed and finite and let $p_\infty(\Xb^n)= p(\Xb^n|\nu=\infty)$ stand for the pre-change joint density (when the change never occurs). 
Let $\{f_\zs{0,n}(X_\zs{n}|\Xb^{n-1})\}_{n\ge 1}$ and $\{f_\zs{1,n}(X_\zs{n}|\Xb^{n-1})\}_{n\ge 1}$ be two sequences of conditional densities of $X_n$ 
given $\Xb^{n-1}$ ($f_{i,1}(X_1|X_0) \equiv f_{i,1}(X_1)$), $i=0, 1$) with respect to some non-degenerate sigma-finite measure $\mu(x)$. We are interested in the general non-iid case that
\begin{equation}\label{noniidmodel}
\begin{aligned}
p_\nu(\Xb^n) & = p_\infty(\Xb^n) = \prod_{i=1}^n f_{0,i}(X_i|\Xb_1^{i-1}) \quad \text{for}~~ \nu \ge n ,
\\
p_\nu(\Xb^n) & =  \prod_{i=1}^{\nu}  f_{0,i}(X_i|\Xb^{i-1}) \times \prod_{i=\nu+1}^{n}  f_{1,i}(X_i|\Xb^{i-1})  \quad \text{for}~~ \nu < n.
\end{aligned}
\end{equation}
In other words,  $\{f_\zs{0,n}(X_\zs{n}|\Xb^{n-1})\}_{n\ge 1}$ and $\{f_\zs{1,n}(X_\zs{n}|\Xb^{n-1})\}_{n\ge 1}$ are the pre-change and post-change conditional densities, respectively, so that if the change occurs
at time $\nu=k$, then the conditional density of the $(k+1)$-th observation changes from $f_{0,k+1}(X_{k+1}|\Xb^{k})$ to $f_{1,k+1}(X_{k+1}|\Xb^{k})$. Note that the post-change densities may depend on 
the change point~$\nu$, i.e., $f_{1,n}(X_n|\Xb^{n-1})= f_{1,n}^{(\nu)}(X_n|\Xb^{n-1})$ for $n > \nu$. We omit the superscript $\nu$ for brevity. 

Let $\Pb_k$ and~$\EV_k$ denote the probability and expectation when $\nu=k<\infty$, and let $\Pb_\infty$ and~$\EV_\infty$
denote the same when there is no change, i.e., $\nu=\infty$. Obviously, the general non-iid model given by \eqref{noniidmodel} implies that under the measure~$\Pb_\infty$  
the conditional density of~$X_n$ given $\Xb^{n-1}$
is $f_{0,n}(X_n|\Xb^{n-1})$ for all $n \ge 1$ and under~$\Pb_{k}$, for any $0\le k<\infty$, the conditional density of~$X_n$ is
$f_{0,n}(X_n|\Xb^{n-1})$  if $n \le k$ and is $f_{1,n}(X_n|\Xb^{n-1})$ if $n > k$. 

In the particular iid case, the observed random variables $X_\zs{1}, X_\zs{2},\dots$ are iid until a change with a common density $f_0(x)$ and after the change occurs,
the observations are again iid, but with another density $f_1(x)$. Therefore, in this case, the conditional densities $f_{0,i}(X_i|\Xb^{i-1})$ and $ f_{1,i}(X_i|\Xb^{i-1})$ in \eqref{noniidmodel} 
are replaced by $f_0(X_i)$ and $f_1(X_i)$, respectively.

A sequential detection procedure is a stopping (Markov) time $\tau$ for an observed sequence $\{X_\zs{n}\}_{n\ge 1}$, i.e.,  $\tau$ is an
extended integer-valued random variable, such that the event $\{\tau \le n\}$ belongs to the sigma-algebra $\Fc_\zs{n}=\sigma(X_\zs{1},\dots,X_\zs{n})$.  We denote by $\Mcc$ the set of all stopping times.
A false alarm is raised whenever the detection is declared before the change occurs, i.e.,  when $\tau\le \nu$. (Recall that $X_{\nu+1}$ is the first post-change observation.)
The goal of the quickest change-point detection problem is to develop a detection procedure that guarantees a stochastically small delay to detection 
$\tau-\nu$ provided that there is no false alarm (i.e., $\tau> \nu$) under a given (typically low) risk of false alarms.

Let $\Pb_k^{(n)}=\Pb_k |_{\Fc_n}$ denote a restriction of the probability measure $\Pb_k$ to the sigma-algebra $\Fc_n$. Then the likelihood ratio between the hypotheses 
``$\Hyp_k: \nu=k$'' that the change happens at $k<\infty$ and ``$\Hyp_\infty: \nu=\infty$'' that there is never a change (i.e.,  the Radon--Nikod\'{y}m density $\d \Pb_\zs{k}^{(n)}/\d \Pb_\zs{\infty}^{(n)}$) 
can be represented in the following exponential form
\begin{equation} \label{sec:Prbf.7}
\frac{\d \Pb_\zs{k}^{(n)}}{\d\Pb_\zs{\infty}^{(n)}}(\Xb^n)= 
e^{Z_\zs{n}^{k}}\,,
\end{equation}
where for $k\le n-1$
$$
Z_\zs{n}^{k} = \sum_\zs{j=k+1}^{n}\, \log \frac{f_\zs{1,j}(X_\zs{j}|\Xb^{j-1})}{f_\zs{0,j}(X_\zs{j}|\Xb^{j-1})}\,.
$$
The process $(Z_\zs{n}^{k})_{n \ge k+1}$ is the log-likelihood ratio (LLR) process between the hypotheses $\Hyp_k$ ($k=0,1, \dots$) and $\Hyp_\infty$.

In this paper, we study the Shiryaev--Roberts  (SR) procedure given by the following stopping time
\begin{equation}\label{sec:Prbf.7-00}
T(h)=\inf\left\{n\ge 1\,:\,
\sum_\zs{k=1}^{n}\,e^{Z^{k-1}_\zs{n}}
\ge h\right\}
\,,
\end{equation}
where $h>0$ is some fixed positive threshold which will be specified later. We set $\inf\{\varnothing\}=+\infty$. In the iid case, this procedure has certain interesting strict optimality properties 
(see \citet{PollakTartakovskySS09} and \citet{TNBbook14}).

Another popular change detection procedure is the CUSUM procedure given by the stopping time
\[
T_{\rm {CS}}(a)=\inf\left\{n\ge 1\,:\, \max_\zs{1\le k \le n} \,Z^{k-1}_\zs{n}\ge a\right\}\,, \quad a>0.
\]
 It may be shown that this procedure has essentially the same asymptotic performance as the SR procedure. In fact, using essentially the same line of argument, it can be proved that both procedures 
 are first-order  asymptotically optimal under the same general conditions. For this reason, we consider only the SR procedure.
 
Our main goal is to show that the SR detection procedure $T(h)$ is nearly optimal in two pointwise and minimax problems described below.  
We will also show that this procedure is asymptotically pointwise and minimax optimal in a class of Bayes-type procedures (see Section~\ref{sec:Bay}).  

To describe these problems we introduce for any $0<\beta<1$,  $m^{*}\ge 1$
 and $k^{*}>m^{*}$ the following classes of change detection procedures
\begin{equation}\label{sec:Prbf.3}
\Hc(\beta, k^{*},m^{*})=\left\{\tau\in\Mcc: \sup_\zs{1\le k\le k^{*}-m^{*}}\, \Pb_\zs{\infty}(k\le \tau < k+m^{*}) \le \beta \right\}
\end{equation}
and
\begin{equation}\label{sec:Prbf.4}
\Hc^{*}(\beta, k^{*},m^{*})=\left\{\tau\in\Mcc: \sup_\zs{1\le k\le k^{*}-m^{*}}\, \Pb_\zs{\infty}(\tau < k+m^{*} | \tau > k) \le \beta\right\} \,.
\end{equation}
Note that the probability $\Pb_\zs{\infty}(k\le \tau < k+m)$ is the probability of false alarm in the time interval $[k, k+m-1]$ of the length $m$, which we refer to as the {\em local probability of false alarm}
(LPFA), and the probability $\Pb_\zs{\infty}(\tau < k+m | \tau\ge k)= \Pb_\zs{\infty}(k\le \tau < k+m | \tau\ge k)$ is the corresponding {\em local conditional probability of false alarm} (LCPFA). 

We consider two risks: {\em positive part detection delay risk}
\begin{equation}\label{SUADD}
 \Rc_\nu(\tau)=   \EV_\zs{\nu}\left(\tau-\nu\right)^\zs{+}
\end{equation}
and {\em conditional detection delay risk}
\begin{equation}\label{SCADD}
\Rc^{*}_\nu(\tau)=  \EV_\zs{\nu}\left(\tau-\nu\,|\,\tau> \nu \right)
\end{equation}
(compare with \eqref{SADDdef}) and the following   problems:
the pointwise minimization, i.e., for any $\nu\ge 0$
\begin{equation}\label{sec:Prbf.5-0}
\inf_\zs{\tau\in \Hc(\beta, k^{*},m^{*})}\,
\Rc_\nu(\tau) \quad\mbox{and}\quad \inf_\zs{\tau\in \Hc^{*}(\beta, k^{*},m^{*})}\,\Rc^{*}_\nu(\tau)\,;
\end{equation}
and the minimax optimization
\begin{equation}\label{sec:Prbf.5}
\inf_\zs{\tau\in \Hc(\beta, k^{*},m^{*})}\,\sup_\zs{0 \le \nu < \infty}\,\Rc_\nu(\tau)
\quad\mbox{and}\quad \inf_\zs{\tau\in \Hc^{*}(\beta, k^{*},m^{*})}\,\max_\zs{0\le \nu \le k^{*}}\,\Rc^{*}_\nu(\tau)\,.
\end{equation}
\noindent The parameters $k^{*}$ and $m^{*}$  will be specified later.

In addition, we consider a Bayesian-type problem of minimizing the risks \eqref{SUADD} and \eqref{SCADD} in a class of procedures with the given weighted probability of false alarm. 
This problem is formulated and solved in the next section.

It would be more natural to address the classes of detection procedures with the given LPFA and LCPFA defined as
\[
\LPFA(\tau) = \sup_\zs{1\le k<\infty}\, \Pb_\zs{\infty}(k\le \tau < k+m^{*}) \quad \text{and} \quad \LCPFA(\tau) =\sup_\zs{1\le k<\infty}\, \Pb_\zs{\infty}(k\le \tau < k+m^{*} | \tau \ge k)
\]
and the maximal risks 
\[
 \sup_\zs{0 \le \nu<\infty}\, \Rc_\nu(\tau)   \quad \text{and} \quad \sup_\zs{0 \le \nu < \infty}\, \Rc^{*}_\nu(\tau) \, ,
\]
i.e., the optimality criteria
\begin{equation}\label{sec:Prbf.Lai}
\inf_\zs{\{\tau: \LPFA(\tau)\le \beta\}}\,
\sup_\zs{0 \le \nu<\infty}\, \Rc_\nu(\tau)  \quad \text{and} \quad  \inf_\zs{\{\tau: \LCPFA(\tau) \le \beta\}} \,
 \sup_\zs{0 \le \nu < \infty}\, \Rc^{*}_\nu(\tau) \, ,
\end{equation}
as in \citet{LaiIEEE98}, \citet{TartakovskyIEEECDC05} and \citet{TNBbook14}. However, in this case, one requires much stronger essential supremum conditions 
on the tail probabilities of the log-likelihood ratio, which do not hold in certain interesting examples (see Remark~\ref{Rem:Lai} below for details). For this reason, we modified these more natural optimality criteria.
In the following, we suppose that $k^*$, $m^*$ and $m^*-k^*$ go to infinity as $\beta\to0$, so that for practical purposes the optimality criteria \eqref{sec:Prbf.5}, 
considered in the present paper,  are not too much different
from the criteria \eqref{sec:Prbf.Lai}. At the same time, this allows us to substantially relax the sufficient conditions for asymptotic optimality of the detection procedures.

We need the following definition.

\begin{definition}\label{rcomplete}
For $k=0,1,\dots$ and $r>0$, we say that the normalized LLR  process $n^{-1} Z_{n+k}^k$ converges {\em $r-$completely} to a constant $I$ under the probability measure $\Pb_k$ as $n \to \infty$ if
\begin{equation}\label{rcompletedef}
\sum_{n=1}^\infty n^{r-1} \Pb_k\set{\abs{n^{-1}Z_{n+k}^k - I} > \varepsilon} < \infty \quad \text{for all}~~ \varepsilon > 0.
\end{equation}
If 
\begin{equation}\label{rcompletedefuniform}
\sum_{n=1}^\infty n^{r-1} \sup_{k \ge 0} \Pb_k\set{\abs{n^{-1}Z_{n+k}^k - I} > \varepsilon} < \infty \quad \text{for all}~~ \varepsilon > 0
\end{equation}
we say that $n^{-1} Z_{n+k}^k$ converges to a constant $I$ {\em uniformly $r-$completely} as $n\to\infty$.
\end{definition}

The $r-$complete convergence is an extension (for $ r \neq 1$) of the complete convergence introduced by \citet{HsuRobbins47}. It was introduced and extensively used for various hypothesis testing and change 
detection problems  by \citet{TNBbook14}.

In the following, we mostly deal with the case that $r=1$. In this case, we refer to \eqref{rcompletedef} as $\Pb_k$-complete convergence and to \eqref{rcompletedefuniform} as 
{\em uniform complete convergence}.

Note that, for any $r\ge 1$, $r-$complete convergence implies almost sure convergence of  $n^{-1} Z_{n+k}^k$ to $I$ under $\Pb_k$. Hence it can be interpreted as a rate of convergence 
in the strong law of large numbers. See \citet[Ch 2]{TNBbook14} for further details.


\section{Asymptotic optimality in the Bayesian-type class}\label{sec:Bay} 

We begin with considering a Bayesian-type class of change detection procedures that upper-bounds a weighted probability of false alarm 
$\PFA(\tau) = \sum_{k=0}^\infty \Pb_k(\tau \le k) \Pb(\nu =k)$, assuming that the change point $\nu$ is a random variable independent of the observations with prior distribution 
$\Pb(\nu=k)$, $k=0,1,2,\dots$. However, instead of considering a Bayes risk (weighted average delay to detection)
\begin{equation} \label{WADD}
\EV(\tau -\nu | \tau > \nu)= \frac{\sum_{k=0}^\infty \Pb(\nu=k) \Rc_k(\tau)}{1-\PFA(\tau)}\,,
\end{equation}
as it was done by \citet{TV2004}, we are interested in risks \eqref{SUADD} and \eqref{SCADD}, i.e., in the optimization problems 
\begin{equation}\label{sec:PrbfBayes1}
\inf_\zs{\{\tau: \PFA(\tau) \le \alpha\}}\,\Rc_k(\tau)  \quad \text{and} \quad   \inf_\zs{\{\tau: \PFA(\tau) \le \alpha\}}\, \Rc^*_k(\tau) \quad \text{for all}~ k \ge 0 \, ,
\end{equation}
and
\begin{equation}\label{sec:PrbfBayes}
\inf_\zs{\{\tau: \PFA(\tau) \le \alpha\}}\,\sup_\zs{k\ge 0}\Rc_k(\tau)  \quad \text{and} \quad   \inf_\zs{\{\tau: \PFA(\tau) \le \alpha\}}\,\max_\zs{0\le k\le k^{*}}  \Rc^*_k(\tau),
\end{equation}
where $0<\alpha <1$ is a prespecified (small) number.

In what follows, for simplicity of the presentation, assume that the prior probability distribution $\Pb(\nu=k)$ of the change point $\nu$ is geometric with the parameter $0<\varrho <1$, i.e., 
\begin{equation}\label{sec:Prbf.1}
\Pb(\nu=k)=\pi_\zs{k}(\varrho)=\varrho\,\left(1-\varrho\right)^{k}\,, \quad k =0,1,2,\dots
\end{equation}
Using this distribution we introduce the probability measure on the Borel $\sigma$-algebra in $\bbr^{\infty}\times \bbn$ as
$$
\Qb_\varrho(A\times J)=\sum_\zs{k\in J}\,\pi_\zs{k}(\varrho)\, \Pb_\zs{k}\left(A\right)\,, \quad A\in \Bc(\bbr^{\infty}) \, , ~~ J\subseteq\bbn\,.
$$ 
Now, for some fixed $0<\varrho,\alpha<1$, we define the following {\em Bayesian class} of change-point detection procedures with the weighted PFA $\Qb_\varrho(\tau \le \nu)=\PFA(\tau)$ 
not greater that the given number $\alpha$:
\begin{equation}\label{sec:Prbf.2}
\Delta(\alpha,\varrho)= \left\{\tau\in\Mcc:\, \Qb_\varrho\left(\tau \le \nu\right)\le \alpha\right\}
= \left\{\tau\in\Mcc:\, \sum_\zs{k\ge 1}\,\pi_\zs{k}(\varrho)\, \Pb_\zs{\infty}\left(\tau \le k\right)\le \alpha \right\}\,,
\end{equation}
where we took into account that $\Pb_k(\tau \le k) = \Pb_\infty (\tau\le k)$.

It follows from \cite{TV2004} that in the Bayesian setting, when one wants to minimize the weighted average delay to detection \eqref{WADD}, 
the asymptotically (as $\alpha\to 0$) optimal detection procedure in the class \eqref{sec:Prbf.2} is the Shiryaev detection procedure that raises an alarm at the first time such that the
posterior probability $g_\zs{n}(\varrho)=\Qb\left(\nu < n\mid \Fc_\zs{n}\right)$ exceeds threshold $1-\alpha$, i.e.,
$$
\wt{\tau}_\zs{b}(\alpha,\varrho)=\inf\{n\ge 1: g_\zs{n}(\varrho) \ge 1-\alpha\}\,.
$$
Note that it is easy to show  \cite{TV2004} that $\wt{\tau}_\zs{b}\in \Delta(\alpha,\varrho)$ for any $0<\alpha,\varrho<1$.  

Using the LLR process $(Z^{k}_\zs{n})_\zs{0\le k\le n-1}$ defined in \eqref{sec:Prbf.7}, the posterior probability $g_\zs{n}(\varrho)$ can be represented as
$g_\zs{n}(\varrho)=\Lambda_\zs{n}(\varrho)/[1/\varrho+\Lambda_\zs{n}(\varrho)]$, where 
\begin{equation}\label{sec:Prbf.9}
\Lambda_\zs{n}(\varrho)= \sum_{k=0}^{n-1}\,\left(1-\varrho\right)^{-(n-k)}\,e^{Z^{k}_{n}}\,.
\end{equation}
Therefore, the Shiryaev procedure can be also written as
\begin{equation}\label{sec:Prbf.8}
\wt{\tau}(\alpha, \varrho) =\inf \left\{n\ge 1: \Lambda_\zs{n}(\varrho)\ge (1-\alpha)/(\varrho \alpha) \right\}\,.
\end{equation}

Note first that, as $\rho\to0$, the statistic $\Lambda_\zs{n}(\varrho)$ converges to the SR statistic,
\[
\Lambda_\zs{n}(\varrho) \xra[\rho\to0]{} \sum_{k=1}^n e^{Z_n^{k-1}} \, .
\]
Thus, if we are interested in small values of $\varrho$, as it is the case in the following, then the behavior of the Shiryaev procedure  \eqref{sec:Prbf.8} is similar to that of the SR procedure $T(h)$ so long as
we can define the threshold $h=h_\alpha$ is such a way that $\PFA(T(h_\alpha)) \le \alpha$. Hence, instead of considering the procedure 
$\wt{\tau}(\alpha, \varrho)$,  which is shown in \cite{TV2004} to be asymptotically optimal in the Bayesian context with respect to the weighted average delay to detection \eqref{WADD}, 
we will focus on the SR procedure.

\subsection{Asymptotic lower bounds} \label{ssec:ALBBayes}

We do not assume any particular model or even class of models for the observations, and as a result, there is no
``structure" of the LLR process. We therefore have to impose some conditions on the behavior of the LLR process
at least for large $n$. It is natural to assume that there exists a positive finite
number $I$ such that $Z_\zs{n}^k/(n-k)$ converges almost surely to $I$ under $\Pb_k$, i.e., \vspace{3mm}

\noindent $(\A_\zs{1}$)
{\em  Assume that there exists a number $I>0$ such that for any $k\ge 0$}
\begin{equation}\label{sec:MaRe.1}
\frac{1}{n}Z_\zs{k+n}^{k} \xra[n\to\infty]{\Pb_\zs{k}-\text{a.s.}} I \,.
\end{equation}
\noindent This is always true for iid data models with 
\[
I=I(f_\zs{1},f_\zs{0})=\EV_\zs{0} Z_\zs{0}^1 = \int \log \brcs{\frac{f_1(x)}{f_0(x)}} f_1(x) \d\mu(x)
\]
being the Kullback--Leibler information number. It turns out that the a.s.\
convergence condition \eqref{sec:MaRe.1} is sufficient for obtaining lower bounds for all
positive moments of the detection delay.

The following theorem establishes  asymptotic lower bounds for the optimization problems \eqref{sec:PrbfBayes1} and \eqref{sec:PrbfBayes}.  
We write $\Delta(\alpha)$ for the class $\Delta(\alpha,\varrho_\alpha)$ 
when the parameter $\varrho=\varrho_\alpha$ depends on $\alpha$.

\begin{theorem} \label{Th.sec:Bay.1} 
 Assume that the almost sure convergence condition $(\A_\zs{1})$  holds and in \eqref{sec:Prbf.1} the parameter of the geometric prior distribution
$\varrho=\varrho_\zs{\alpha}\to 0$ as $\alpha\to 0$. Then, for any $\nu \ge 0$, 
\begin{equation} \label{sec:Bay.1}
\liminf_\zs{\alpha\to 0} \frac{1}{|\log\alpha|}\, \inf_\zs{\tau\in \Delta(\alpha)} \, \sup_{\nu\ge 0} \Rc_\nu(\tau) \ge  
\liminf_\zs{\alpha\to 0}  \frac{1}{|\log\alpha|}\, \inf_\zs{\tau\in \Delta(\alpha)} \, \Rc_\nu(\tau) \ge \frac{1}{I} 
\end{equation}
and
\begin{equation} \label{sec:Bay.2}
\liminf_\zs{\alpha\to 0} \frac{1}{|\log\alpha|}\, \inf_\zs{\tau\in \Delta(\alpha)} \, \sup_{\nu\ge 0} \Rc_\nu^{*}(\tau) \ge 
\liminf_\zs{\alpha\to 0} \frac{1}{|\log\alpha|}\, \inf_\zs{\tau\in \Delta(\alpha)}\, \Rc^{*}_\nu(\tau) \,\ge \frac{1}{I}\,.
\end{equation}
\end{theorem}

\proof
First, note that it is not difficult to prove that condition $(\A_\zs{1}$) implies  that
for any $\vae>0$ and $k\ge 0$
\begin{equation} \label{sec:Bay.2-0}
U_\zs{k,M}=\Pb_\zs{k}\set{\frac{1}{M}\max_{1\le n \le M}Z_\zs{k+n}^{k}\ge (1+\vae) I}\xra[M\to\infty]{}0\,.
\end{equation}
Next, it is clear that, for any $k\ge 0$, $\EV_\zs{k}\left(\tau-k\,|\,\tau> k\right)\ge \EV_\zs{k}\left(\tau-k\right)^+$,
i.e.,  the assertion \eqref{sec:Bay.2} follows from the assertion \eqref{sec:Bay.1}, and hence it suffices to prove
only inequality \eqref{sec:Bay.1}.

Define $\gamma^{(k)}_\zs{\varepsilon,\alpha}(\tau)=\Pb_\zs{k}(k\le \tau\le k+M_\zs{\varepsilon,\alpha})$ and $M_\zs{\varepsilon,\alpha}= (1-\varepsilon)|\log\alpha|/(I+d)$,
where $d=-\log(1-\varrho)$. Let us show that for any $k \ge 0$ and $0<\varepsilon<1$
\begin{equation} \label{sec:Bay.3}
\lim_\zs{\alpha\to 0}\sup_\zs{\tau\in \Delta(\alpha,\varrho)}\,\,\gamma^{(k)}_\zs{\varepsilon,\alpha}(\tau)\,=0\,.
\end{equation}
 Indeed, using the change of measure trick, similarly to \cite{TV2004} we  obtain
\begin{equation} \label{sec:Bay.4}
\begin{aligned} 
 \gamma^{(k)}_\zs{\varepsilon,\alpha}(\tau) \le 
e^{(1+\varepsilon)I  M_\zs{\varepsilon,\alpha}}\,
\Pb_\zs{\infty}\left(k\le \tau\le k+M_\zs{\varepsilon,\alpha}\right)
 +\Pb_\zs{k}\left(\max_\zs{1\le n\le M_\zs{\varepsilon,\alpha}}Z_\zs{k+n}^{k}\ge (1+\varepsilon)I  M_\zs{\varepsilon,\alpha} \right)\,.
\end{aligned}
\end{equation}

The definition of the class $\Delta(\alpha,\varrho)$ in \eqref{sec:Prbf.2} implies that for any $0<\alpha,\varrho<1$ 
and for any $l\ge 1$
$$
\alpha\ge\varrho \sum_\zs{k\ge l}\,(1-\varrho)^{k}\,\Pb_\zs{\infty}\left(\tau\le k\right)\,
\ge 
\varrho \Pb_\zs{\infty}\left(\tau\le l\right)
 \sum_\zs{k\ge l}\,(1-\varrho)^{k}\,=\,
 \Pb_\zs{\infty}\left(\tau\le l\right)\,(1-\varrho)^{l}\,,
$$
i.e.
\begin{equation}\label{sec:Bay.5}
\sup_\zs{\tau\in\Delta(\alpha,\varrho)}\,
\Pb_\zs{\infty}\left(\tau\le l\right)\le \alpha\,(1-\varrho)^{-l}\,.
\end{equation}

Therefore, the first term in the right side of the inequality \eqref{sec:Bay.4}
may be estimated as
$$
e^{(1+\varepsilon)I  M_\zs{\varepsilon,\alpha}-|\log\alpha|
+dk+dM_\zs{\varepsilon,\alpha}}
\le e^{-\varepsilon^{2} |\log\alpha|+dk}
$$
and it goes to zero for any fixed $0\le k<\infty$. By \eqref{sec:Bay.2-0}, the second term in \eqref{sec:Bay.4}, 
\[
\Pb_\zs{k}\left(\max_\zs{1\le n \le M_\zs{\varepsilon,\alpha}}Z_\zs{k+n}^{k}\ge (1+\varepsilon)I  M_\zs{\varepsilon,\alpha} \right) = U_{k,M_\alpha} \, ,
\]
approaches zero as $\alpha\to0$ for all $k \ge 0$, and we obtain \eqref{sec:Bay.3}.  By the Chebyshev inequality, 
$$
\EV_\zs{k}\left(\tau-k\right)^+\ge M_\zs{\varepsilon,\alpha}\,\,\Pb_\zs{k}\left( \tau > k+M_\zs{\varepsilon,\alpha}\right)=
M_\zs{\varepsilon,\alpha}\left(\Pb_\zs{k}\left( \tau \ge k\right) -\gamma^{(k)}_\zs{\varepsilon,\alpha}(\tau)\right)\,.
$$
From \eqref{sec:Bay.5} we obtain immediately that for any fixed $k\ge 0$
$$
\Pb_\zs{k}\left( \tau \ge k\right)=\Pb_\zs{\infty}\left( \tau \ge k\right)\ge 1-\alpha\,(1-\varrho)^{-k+1}
\to 1 \quad\mbox{as}\quad \alpha\to 0\,.
$$
Therefore,  for any $\nu\ge 0$ and for any small $\varepsilon$
\[
\lim_{\alpha\to0} \frac{\inf_\zs{\tau\in \Delta(\alpha)}\, \Rc_\nu(\tau)}{|\log\alpha|/(I+d_\alpha)} \ge 1-\varepsilon.
\]
Taking into account that $d_\alpha=-\log(1-\varrho_\alpha)\to0$ as $\alpha\to0$ and letting $\varepsilon\to 0$,  we obtain the lower bounds \eqref{sec:Bay.1}.  Hence Theorem~\ref{Th.sec:Bay.1}. 
\endproof

Observe that the lower bounds \eqref{sec:Bay.1} and \eqref{sec:Bay.2} can be generalized for all positive moments of the detection delay
$\EV_\zs{\nu}\left[(\tau-\nu)^r\,|\,\tau> \nu \right]$ and  $\EV_\zs{\nu}\left[(\tau-\nu)^+ \right]^r$, $r>1$. Indeed, using Jensen's inequality
$\EV_\zs{\nu}\left[(\tau-\nu)^+ \right]^r \ge \left[\EV_\zs{\nu}(\tau-\nu)^+ \right]^r$, we immediately obtain that under the conditions of Theorem~\ref{Th.sec:Bay.1}, for any $\nu \ge 0$,
\begin{equation} \label{sec:BayrLB}
\liminf_\zs{\alpha\to 0} \frac{1}{|\log \alpha|^r} \inf_\zs{\tau\in \Delta(\alpha)}\, \EV_\zs{\nu}\left[(\tau-\nu)^+ \right]^r\, \ge \frac{1}{I^r}
\end{equation}
and analogously
\begin{equation} \label{sec:CBayrLB}
\liminf_\zs{\alpha\to 0} \frac{1}{|\log \alpha|^r} \inf_\zs{\tau\in \Delta(\alpha)}\, \EV_\zs{\nu}\left[(\tau-\nu)^r \,|\,\tau> \nu \right]\, \ge \frac{1}{I^r} \, .
\end{equation}
Since higher moments of the detection delay may also be of interest, the asymptotic lower bounds \eqref{sec:BayrLB} and \eqref{sec:CBayrLB} can be useful for establishing asymptotic optimality properties 
of the SR procedure with respect to the risks $\EV_\zs{\nu}\left[(\tau-\nu)^r\,|\,\tau> \nu \right]$ and  $\EV_\zs{\nu}\left[(\tau-\nu)^+ \right]^r$ for $r>1$ uniformly for all $\nu\ge 0$, as well as with respect to the 
maximal risks.

\subsection{Asymptotic optimality of the Shiryaev--Roberts procedure}\label{ssec:AOBayes}

In order to study asymptotics for the average detection delay of the SR procedure and for establishing its asymptotic optimality, we impose the following constraint on the rate of convergence for
\begin{equation}\label{sec:MaRe.2}
\wt{Z}_\zs{k,n}= \frac{1}{n}Z_\zs{k+n}^{k} - I \,.
\end{equation}

\noindent $(\A_\zs{2}$) {\em 
Assume that $\wt{Z}_\zs{k,n}$ converges uniformly completely to $0$ as $n \to \infty$, i.e.,
for any  $\vae>0$} 
\begin{equation} \label{sec:MaRe.3}
\Upsilon^{*}(\varepsilon)= \sum_{n=1}^\infty\,\sup_\zs{k\ge 0}\, \Pb_\zs{k}\Bigl\{\Bigl|\wt{Z}_\zs{k,n}\Bigr|>\vae \Bigr\}<\infty \,.
\end{equation}

Write $R_n = \sum_{k=1}^n e^{Z_n^{k-1}}$ for the SR statistic and denote as $T(h)=T(\alpha,\varrho)$ the SR procedure when the threshold $h=h(\alpha,\varrho)$ is 
selected as $h(\alpha,\varrho)=(1-\alpha)/\varrho\alpha$, i.e.,
\begin{equation}\label{SRrhoalpha}
T(\alpha,\varrho) =  \inf\set{n \ge 1:  R_n \ge \frac{1-\alpha}{\varrho\alpha}}.
\end{equation}

\begin{lemma}\label{Lem:PFASR} 
The SR procedure $T(\alpha,\varrho)$ given by \eqref{SRrhoalpha} belongs to the class $\Delta(\alpha,\varrho)$ for any $0<\alpha, \varrho<1$.
\end{lemma}

\proof
Note that the stopping time \eqref{sec:Prbf.8} can be written as
$$
\wt{\tau}_\zs{b}(\alpha,\varrho) =\inf \left\{n\ge 1: \,\sum_{k=1}^{n}\,\left(1-\varrho\right)^{-[n-(k-1)]}\,e^{Z^{k-1}_\zs{n}}\,\ge h(\alpha,\varrho)\right\}
$$
(see \eqref{sec:Prbf.9}). Obviously,  $\wt{\tau}_\zs{b}(\alpha,\varrho) \le T(\alpha,\varrho)$ almost surely for any $0<\alpha, \varrho <1$.
Since by (2.11) in \cite{TV2004} $\PFA(\wt{\tau}_\zs{b}(\alpha,\varrho)) \le \alpha$, it follows that,  for any $0<\alpha,\varrho<1$, $T(\alpha,\varrho) \in\Delta(\alpha,\varrho)$
and the proof is complete. 
\endproof

In what follows, we assume that the parameter $\varrho$ is a function of $\alpha$, i.e. $\varrho=\varrho_\zs{\alpha}$, such that
\begin{equation}\label{sec:Bay.6}
\lim_\zs{\alpha\to 0}\,\varrho_\zs{\alpha}\,=0 \quad\mbox{and}\quad \lim_\zs{\alpha\to 0}\,\frac{|\log\varrho_\zs{\alpha}|}{|\log\alpha|}\,=0\,.
\end{equation}

\noindent Moreover, let  $k^{*}$ be a function of $\alpha$, i.e. $k^{*}=k^{*}_\alpha$, such that
\begin{equation}\label{sec:Bay.7}
\lim_\zs{\alpha\to 0} k^{*}_\alpha=\infty \quad\mbox{and}\quad \lim_\zs{\alpha\to 0} \left(|\log\alpha|+  k^{*}_\alpha\log(1-\varrho_\alpha)\right)=+\infty\,.
\end{equation}

Denote as $T_\alpha=T(h_\alpha)$ the SR procedure defined in \eqref{SRrhoalpha} when the threshold $h(\alpha,\varrho_\alpha)=h_\alpha$ is selected as 
$h_\alpha= (1-\alpha)/(\varrho_\alpha \alpha)$.
Note that if conditions \eqref{sec:Bay.6} hold, then $h_\alpha\to\infty$ as $\alpha\to0$. 
Clearly, we need the threshold to become large for small $\alpha$; otherwise the problem is degenerate. By Lemma~\ref{Lem:PFASR}, this choice of the threshold guarantees that 
$T_\alpha \in \Delta(\alpha, \varrho_\alpha) = \Delta(\alpha)$ for every $0<\alpha<1$.

The following theorem identifies the asymptotic upper bounds for the risks of the SR procedure.

\begin{theorem}\label{Th.sec:Bay.2} 
{\rm \bf(i)} Assume that the uniform complete convergence condition $(\A_\zs{2})$  holds for some  $0<I<\infty$, and the parameter $0<\varrho=\varrho_\alpha<1$  in the SR procedure \eqref{SRrhoalpha}
satisfies conditions \eqref{sec:Bay.6}. Then
\begin{equation}\label{sec:Bay.8}
\limsup_\zs{\alpha\to 0}\,\frac{1}{|\log\alpha|}\,\sup_\zs{\nu\ge 0}\,\Rc_\nu(T_\alpha)\,\le \frac{1}{I}\,.
\end{equation}

\noindent {\rm \bf(ii)}  Assume that in addition to conditions $(\A_\zs{2})$ and \eqref{sec:Bay.6}, conditions \eqref{sec:Bay.7} hold for $k^*=k^*_\alpha$. Then
\begin{equation} \label{sec:Bay.9}
\limsup_\zs{\alpha\to 0}\,\frac{1}{|\log\alpha|} \,\max_\zs{0 \le \nu \le k^{*}_\alpha}\,\Rc^*_\nu(T_\alpha) \,\le \frac{1}{I}\,.
\end{equation}
\end{theorem}

\proof Proof of (i). Taking into account that $R_{\nu+n}=\sum_{i=1}^{\nu+n} e^{Z^{i-1}_{\nu+n}} \ge e^{Z^{\nu}_{\nu+n}}$, we obtain
\[
\begin{aligned}
\Pb_\zs{\nu}\brc{T_\alpha>\nu+n}  &= \Pb_\nu\brc{R_\ell < h_\alpha, \ell =1,\dots,\nu+n} \le \Pb_\nu\brc{R_{\nu+n} < h_\alpha} 
\\
& \le  \Pb_\nu\brc{Z_{\nu+n}^\nu < \log h_\alpha} = \Pb_\nu\brc{\wt{Z}_\zs{\nu,n} < -I+\log h_\alpha/n} \,.
\end{aligned}
\]
Evidently, for any $0<\varepsilon<I$  and any $n \ge n_\zs{0}=\lfloor \log h_\alpha/(I-\varepsilon) \rfloor+1$, the last probability can be bounded as
\begin{equation}\label{upperprobnu}
\Pb_\nu\brc{\wt{Z}_\zs{\nu,n} < -I+\log h_\alpha/n} \le \Pb_\nu\brc{\wt{Z}_\zs{\nu,n} < -\varepsilon} \le 
\Pb_\nu\brc{\vert \wt{Z}_\zs{\nu,n}\vert > \varepsilon} \,,
\end{equation}
and hence, for any $0<\varepsilon<I$,
\begin{equation}\label{Upper}
\begin{aligned}
 \EV_\zs{\nu}\left(T_\alpha-\nu\right)^+
 &=\sum_\zs{n\ge 0}\Pb_\zs{\nu}\left(T_\alpha>\nu+n\right) \le  
 \sum_\zs{n\ge 0}\,\Pb_\nu\brc{\wt{Z}_\zs{\nu,n} < -I+\log h_\alpha/n}
 \\
&= \sum_\zs{0\le n< n_\zs{0}}\,\Pb_\nu\brc{\wt{Z}_\zs{\nu,n} < -I+\log h_\alpha/n}
 +\sum_\zs{n\ge n_\zs{0}}\,\Pb_\nu\brc{\wt{Z}_\zs{\nu,n} < -I+\log h_\alpha/n}\\
 &\le  n_\zs{0}+ \sum_\zs{n\ge n_\zs{0}}\,\Pb_\nu\brc{\vert \wt{Z}_\zs{\nu,n}\vert >\varepsilon}\le \frac{|\log h_\alpha|}{I-\varepsilon}+1+\Upsilon^{*}(\varepsilon)
 \\
 & = \frac{|\log\alpha| +|\log\varrho _\alpha| +\log(1-\alpha)}{I-\varepsilon}+1+\Upsilon^{*}(\varepsilon)\,.
\end{aligned}
\end{equation}
Hereafter $\lfloor x \rfloor$ denotes the integer number less than or equal to $x$. Since the right-hand side does not depend on $\nu$, 
we have (for any $\nu \ge 0$ and $0<\varepsilon < I$)
\[
\Rc_\nu(T_\alpha) \le \sup_{\nu \ge 0} \Rc_\nu(T_\alpha) \le \frac{|\log\alpha|}{I-\varepsilon} \brc{1+ \frac{|\log\varrho_\alpha| +\log(1-\alpha)}{|\log \alpha|}}+1+\Upsilon^{*}(\varepsilon) \, .
\]
By condition $(\A_\zs{2}$),  $\Upsilon^{*}(\varepsilon)<\infty$, so using condition  \eqref{sec:Bay.6} and the fact that $\varepsilon$ is arbitrary 
yields the upper bound \eqref{sec:Bay.8} and the assertion (i) follows.

Proof of (ii). In view of inequality  \eqref{sec:Bay.5}, for any $0 \le \nu \le k^*_\alpha$,
\[
\Pb_\zs{\nu}\left(T_\alpha > \nu\right)= \Pb_\zs{\infty}\left(T_\alpha > \nu\right)\ge  \Pb_\zs{\infty}\left(T_\alpha > k_\alpha^{*} \right) \ge 1-  \alpha(1-\varrho_\alpha)^{-k^*_\alpha}.
\]
Evidently, under conditions \eqref{sec:Bay.6}  and \eqref{sec:Bay.7} the right-hand side approaches $1$ as $\alpha \to 0$,
which implies that $\Pb_{\nu}\left(T_\alpha > \nu\right) \to 1$ as $\alpha\to0$  for all $0 \le \nu\le k^*_\alpha$. Since 
$\EV_\nu(T_\alpha -\nu | T_\alpha >\nu) =\EV_\nu(T_\alpha -\nu)^+/\Pb_\zs{\nu}\left(T_\alpha > \nu\right)$,
inequality \eqref{sec:Bay.8} implies \eqref{sec:Bay.9} for $k^*=k^*_\alpha$ satisfying conditions \eqref{sec:Bay.7} and the proof is complete.\endproof

Finally, combining Theorem~\ref{Th.sec:Bay.1} and Theorem~\ref{Th.sec:Bay.2}, we conclude that the SR procedure is first-order asymptotically  uniformly pointwise optimal and minimax
in the class $\Delta(\alpha)$, which is formalized in the next theorem. 

\begin{theorem}\label{Th.sec:Bay.3} 
{\rm \bf(i)} Assume that the uniform complete convergence condition $(\A_\zs{2})$  holds for some  $0<I<\infty$, and the parameter $0<\varrho=\varrho_\alpha<1$  in the SR procedure \eqref{SRrhoalpha}
satisfies conditions \eqref{sec:Bay.6}. Then
\begin{equation}\label{sec:Bay.AO1}
\lim_\zs{\alpha\to 0}\,\frac{\inf_{\tau\in\Delta(\alpha)} \Rc_\nu(\tau)}{\Rc_\nu(T_\alpha)} = 1 \quad \text{for all}~ \nu \ge 0 
\end{equation}
and 
\begin{equation}\label{sec:Bay.AO}
\lim_\zs{\alpha\to 0}\,\frac{\inf_{\tau\in\Delta(\alpha)}\sup_\zs{\nu\ge 0}\Rc_\nu(\tau)}{\sup_\zs{\nu\ge 0}\Rc_\nu(T_\alpha)} = 1 \, .
\end{equation}
Moreover, as $\alpha\to0$,
\[
\inf_{\tau\in\Delta(\alpha)}\Rc_\nu(\tau) ~ \sim ~ \Rc_\nu(T_\alpha) ~ \sim ~ \frac{|\log \alpha|}{I} \quad \text{for all}~ \nu \ge 0
\]
and
\[
\inf_{\tau\in\Delta(\alpha)}\sup_{\nu\ge 0} \Rc_\nu(\tau) ~ \sim ~ \sup_{\nu\ge 0} \Rc_\nu(T_\alpha) ~ \sim ~ \frac{|\log \alpha|}{I} \, .
\]

\noindent{\rm \bf(ii)}  Assume that in addition to conditions $(\A_\zs{2})$ and \eqref{sec:Bay.6} conditions \eqref{sec:Bay.7} hold for $k^*=k^*_\alpha$. Then 
\begin{equation}\label{sec:Bay.AOcond1}
\lim_{\alpha\to 0}\,\frac{\inf_{\tau\in\Delta(\alpha)} \Rc_{\nu}^{*}(\tau)}{\Rc_{\nu}^{*}(T_\alpha)}= 1 \quad \text{for all fixed}~ \nu \ge 0
\end{equation}
and
\begin{equation}\label{sec:Bay.AOcond}
\lim_{\alpha\to 0}\,\frac{\inf_{\tau\in\Delta(\alpha)}\max_{0\le \nu \le k_\alpha^*} \Rc^*_\nu(\tau)}{\max_{0\le \nu \le k_\alpha^*}\Rc^*_\nu(T_\alpha)}= 1 \, .
\end{equation}
Moreover, as $\alpha\to0$, 
\[
\inf_{\tau\in\Delta(\alpha)}\Rc_\nu^*(\tau) ~ \sim ~ \Rc_\nu^*(T_\alpha) ~ \sim ~ \frac{|\log \alpha|}{I} \quad \text{for all fixed}~ \nu \ge 0
\]
and
\[
\inf_{\tau\in\Delta(\alpha)} \max_{0\le \nu \le k_\alpha^*} \Rc^*_\nu(\tau) ~ \sim ~  \max_{0\le \nu \le k_\alpha^*} \Rc^*_\nu(T_\alpha) ~ \sim ~ \frac{|\log \alpha|}{I}.
\]
\end{theorem}

\proof
All assertions follow from Theorem~\ref{Th.sec:Bay.1} and Theorem~\ref{Th.sec:Bay.2} in an obvious manner.
\endproof

The above asymptotic optimality results can be generalized for higher moments of the detection delay if the uniform complete convergence condition $(\A_\zs{2}$) is strengthened into the uniform 
$r-$complete convergence condition for some $r>1$. In particular, the following result holds true.

\begin{theorem}\label{Th.sec:Bay.rcomp} 
Let conditions \eqref{sec:Bay.6} and \eqref{sec:Bay.7} hold and, for some $r>1$ and all  $\vae>0$,
\begin{equation} \label{sec:UCCCr}
\sum_{n=1}^\infty\, n^{r-1} \, \sup_\zs{k\ge 0}\, \Pb_\zs{k}\Bigl\{\Bigl|\wt{Z}_\zs{k,n}\Bigr|>\vae \Bigr\}<\infty \,.
\end{equation}
Then the SR procedure $T_\alpha$ is first-order asymptotically uniformly pointwise optimal and minimax in the class $\Delta(\alpha,\varrho_\alpha)=\Delta(\alpha)$ with respect to the moments of the 
detection delay up to order $r$: for all $1\le \ell \le r$ as $\alpha\to0$
\begin{equation}\label{sec:Bay.SRAO}
\begin{aligned}
\EV_\zs{\nu}\left[(T_\zs{\alpha}-\nu)^\ell | T_\zs{\alpha}>\nu\right] & ~ \sim~  \inf_\zs{\tau\in \Delta(\alpha)}\, \EV_\zs{\nu}\left[(\tau-\nu)^\ell | \tau>\nu \right]
\\
& \sim~ \brc{\frac{|\log \alpha|}{I}}^\ell \quad \text{for all fixed}~ \nu \ge 0
\end{aligned}
\end{equation}
and
\begin{equation}\label{sec:Bay.SRAO1}
\sup_\zs{0 \le \nu \le k^{*}_\alpha}\, \EV_\zs{\nu}\left[(T_\zs{\alpha}-\nu)^\ell | T_\zs{\alpha}>\nu\right] ~ \sim~  \inf_\zs{\tau\in \Delta(\alpha)}\, 
 \sup_\zs{0 \le \nu \le k^{*}_\alpha}\,  \EV_\zs{\nu}\left[(\tau-\nu)^\ell | \tau>\nu \right]~ \sim~ \brc{\frac{|\log \alpha|}{I}}^\ell\,.
\end{equation}
\end{theorem}

\proof
We give only a sketch of the proof and omit certain details. By \eqref{sec:CBayrLB}, for all $r \ge 1$ and for all $\nu\ge 0$,
\begin{equation}\label{LBnew}
\liminf_{\alpha\to 0} \frac{1}{|\log \alpha|^r} \inf_{\tau\in \Delta(\alpha)}\, \EV_{\nu}\left[(\tau-\nu)^r \,|\,\tau> \nu \right]\, \ge \frac{1}{I^r} \, .
\end{equation}
(This bound holds since condition \eqref{sec:UCCCr} obviously implies the a.s.\ convergence ($\A_1$), which in turn implies \eqref{sec:Bay.2-0} for all $k\ge0$.)

Next, using the reasoning similar to that used in the proof of Proposition~\ref{Prop: A1} in the Appendix, which has lead to inequality \eqref{IneqExpTb}, we obtain
\begin{equation}\label{IneqExpTalpha}
\EV_{\nu}\brcs{(T_\alpha-\nu)^+}^r  \le  n_0^{r} + r 2^{r-1} \sum_{n=n_0}^{\infty}  n^{r-1}   \Pb_\nu (T_\alpha > \nu+ n) \, ,
\end{equation}
where $n_0=1+\lfloor \log h_\alpha/(I-\varepsilon) \rfloor$ and $h_\alpha=(1-\alpha)/\alpha\varrho_\alpha$. By \eqref{upperprobnu}, for any $0<\varepsilon<I$  and $n \ge n_0$,
$\Pb_\nu (T_\alpha > \nu+ n) \le \Pb_\nu(\vert \wt{Z}_{\nu,n}\vert > \varepsilon)$, and we obtain
\[
\begin{aligned}
\EV_{\nu}\brcs{(T_\alpha-\nu)^+}^r  & \le  n_0^{r} + r 2^{r-1} \sum_{n=n_0}^{\infty}  n^{r-1} \Pb_\nu\brc{\vert \wt{Z}_{\nu,n}\vert > \varepsilon} 
\\
&\le \brc{1+ \frac{\log h_\alpha}{I-\varepsilon}}^r +  r 2^{r-1} \sum_{n=1}^{\infty}  n^{r-1} \sup_{\nu\ge0} \Pb_\nu\brc{\vert \wt{Z}_{\nu,n}\vert > \varepsilon} \, ,
\end{aligned}
\]
where the last term is finite due to condition \eqref{sec:UCCCr} and $\log h_\alpha \sim |\log\alpha|$ as $\alpha\to0$ due to condition \eqref{sec:Bay.6}. Thus, for an arbitrary
$0<\varepsilon<I$,
\begin{equation}\label{UBnew}
\EV_{\nu}\brcs{(T_\alpha-\nu)^+}^r  \le \brc{\frac{|\log \alpha|}{I-\varepsilon}}^r(1+o(1)) \quad \text{as}~ \alpha\to0,
\end{equation}
and we established the asymptotic upper bound
\[
\limsup_{\alpha\to 0} \frac{1}{|\log \alpha|^r} \EV_{\nu}\left[(T_\alpha-\nu)^r \,|\,\tau> \nu \right]\, \le \frac{1}{I^r} \quad \text{for all}~\nu\ge 0.
\]
Applying this upper bound together with the lower bound \eqref{LBnew} proves asymptotic relations \eqref{sec:Bay.SRAO}. 

The upper bound
\[
\limsup_{\alpha\to 0} \frac{1}{|\log \alpha|^r} \, \sup_{0 \le \nu \le k^{*}_\alpha}\, \EV_{\nu}\left[(T_\alpha-\nu)^r \,|\,\tau> \nu \right]\, \le \frac{1}{I^r} 
\]
can be established similarly to \eqref{sec:Bay.9}, using \eqref{UBnew} and the fact that $\max_{0 \le \nu\le k^*_\alpha} \Pb_{\nu}(T_\alpha > \nu) \to 1$ as $\alpha\to0$ (see the proof of Theorem~\ref{Th.sec:Bay.2}(ii)) 
and that 
\[
\EV_\nu[(T_\alpha -\nu)^r | T_\alpha >\nu) =\EV_\nu[(T_\alpha -\nu)^+]^r/\Pb_{\nu}(T_\alpha > \nu) \, .
\]
This upper bound and the lower bound \eqref{LBnew} imply  \eqref{sec:Bay.SRAO1}.
\endproof

\begin{remark}
While for the sake of simplicity we consider the geometric prior distribution with the small parameter 
$\varrho_\alpha\to 0$ as $\alpha\to 0$, all the asymptotic results hold true for an arbitrary prior distribution 
$\pi_k^\alpha$ such that  the mean value of the change point $\EV \, \nu = \sum_{k=1}^\infty k \pi_k^\alpha$ approaches infinity as $\alpha\to 0$,
assuming that conditions \eqref{sec:Bay.6} and \eqref{sec:Bay.7} hold with $\varrho_\alpha$ replaced by $(\sum_{k=1}^\infty k \pi_k^\alpha)^{-1}$.
\end{remark}

\begin{remark}
Analogous asymptotic optimality results hold for the Shiryaev procedure $\wt{\tau}(\alpha)$ defined in  \eqref{sec:Prbf.8}. The proofs are essentially similar.
\end{remark}

\section{Asymptotic optimality in classes with given local probabilities of false alarm} \label{sec:MaRe}

We now proceed with tackling the pointwise and minimax problems \eqref{sec:Prbf.5-0} and \eqref{sec:Prbf.5} in the classes of procedures with given LPFA and LCPFA.  
The method of establishing asymptotic optimality of the SR procedure is again based on the lower-upper bounding technique. Specifically, we first obtain asymptotic lower bounds for the 
risk $\Rc_\nu(\tau)$ in the class $\Hc\left(\beta,k^{*},m^{*}\right)$ and for the risk $\Rc^{*}_\nu(\tau)$ in the class $\Hc^{*}\left(\beta,k^{*},m^{*}\right)$,  and then we show that these asymptotic lower 
bounds are attained for the SR procedure $T(h)$ with a certain threshold $h=h_\beta$. Note that the asymptotic optimality results of the previous section are essential, since asymptotic optimality 
in classes $\Hc\left(\beta,k^{*},m^{*}\right)$ and $\Hc^{*}\left(\beta,k^{*},m^{*}\right)$ is obtained by imbedding these classes in the class $\Delta(\alpha,\rho)$ with specially selected 
parameters $\rho$ and $\alpha$.

\subsection{Asymptotic lower bounds} \label{ssec:ALBminimax}

For any  $0<\beta<1$, $m^{*}\ge 1$ and $k^{*}>m^{*}$, define
\begin{equation}\label{sec:Absrsk.1-0}
\alpha_\zs{1}=
\alpha_\zs{1}(\beta,m^{*})
=\beta+(1-\varrho_\zs{1,\beta})^{m^{*}+1}
\end{equation}
and
\begin{equation}\label{sec:Absrsk.1}
\alpha_\zs{2}=\alpha_\zs{2}(\beta,k^{*})= \beta(1-\varrho_\zs{2,\beta})^{k^{*}}
\,,
\end{equation}
\noindent where 
\begin{equation}\label{sec:Absrsk.1-01}
\varrho_\zs{1,\beta}=\frac{1}{1+|\log \beta|}
\quad\mbox{and}\quad 
\varrho_\zs{2,\beta}=
\frac{\varrho_\zs{1,\beta}}{1+|\log\,|\log \beta||}\,.
\end{equation}

To find asymptotic lower bounds for the problems  \eqref{sec:Prbf.5-0}  and \eqref{sec:Prbf.5} in addition to condition $(\A_\zs{1}$) we impose the following condition related to the growth of the window 
size $m^*$ in the LPFA: \vspace{3mm}

\noindent $(\H_\zs{1})$ {\em  The size of the window $m^{*}$ in \eqref{sec:Absrsk.1-0}  is a function of $\beta$, i.e. $m^{*}=m^{*}_\zs{\beta}$,  such that
\begin{equation}\label{sec:Absrsk.4}
\lim_\zs{\beta\to 0}\, \frac{|\log \alpha_\zs{1,\beta}|}{|\log\beta|}=1\,,
\end{equation}
where 
$\alpha_\zs{1,\beta}=\alpha_\zs{1}(\beta,m^{*}_\zs{\beta})$.}

For example, we can take $m^{*}_\zs{\beta}=1+\lfloor (1+|\log\beta|)^{2} \rfloor$. 

The following theorem establishes asymptotic lower bounds. 

\begin{theorem} \label{Th.sec:Cnrsk.1} 
 Assume  that  conditions $(\A_\zs{1})$ and $(\H_\zs{1})$ hold. Then, for any $k^*> m^*$ and $\nu\ge 0$, 
\begin{equation} \label{sec:Absrsk.5}
\liminf_{\beta\to 0}  \frac{1}{|\log \beta|}  \inf_{\tau\in \Hc\left(\beta,k^{*},m^{*}\right)}\, \sup_{\nu \ge 0} \Rc_\nu(\tau)\, \ge 
\liminf_{\beta\to 0}  \frac{1}{|\log \beta|}  \inf_{\tau\in \Hc\left(\beta,k^{*},m^{*}\right)}\, \Rc_\nu(\tau)\, \ge \frac{1}{I}
\end{equation}
and
\begin{equation} \label{sec:Cnrsk.5}
\liminf_{\beta\to 0} \frac{1}{|\log\beta|} \inf_{\tau\in \Hc^{*}\left(\beta,k^{*},m^{*}\right)}\, \sup_{\nu\ge 0} \Rc^{*}_\nu(\tau)\, \ge 
\liminf_{\beta\to 0} \frac{1}{|\log\beta|} \inf_{\tau\in \Hc^{*}\left(\beta,k^{*},m^{*}\right)}\, \Rc^{*}_\nu(\tau)\, \ge \frac{1}{I}\,.
\end{equation}
\end{theorem}

\proof
By Proposition~\ref{Pr.sec:Absrsk.1} (see Appendix~\ref{A}), for all $\nu\ge 0$ and for sufficient small $\beta>0$
(for which the conditions of this proposition hold)
$$
\inf_\zs{\tau\in \Hc\left(\beta,k^{*},m^{*}\right)}\, \Rc_\nu(\tau)\,\ge \inf_\zs{\tau\in \Delta(\alpha_\zs{1,\beta},\varrho_\zs{1,\beta})}\,\Rc_\nu(\tau)\,.
$$ 
Now inequality \eqref{sec:Bay.1} and condition $(\H_\zs{1}$) imply immediately \eqref{sec:Absrsk.5}.

Proposition~\ref{Pr.sec:Cnrsk.1} (see Appendix~\ref{A}) implies that for all $\nu\ge 0$  and for a sufficient small $\beta>0$ (for which the conditions of this proposition hold)
 $$
\inf_\zs{\tau\in \Hc\left(\beta,k^{*},m^{*}\right)}\,\Rc^{*}_\nu(\tau)\,\ge \inf_\zs{\tau\in \Delta(\alpha^{*}_\zs{1,\beta},\varrho_\zs{1,\beta})}\,\Rc^{*}_\nu(\tau)\,.
$$ 
Inequality \eqref{sec:Bay.2} and  condition $(\H_\zs{1}$) imply immediately \eqref{sec:Cnrsk.5}. 
\endproof

\subsection{Asymptotic optimality of the Shiryaev--Roberts procedure} \label{ssec:AOminimax}

To establish asymptotic optimality properties of the SR procedure with respect to the risks $\Rc_\nu(\tau)$ (for all $\nu\ge 0$) and $\sup_{\nu\ge0}\Rc_\nu(\tau)$ in the class $\Hc\left(\beta,k^{*},m^{*}\right)$ we 
need the uniform  complete convergence condition $(\A_\zs{2}$) as well as  the following condition. \vspace{3mm}

\noindent $(\H_\zs{2})$  {\em 
Parameter $k^{*}$ in \eqref{sec:Absrsk.1} is a function of $\beta$, i.e. $k^{*}=k^{*}_\zs{\beta}$, such that
\begin{equation}\label{sec:Absrsk.6}
\lim_\zs{\beta\to 0}\,\frac{|\log\alpha_\zs{2,\beta}|}{|\log\beta|}=1\,,
\end{equation}
where $\alpha_\zs{2,\beta}=\alpha_\zs{2}(\beta,k^{*}_\zs{\beta})$.}

\noindent We can choose, for example, 
\begin{equation}\label{sec:Absrsk.9}
m^{*}_\zs{\beta}=1+\lfloor (1+|\log\beta|)^{2} \rfloor\quad\mbox{and}\quad k^{*}_\zs{\beta}=2m^{*}_\zs{\beta}\,.
\end{equation}

Next, denote by $T_\beta$ the SR procedure $T(h_\beta)$ defined in \eqref{sec:Prbf.7-00} with the threshold $h_\beta$ given by
\begin{equation}\label{sec:Absrsk.7}
h_\zs{\beta}= \frac{1-\alpha_\zs{2,\beta}}{\varrho_\zs{2,\beta} \alpha_\zs{2,\beta}} \, . 
\end{equation}

The following theorem establishes first-order asymptotic optimality of the SR procedure $T_\beta$ with respect to the risks $\Rc_\nu(\tau)$ and $\sup_{\nu\ge0}\Rc_\nu(\tau)$ 
in the class $\Hc\left(\beta,k^{*},m^{*}\right)$ as $\beta\to0$, i.e., $T_\beta$ is an asymptotic solution of the problems \eqref{sec:Prbf.5-0}  and \eqref{sec:Prbf.5} as the LPFA vanishes. 

\begin{theorem}\label{Th.sec:Absrsk.2} 
If conditions $(\H_\zs{1})$ and $(\H_\zs{2})$ hold, then, for any $0<\beta<1$, the SR procedure $T_\beta$ with the threshold $h_\beta$ given by \eqref{sec:Absrsk.7} belongs to the class $\Hc\left(\beta,k^{*},m^{*}\right)$. 
If, in addition, condition $(\A_\zs{2})$ is satisfied,  then the SR procedure $T_\beta$  is first-order 
asymptotically uniformly pointwise  optimal and minimax in the class $\Hc\left(\beta,k^{*},m^{*}\right)$, i.e.,
\begin{equation}\label{sec:Absrsk.8-001}
\lim_\zs{\beta\to 0}\,\dfrac{\inf_\zs{\tau\in \Hc(\beta, k^{*},m^{*})}\Rc_\nu(\tau)}{\Rc_\nu(T_\zs{\beta})} =1 \quad \text{for all}~ \nu \ge 0
\end{equation}
and
\begin{equation}\label{sec:Absrsk.8-00}
\lim_\zs{\beta\to 0}\,\dfrac{\inf_\zs{\tau\in \Hc(\beta, k^{*},m^{*})}\sup_\zs{\nu\ge 0} \Rc_\nu(\tau)}{\sup_\zs{\nu\ge 0}\Rc_\nu(T_\zs{\beta})} =1 \,.
\end{equation}
\noindent Also, as $\beta\to0$, the following first-order asymptotic approximations hold for the pointwise and maximal risks:
\begin{equation}\label{sec:SRAO1}
 \Rc_\nu(T_\zs{\beta}) ~  \sim~   \inf_\zs{\tau\in \Hc\left(\beta,k^{*},m^{*}\right)}\, \Rc_\nu(\tau)~ \sim~ \frac{|\log \beta|}{I} \quad \text{for any}~ \nu\ge 0
\end{equation}
and
\begin{equation}\label{sec:SRAO1b}
\sup_{\nu\ge 0} \Rc_\nu(T_\zs{\beta}) ~  \sim~   \inf_\zs{\tau\in \Hc\left(\beta,k^{*},m^{*}\right)}\, \sup_{\nu\ge 0} \Rc_\nu(\tau)~ \sim~ \frac{|\log \beta|}{I} \, .
\end{equation}
\end{theorem}

\proof
By Lemma~\ref{Lem:PFASR}, the SR procedure $T(\alpha,\varrho) \in\Delta(\alpha,\varrho)$ for any $0<\alpha, \varrho<1$.
Moreover, note that the definition \eqref{sec:Absrsk.7} yields $T_\zs{\beta}=T(\alpha_\zs{2,\beta},\varrho_\zs{2,\beta})$,
i.e., $T_\zs{\beta}\in \Delta(\alpha_\zs{2,\beta},\varrho_\zs{2,\beta})$. Using Proposition~\ref{Pr.sec:Absrsk.1},
we obtain that  $T_\zs{\beta}\in \Hc\left(\beta,k^{*},m^{*}\right)$ for any $0<\beta<1$. Furthermore,  condition $(\H_\zs{2}$)
and the definition of $\varrho_\zs{2,\beta}$ in \eqref{sec:Absrsk.1-01} imply directly that $\lim_\zs{\beta\to 0}\,\log h_\zs{\beta}/|\log \beta|=1$.
Thus, the asymptotic upper bound \eqref{AAAproxa} (with $r=1$) in Proposition~\ref{Prop: A1} implies the following upper bound
\begin{equation*}
\limsup_{\beta\to 0}\, \frac{1}{|\log\beta|}\sup_{\nu\ge 0}\, \Rc_\nu(T_\zs{\beta})\,\le \frac{1}{I}\, .
\end{equation*}
The asymptotic equalities \eqref{sec:Absrsk.8-001}  and \eqref{sec:Absrsk.8-00} follow immediately from this upper bound and the lower bounds \eqref{sec:Absrsk.5} 
in Theorem~\ref{Th.sec:Cnrsk.1}. The asymptotic expansions \eqref{sec:SRAO1} and \eqref{sec:SRAO1b}  are obvious.
\endproof

Now we define
\begin{equation}\label{sec:Cnrsk.1}
\alpha_\zs{3}=
\alpha_\zs{3}(\beta,k^{*})=
\frac{\beta(1-\varrho_\zs{2,\beta})^{k^{*}}}
{1+\beta}
\,,
\end{equation}
where  the function  $\varrho_\zs{2,\beta}$ is defined in 
\eqref{sec:Absrsk.1-01}.

To prove asymptotic optimality in the class $\Hc^*\left(\beta,k^{*},m^{*}\right)$ with respect to the risk $\Rc^*_\nu(\tau)$ we need the following condition. \vspace{3mm}

\noindent $(\H_\zs{3})$  {\em 
Parameters $k^{*}$ and $m^{*}$ are functions of $\beta$, i.e. $k^{*}=k^{*}_\zs{\beta}$ and $m^{*}=m^{*}_\zs{\beta}$, such that
\begin{equation}\label{sec:Cnrsk.6}
\lim_\zs{\beta\to 0}\,
\left( |\log\alpha_\zs{3,\beta}|+k^{*}_\zs{\beta} \log(1-\varrho_\zs{2,\beta})\right)
=+\infty
\quad\mbox{and}\quad
\lim_\zs{\beta\to 0}\frac{ |\log\alpha_\zs{3,\beta}|}{|\log\beta|}=1\,.
\end{equation}
where $\alpha_\zs{3,\beta}=\alpha_\zs{3}(\beta,k^{*}_\zs{\beta})$.} 
 
 \noindent We can take, for example, the parameters $k^{*}=k^{*}_\zs{\beta}$ and $m^{*}=m^{*}_\zs{\beta}$ as in \eqref{sec:Absrsk.9}.

Denote by $T_\beta^*$ the SR procedure $T(h_\beta^*)$ defined in \eqref{sec:Prbf.7-00} with the threshold $h_\beta^*$ given by
\begin{equation}\label{sec:Cnrsk.7}
 h^{*}_\zs{\beta}= \frac{1-\alpha_\zs{3,\beta}}{\varrho_\zs{2,\beta}\alpha_\zs{3,\beta}}\,. 
\end{equation}

\begin{theorem}\label{Th.sec:Cnrsk.2} 
If conditions  $(\H_\zs{1})$ and $(\H_\zs{3})$ hold, then, for any $0<\beta<1$, the SR procedure $T^*_\beta$  with the threshold $h_\beta^*$ given by \eqref{sec:Cnrsk.7} belongs to the 
class $\Hc^{*}\left(\beta,k^{*},m^{*}\right)$. Assume that in addition condition $(\A_\zs{2})$ is satisfied. Then the SR procedure $T_\beta^*$ is first-order 
asymptotically uniformly poitwise optimal and minimax in the class $\Hc^*\left(\beta,k^{*},m^{*}\right)$, i.e.,
\begin{equation}\label{sec:Cbsrsk.8-00a}
\lim_\zs{\beta\to 0}\,\dfrac{\inf_\zs{\tau\in \Hc^*(\beta, k^{*},m^{*})} \Rc_\nu^*(\tau)}{\Rc^*_\nu(T_\zs{\beta}^*)} =1 \quad \text{for all fixed}~\nu\ge 0\,.
\end{equation}
and
\begin{equation}\label{sec:Cbsrsk.8-00}
\lim_\zs{\beta\to 0}\,\dfrac{\inf_\zs{\tau\in \Hc^*(\beta, k^{*},m^{*})} \max_\zs{0\le \nu\le k^{*}_\beta}\Rc_\nu^*(\tau)}{\max_\zs{0\le \nu\le k^{*}_\beta}\Rc^*_\nu(T_\zs{\beta}^*)} =1 \,.
\end{equation}
\noindent 
\noindent Also, as $\beta\to0$, the following first-order asymptotic approximations hold for the pointwise and maximal risks:
\begin{equation}\label{sec:SRAO2}
 \Rc_\nu^*(T_\zs{\beta}^*) ~  \sim~   \inf_\zs{\tau\in \Hc\left(\beta,k^{*},m^{*}\right)}\, \Rc_\nu^*(\tau)~ \sim~ \frac{|\log \beta|}{I} \quad \text{for any}~ \nu\ge 0
\end{equation}
and
\begin{equation}\label{sec:SRAO2b}
\sup_{0\le \nu\le k^*_\beta} \Rc_\nu^*(T_\zs{\beta}^*) ~  \sim~   \inf_\zs{\tau\in \Hc\left(\beta,k^{*},m^{*}\right)}\, \sup_{0\le \nu\le k^*_\beta} \Rc_\nu^*(\tau)~ \sim~ \frac{|\log \beta|}{I} \, .
\end{equation}
\end{theorem}

\proof
By Lemma~\ref{Lem:PFASR}, the SR procedure $T(\alpha,\varrho) \in\Delta(\alpha,\varrho)$ for any $0<\alpha, \varrho<1$.
Now, note that the definition \eqref{sec:Cnrsk.7} yields $T^*_\zs{\beta}=T(\alpha_\zs{3,\beta},\varrho_\zs{2,\beta})$,
i.e., $T^*_\zs{\beta}\in \Delta(\alpha_\zs{3,\beta},\varrho_\zs{2,\beta})$. Using Proposition~\ref{Pr.sec:Cnrsk.1},
we obtain that  the stopping time $T^*_\zs{\beta}$ belongs to $\Hc^{*}\left(\beta,k^{*},m^{*}\right)$ for any $0<\beta<1$. 

Next, in view of the definition of $h^{*}_\zs{\beta}$ in \eqref{sec:Cnrsk.7} and of the form of the function $\varrho_\zs{2,\beta}$ in \eqref{sec:Absrsk.1-01}
we obtain, using condition $(\H_\zs{3}$), that $\lim_\zs{\beta\to 0}\,\log h^{*}_\zs{\beta}/|\log\beta|=1$. Thus, by \eqref{AAAproxb} (with $r=1$) in Proposition~\ref{Prop: A1},
\[
\limsup_{\beta\to\infty}\, \frac{1}{|\log \beta|} \, \Rc_\nu^*(T^*_\beta) \le \frac{1}{I} \quad \text{for all}~ \nu \ge 0.
\]
Comparing to the reverse inequality \eqref{sec:Absrsk.5} implies \eqref{sec:Cbsrsk.8-00a}. Asymptotic approximations \eqref{sec:SRAO2} are obvious from \eqref{sec:Absrsk.5} and  \eqref{sec:Cbsrsk.8-00a}.

Using inequality \eqref{sec:Bay.5} and  condition $(\H_\zs{3}$), we obtain
$$
\Pb_\zs{\infty}\left(T^{*}_\zs{\beta}\le k^{*}_\zs{\beta}\right) \le e^{\log \alpha_\zs{3,\beta}-k^{*}_\zs{\beta}\log(1-\varrho_\zs{2,\beta})}
\to 0 \quad\mbox{as}\quad \beta\to 0\,.
$$
Therefore,
\[
\begin{aligned}
 \min_\zs{0\le k\le k^{*}_\beta}\,\Pb_\zs{k}\left(T^{*}_\zs{\beta} > k \right)
&=\min_\zs{0\le k\le k^{*}_\beta}\,\Pb_\zs{\infty}\left(T^{*}_\zs{\beta}>  k \right)\\
&=\Pb_\zs{\infty}\left(T^{*}_\zs{\beta}> k^{*}_\zs{\beta} \right)
=1-\Pb_\zs{\infty}\left(T^{*}_\zs{\beta} \le k^{*}_\zs{\beta} \right) \to 1\quad\mbox{as}\quad \beta\to 0\,.
\end{aligned}
\]
Note that the maximal risk $\max_\zs{0\le \nu \le k^{*}_\beta}\,\Rc^{*}_\nu(T^{*}_\zs{\beta})$ can be estimated as
$$
\max_{0\le \nu \le k^{*}_\beta}\,\Rc^{*}_\nu(T^{*}_{\beta})\le \dfrac{\max_{0\le \nu \le k^{*}_\beta}\, \Rc_\nu(T^{*}_{\beta})}{\min_{0\le \nu \le k^{*}_\beta}\,\Pb_{\infty}\left(T^{*}_{\beta}> \nu \right)}\,.
$$
Asymptotic equality \eqref{AAAprox} with $r=1$ in Proposition~\ref{Prop: A1} implies that 
\[
\max_{0\le \nu \le k^{*}_\beta}\, \Rc_\nu(T^{*}_{\beta}) \le \sup_{0\le \nu <\infty}\, \Rc_\nu(T^{*}_{\beta}) = \frac{\log h_\beta^*}{I} (1+o(1)) \quad \text{as}~ \beta\to0.
\]
Since, as we mentioned above, $\lim_{\beta\to 0}\,\log h^{*}_{\beta}/|\log\beta|=1$, we obtain the upper bound 
\begin{equation}\label{UpperTbeta}
\limsup_{\beta\to\infty}\, \frac{1}{|\log \beta|} \, \max_{0\le \nu \le k^{*}_\beta} \Rc_\nu^*(T^*_\beta) \le \frac{1}{I}.
\end{equation}
Asymptotic equalities  \eqref{sec:Cbsrsk.8-00} now follow from the upper bound \eqref{UpperTbeta} and the lower bound \eqref{sec:Cnrsk.5}. 
Asymptotic approximations \eqref{sec:SRAO2b} are obvious from \eqref{sec:Cnrsk.5} and  \eqref{sec:Cbsrsk.8-00}. The proof is complete.
\endproof

\begin{remark} \label{Rem:Lai}
We recall that Lai's condition (6) in \cite{LaiIEEE98} for the asymptotic lower bound 
\[
\liminf_{\gamma \to\infty} \frac{1}{\log \gamma} \inf_{\tau \in \Hc_\gamma}  \ESADD(\tau) \ge \frac{1}{I}
\]
in the class $\Hc_\gamma=\{\tau: \EV \tau \ge \gamma\}$ is the following:  
\begin{equation}
\label{sec:LaiCn.0-0}
\lim_\zs{n\to\infty}\,\sup_\zs{\nu \ge 0} \, \esssup\, \Pb_\zs{\nu} \left(\max_{1 \le i \le n}Z^{\nu}_\zs{\nu+i} \ge I(1+\varepsilon)n\, \vert\, \Fc_\nu \right) =0 \quad \text{for all}~\varepsilon>0\,,
\end{equation}
where the parameter $I$ is given in condition $(\A_\zs{1}$). Clearly, condition \eqref{sec:LaiCn.0-0} is much stronger than the a.s.\ convergence condition $(\A_\zs{1}$) 
required in Theorem~\ref{Th.sec:Cnrsk.1}, and it does not hold in many important practical cases. 
Also, Lai's condition (24) in \cite{LaiIEEE98} for asymptotic optimality of the CUSUM procedure in the classes $\Hc_\gamma$ and 
$\Hc(\beta)=\{\tau: \sup_{k \ge 1} \Pb_\infty (k \le \tau < k +m^*_\beta) \le \beta\}$ is: 
\begin{equation}
\label{sec:LaiCn.0-1}
\lim_\zs{n\to\infty}\,\sup_\zs{\ell \ge \nu}\, \esssup\, \Pb_\zs{\nu} \left(\wt{Z}_\zs{\ell,n} \le -\varepsilon \,\vert\, \Fc_\ell \right)=0  \quad \text{for all}~\varepsilon>0\, .
\end{equation}
Typically this condition is more difficult to check than the uniform complete convergence condition $(\A_\zs{2}$) required in Theorem~\ref{Th.sec:Cnrsk.2}, which in fact can be relaxed to
\[
\sum_{n=1}^\infty\, \sup_\zs{\nu\ge 0}\, \Pb_\zs{\nu}\Bigl\{\wt{Z}_\zs{\nu,n} < -\vae \Bigr\}<\infty 
\] 
(see Remark~\ref{Rem:Lefttail}). In addition, for certain models condition \eqref{sec:LaiCn.0-1} does not hold, while condition 
$(\A_\zs{2}$)  holds (see, e.g., an example in Subsection~\ref{ssec:LaiCn} below). On the other hand, in the iid case condition \eqref{sec:LaiCn.0-1} is less stringent than $(\A_\zs{2}$). 
\end{remark}

As in Theorem~\ref{Th.sec:Bay.rcomp}, the results of Theorem~\ref{Th.sec:Absrsk.2} and Theorem~\ref{Th.sec:Cnrsk.2} can be extended to higher moments of the detection delay by strengthening the 
complete convergence with the uniform $r-$complete convergence \eqref{sec:UCCCr}. More specifically, the following asymptotic optimality result holds true. 

\begin{theorem}\label{Th.sec:Cnrskr.2} 
Assume that conditions $(\H_\zs{1})$ and $(\H_\zs{3})$ hold, and in addition, for some $r>1$ the uniform $r$-complete convergence condition \eqref{sec:UCCCr} is satisfied.
Then, for any $0<\beta<1$, the SR procedure $T^*_\beta$ with the threshold $h_\beta^*$ given by \eqref{sec:Cnrsk.7} belongs to the class $\Hc^{*}\left(\beta,k^{*},m^{*}\right)$ and as $\beta\to0$ for any $0<\ell \le r$
\begin{equation}\label{sec:Cbsrskr.8-00}
\begin{aligned}
 \EV_{\nu}\left[(T^{*}_{\beta}-\nu)^\ell | T^{*}_{\beta}>\nu\right] & ~ \sim~  \inf_{\tau\in \Hc^*\left(\beta,k^{*},m^{*}\right)}\, \EV_{\nu}\left[(\tau-\nu)^\ell | \tau>\nu \right]
 \\ 
 &~ \sim~ \brc{\frac{|\log \beta|}{I}}^\ell \quad \text{for all}~\nu\ge 0
 \end{aligned}
\end{equation}
and
\begin{equation}\label{sec:Cbsrskr.8-00b}
\begin{aligned}
\max_{0 \le \nu \le k^{*}_\beta}\, \EV_{\nu}\left[(T^{*}_{\beta}-\nu)^\ell | T^{*}_{\beta}>\nu\right] & ~ \sim~  
\inf_{\tau\in \Hc^*\left(\beta,k^{*},m^{*}\right)}\, \max_{0 \le \nu \le k^{*}_\beta}\,  \EV_{\nu}\left[(\tau-\nu)^\ell | \tau>\nu \right]
\\
& ~ \sim~ \brc{\frac{|\log \beta|}{I}}^\ell\,.
\end{aligned}
\end{equation}
Therefore, the SR procedure $T^*_\beta$ is first-order asymptotically uniformly  pointwise optimal and also minimax in the class $\Hc^*\left(\beta,k^{*},m^{*}\right)$ with respect to the moments of the detection delay 
up to order $r$.
\end{theorem}

\proof
The facts that $T^*_{\beta}\in\Hc^{*}\left(\beta,k^{*},m^{*}\right)$ for any $0<\beta<1$ and that $\log h^{*}_{\beta}\sim |\log\beta|$ as $\beta\to0$ were established in Theorem~\ref{Th.sec:Cnrsk.2}. 
Now, using \eqref{AAAproxb} in Proposition~\ref{Prop: A1} (along with the equality $\lim_{\beta\to0} \Pb_\infty(T^*_\beta>\nu)=1$, $\nu\ge0$), we obtain the upper bound
\[
 \EV_\zs{\nu}\left[(T^{*}_\zs{\beta}-\nu)^r | T^{*}_\zs{\beta}>\nu\right] \le \brc{\frac{|\log \beta|}{I}}^r(1+o(1)) \quad \text{as}~\beta\to0.
\]
Jensen's inequality and the lower bound \eqref{sec:Cnrsk.5} yield, for any $r\ge 1$ and $\nu\ge 0$,
\[
\inf_\zs{\tau\in \Hc^*\left(\beta,k^{*},m^{*}\right)}\, \sup_\zs{\nu \ge 0}\,  \EV_\zs{\nu}\left[(\tau-\nu)^r | \tau>\nu \right] \ge 
\inf_\zs{\tau\in \Hc^*\left(\beta,k^{*},m^{*}\right)}\,\EV_\zs{\nu}\left[(\tau-\nu)^r | \tau>\nu \right] \ge \brc{\frac{|\log \beta|}{I}}^r(1+o(1)), 
\]
which along with the previous upper bound proves \eqref{sec:Cbsrskr.8-00}.

To prove \eqref{sec:Cbsrskr.8-00b} it suffices to show that
\begin{equation}\label{UBbeta}
\limsup_{\beta\to0} \frac{\max_{0\le \nu \le k^{*}_\beta}\,\EV_{\nu}\left[(T^{*}_{\beta}-\nu)^r | T^{*}_{\beta}>\nu\right]}{|\log\beta|^r} \le \frac{1}{I^r} \, .
\end{equation}
Note that
$$
\max_{0\le \nu \le k^{*}_\beta}\, \EV_{\nu}\left[(T^{*}_{\beta}-\nu)^r | T^{*}_{\beta}>\nu\right] \le 
\dfrac{\max_{0\le \nu \le k^{*}_\beta}\,  \EV_{\nu}\left[(T^{*}_{\beta}-\nu)^+]^r\right]}{\min_{0\le \nu \le k^{*}_\beta}\,\Pb_{\infty}\left(T^{*}_\zs{\beta}> \nu \right)},
$$
where 
\[
\min_{0\le \nu \le k^{*}_\beta}\,\Pb_{\infty}\left(T^{*}_{\beta}> \nu \right) = \Pb_{\infty}\left(T^{*}_{\beta}> k^*_\beta \right) \to 1 \quad \text{as}~\beta\to0 \, .
\]
As a result, using \eqref{AAAprox} in Proposition~\ref{Prop: A1}, we obtain
\[
\begin{aligned}
\max_{0\le \nu \le k^{*}_\beta}\,\EV_{\nu}\left[(T^{*}_{\beta}-\nu)^r | T^{*}_{\beta}>\nu\right] & \le \dfrac{\sup_{0\le \nu <\infty}\,  \EV_{\nu}\left[(T^{*}_{\beta}-\nu)^+]^r\right]}{\Pb_{\infty}\left(T^{*}_{\beta}> k^*_\beta \right)}
\\
&=\dfrac{(\log h^*_\beta /I)^r(1+o(1))}{\Pb_{\infty}\left(T^{*}_{\beta}> k^*_\beta \right)} = \brc{\frac{|\log \beta|}{I}}^r(1+o(1)).
\end{aligned}
\]
This obviously yields  the upper bound \eqref{UBbeta} and the proof is complete.
\endproof

\begin{remark}\label{Rem:Lefttail}
The uniform $r-$complete convergence condition \eqref{sec:UCCCr} can be relaxed to the following one-sided version: for some $r>1$ and any  $\vae>0$
\[
\sum_{n=1}^\infty\, n^{r-1} \, \sup_\zs{\nu\ge 0}\, \Pb_\zs{\nu}\Bigl\{\wt{Z}_\zs{\nu,n} < -\vae \Bigr\}<\infty .
\]
In this case, one needs to additionally require the almost sure convergence condition $(\A_\zs{1}$), which guarantees condition \eqref{Pmax} in Lemma~\ref{Prop: A1}.
\end{remark}

%
\section{Concentration inequalities for functions of  homogeneous Markov processes}
\label{sec:Mrk}
In this section, we obtain certain sufficient conditions for homogeneous Markov processes in order to verify condition $(\A_\zs{2})$ for this class of processes.

Let  $(X_\zs{n})_\zs{n\ge 1}$ be a time homogeneous  Markov process with values in a measurable space $(\Xc,\Bc)$ with the transition probability 
$P(x,A)$ defined in \eqref{subsec:AMc.1-1}.  In the sequel, we denote by $\EV_\zs{x}(\cdot)$   the  expectation 
with respect to this probability. In addition, we assume that this process is geometrically ergodic, i.e., \vspace{3mm}

\noindent $(\B_\zs{1})$
{\em Assume that there exist positives constants $0<R<\infty$, $\kappa>0$, probability measure $\lambda$ on
$(\Xc,\Bc)$
  and the  Lyapunov $\Xc\to [1,\infty)$
function $\V$ with $\lambda(\V)<\infty$,  such that
$$
\sup_\zs{n\ge 0}\,
e^{\kappa n}\,
\sup_\zs{0<f\le \V}\,
\sup_\zs{x\in\bbr}\frac{1}{\V(x)}\,
\left|
\EV_\zs{x}\,f(X_\zs{n})
-\lambda(f)
\right|
\le\, 
R\,.
$$
}
Now, for some $\r>0$, we set
\begin{equation}\label{sec:Mrk.3}
\upsilon^{*}_\zs{\r}(x)=\sup_\zs{n\ge 0}\,\EV_\zs{x}\,\left(\V(X_\zs{n})\right)^{\r}\,.
\end{equation}
\noindent
Let $g$ be a measurable $\Xc\times\Xc\to\bbr$ function such that the following integrals exist
\begin{equation}\label{sec:Mrk.1}
\wt{g}(u)=\int_\zs{\Xc}\,g(v,u)\,P(u,\d v)
\quad\mbox{and}\quad
\lambda(\wt{g})
=\int_\zs{\Xc\times\Xc}\,g(v,u)\,P(u,\d v)\,\lambda(\d u)
\,.
\end{equation}

\noindent $(\B_\zs{2})$
{\em Assume that the function $g$ is such that $|\wt{g}(x)|\le \V(x)$ for all $x\in\Xc$.}

We study the concentration properties for the process
$W_\zs{n}(g)=\sum^{n}_\zs{j=1}\,g(X_\zs{j},X_\zs{j-1})$,
or equivalently the properties of the deviation
$\wt{W}_\zs{n}(g)=n^{-1}W_\zs{n}(g) -\lambda(\wt{g})$.

Similarly to \eqref{sec:Mrk.3}, we define for some $\r>0$ 
\begin{equation}\label{sec:Mrk.2}
g^{*}_\zs{\r}(x)=\sup_\zs{j\ge 1}\,
\EV_\zs{x}\,|g(X_\zs{j},X_\zs{j-1})|^{\r}
\,.
\end{equation}

\noindent
Now we set
\begin{equation}\label{sec:Mrk.4}
W^{*}_\zs{\r}=4^{\r-1}\,
\left(
1+\u^{*}_\zs{\r}+|\lambda(\wt{g})|^{\r}+
(16\r)^{\r/2}
\right)
\quad\mbox{and}\quad
\u^{*}_\zs{\r}=\left(\frac{2\r R^{2}\,e^{\kappa}}{e^{\kappa}-1}\right)^{\r/2}
\,.
\end{equation}

\begin{proposition} \label{Pr.sec:Mrk.1} 
Assume that conditions 
$(\B_\zs{1})$ and $(\B_\zs{2})$ hold. Then for any $x\in\Xc$ and  $\r\ge 2$, for which
$\upsilon^{*}_\zs{\r}(x)<\infty$ and $g^{*}_\zs{\r}(x)<\infty$, one has
\begin{equation}\label{sec:Mrk.4E}
\EV_\zs{x}|\wt{W}_\zs{n}(g)|^{\r}\,\le\,W^{*}_\zs{\r}\,
\frac{\left(1+\upsilon^{*}_\zs{\r}(x)+g^{*}_\zs{\r}(x)\right)}{n^{\r/2}} \quad \text{for any}~ n\ge 2.
\end{equation}
\end{proposition}

\proof
Note that  we can represent the term $W_\zs{n}(g)$ as 
\begin{equation}\label{sec:Mrk.4-1n}
W_\zs{n}(g)
=(n-1)\,\lambda(\wt{g})+\wt{g}(x)+\U_\zs{n-1}+\M_\zs{n}\,, 
\end{equation}

where
$$
\U_\zs{n}=\sum^{n}_\zs{j=1}\,\left(\wt{g}(X_\zs{j})-\lambda(\wt{g})\right)
:=\sum^{n}_\zs{j=1}\,u_\zs{j}
\quad\mbox{and}\quad
\M_\zs{n}=\sum^{n}_\zs{j=1}\,\left(g(X_\zs{j},X_\zs{j-1})-\wt{g}(X_\zs{j-1})\right)
\,.
$$
\noindent 
To estimate the powers of the $U_\zs{n}$
we need to estimate the corresponding coefficient $b_\zs{j,n}(\r)$
from Proposition~\ref{Pr.sec:Bi.1} (see Appendix~\ref{B}). To this end, note that for $j\ge l$
$$
\EV_\zs{x}\left(u_\zs{j}|\Fc_\zs{l}\right)=
\EV_\zs{x}\left(u_\zs{j}|X_\zs{1},\ldots,X_\zs{l}\right)
=\wt{\omega}_\zs{l-j}(X_\zs{l}) \, ,
$$
where $\wt{\omega}_\zs{m}(x)=\EV_\zs{x}\,\wt{g}(X_\zs{m})-\lambda(\wt{g})$.
Now, by condition $(\B_\zs{1})$,  for any $x\in \Xc$ and any $m\ge 0$, $|\wt{\omega}_\zs{m}(x)|\le R \, \V(x)\,e^{-\kappa m}$,
i.e., for any $j\ge l\ge 1$
$$
|\EV_\zs{x}\left(u_\zs{j}|\Fc_\zs{l}\right)|\le R ~ \V(X_\zs{l})\,e^{-\kappa (j-l)}\,.
$$
In particular, we have $|u_\zs{l}|\le  R ~ \V(X_\zs{l})$. 
Therefore, the  coefficients \eqref{subsec:AMc.1}  can be estimated as
$$
b_\zs{j,n}(\r)\le \frac{R^{2}\,e^{\kappa}}{e^{\kappa}-1}
\left( \upsilon^{*}_\zs{\r}(x)\right)^{2/\r}
$$
and by  Proposition~\ref{Pr.sec:Bi.1} we get $\EV_\zs{x}\,|\U_\zs{n}|^{\r}\le \u^{*}_\zs{\r}\,
\upsilon^{*}_\zs{\r}(x)\,n^{\r/2}$, 
where $\u^{*}_\zs{\r}$ is defined in \eqref{sec:Mrk.3}.
Similarly, to estimate the martingale $\M_\zs{n}$ we make use of 
 Proposition~\ref{Pr.sec:Bi.1}. Note that in this case the coefficient
\eqref{subsec:AMc.1} has the form $b_\zs{j,n}(\r)=(\EV_\zs{x} |g(X_\zs{j},X_\zs{j-1})-\wt{g}(X_\zs{j-1})|^{\r})^{2/\r}$, 
and it can be estimated for $j\ge 1$ as
$$
b_\zs{j,n}(\r)\,
\le 
\left(2^{\r-1}
\left(
\EV_\zs{x}\,
\left|
g(X_\zs{j},X_\zs{j-1})\right|^{\r}
+
\EV_\zs{x}\,
\left|
\wt{g}(X_\zs{j-1})
\right|^{\r}
\right)
 \right)^{2/\r}\,.
$$
Taking into account  Jensen's inequality and the definition \eqref{sec:Mrk.2}, we obtain that $b_\zs{j,n}(\r)\le 4 (g^{*}_\zs{\r}(x))^{2/\r}$, and therefore, from Proposition~\ref{Pr.sec:Bi.1}
it follows that for $n\ge 1$, $\EV_\zs{x}\,|\M_\zs{n}|^{\r}\le (8\r)^{\r/2}\,g^{*}_\zs{\r}(x) n^{\r/2}$.
Therefore, taking into account that
$$
|\wt{W}_\zs{n}(g)|\le 
\frac{|\wt{g}(x)|+
|\lambda(\wt{g})|+
|\U_\zs{n-1}|
+
|\M_\zs{n}|
}{n}\,,
$$
\noindent 
we obtain,  for any $n\ge 1$, 
\begin{align*}
\EV_\zs{x}|\wt{W}_\zs{n}(g)|^{\r}\,
&\le
\,\frac{4^{\r-1}}{n^{\r/2}}\,
\left(|\wt{g}(x)|^{\r}
+|\lambda(\wt{g})|^{\r}
+\u^{*}_\zs{\r}\,
\upsilon^{*}_\zs{\r}(x)+
(8\r)^{\r/2}\,g^{*}_\zs{\r}(x)
\right)\\[2mm]
&\le
\,\frac{4^{\r-1}}{n^{\r/2}}\,
\left((1+\u^{*}_\zs{\r})\,\upsilon^{*}_\zs{\r}(x)
+|\lambda(\wt{g})|^{\r}\,+
(8\r)^{\r/2}\,g^{*}_\zs{\r}(x)
\right)\\[2mm]
&
\le\,
W^{*}_\zs{\r}\,
\frac{\left(1+\upsilon^{*}_\zs{\r}(x)+g^{*}_\zs{\r}(x)\right) }{n^{\r/2}}
\,. 
\end{align*}
\noindent 
Hence Proposition~\ref{Pr.sec:Mrk.1}.
\endproof

As we will see later in Section \ref{sec:Ex}, condition $(\B_\zs{1})$ does not hold directly for some time series. 
For this reason, we introduce the following modification of this condition. \vspace{3mm}
 
\noindent $(\B'_\zs{1})$
{\em Assume that there is  some integer $p\ge 1$ such that for any $0\le \iota\le p-1$
there exist positive constants $0<R_\zs{\iota}<\infty$, $\kappa_\zs{\iota}>0$,
probability measure $\lambda_\zs{\iota}$ on $\Xc$ and the  Lyapunov $\Xc\to [1,\infty)$
function $\V_\zs{\iota}$ with $\lambda_\zs{\iota}(\V_\zs{\iota})<\infty$,  such that
$$
\sup_\zs{l\ge 0}\,
e^{\kappa_\zs{\iota} l}\,
\sup_\zs{0<f\le \V_\zs{\iota}}\,
\sup_\zs{x\in\bbr}\frac{1}{\V_\zs{\iota}(x)}\,
\left|
\EV_\zs{x}\,f(X_\zs{pl+\iota})
-\lambda_\zs{\iota}(f)
\right|
\le\, 
R_\zs{\iota}\,.
$$
}
Similarly to \eqref{sec:Mrk.3} we introduce 
\begin{equation}\label{sec:Mrk.3-1}
\upsilon^{*}_\zs{\r,\iota}(x)=\sup_\zs{j\ge 0}\,
\EV_\zs{x}\,
\left(\V_\zs{\iota}(X_\zs{pj+\iota})
\right)^{\r}
\quad\mbox{and}\quad
\upsilon^{*}_\zs{\r,\max}(x)=\max_\zs{0\le \iota\le p-1}
\upsilon^{*}_\zs{\r,\iota}(x)
\end{equation}
and impose the following condition. \vspace{3mm}

\noindent $(\B'_\zs{2})$
{\em Assume that the function $g$ defined in \eqref{sec:Mrk.1} is such that 
$|\wt{g}(x)|\le \min_\zs{0\le \iota\le p-1}\V_\zs{\iota}(x)$ for all $x\in\Xc$.
}

Now we set $\overline{W}_\zs{n}(g)= n^{-1}\,W_\zs{n}(g) -\overline \lambda(\wt{g})$,
where $\overline \lambda(g)=(1/p)\sum^{p}_\zs{\iota=0}\,\lambda_\zs{\iota}(g)$.

\begin{proposition} \label{Pr.sec:Mrk.1-1} 
Assume that conditions $(\B'_\zs{1})$ and $(\B'_\zs{2})$ hold. Then for any $x\in\Xc$ and any $\r\ge 2$, for which
$\upsilon^{*}_\zs{\r,\max}(x)<\infty$ and $g^{*}_\zs{\r}(x) < \infty$,
there exists  a constant $\overline{W}^{*}_\zs{\r}>0$ such that 
\begin{equation}\label{sec:Mrk.4-1nE}
\EV_\zs{x}|\overline{W}_\zs{n}(g)|^{\r}\,\le\,\overline{W}^{*}_\zs{\r}\,
\frac{\left(1+\upsilon^{*}_\zs{\r,max}(x)+g^{*}_\zs{\r}(x)\right)}{n^{\r/2}} \quad \text{for any}~n\ge 2\,.
\end{equation}
\end{proposition}

\proof
Note that  the term $W_\zs{n}(g)$ can be represented  as $W_\zs{n}(g)=W_\zs{n,1}(g)+\M_\zs{n}$,
where $W_\zs{n,1}(g)=\sum^{n}_\zs{j=1}\,\wt{g}(X_\zs{j-1})$ and $\M_\zs{n}$ is defined in \eqref{sec:Mrk.4-1n}.
Let now $n-1=mp+r$ for some $0\le r\le p-1$. Thus,
$$
W_\zs{n,1}(g)
=
\sum^{p-1}_\zs{\iota=0}\,
\sum^{m}_\zs{l=0}\,\wt{g}(X_\zs{pl+\iota})
-
\sum^{p-1}_\zs{\iota=r+1}\,
\wt{g}(X_\zs{pm+\iota})
=n\overline{\lambda}(\wt{g})
+
\sum^{p-1}_\zs{\iota=0}\,U_\zs{m,\iota}
-r\overline{\lambda}(\wt{g})
-
\sum^{p-1}_\zs{\iota=r+1}\wt{g}(X_\zs{pm+\iota})\,,
$$
where
$
U_\zs{m,\iota}=\sum^{m}_\zs{l=0}\,
\left( \wt{g}(X_\zs{pl+\iota})
-\lambda(\wt{g})
\right)$.
In just the same way as in the proof of  Proposition~\ref{Pr.sec:Mrk.1}, we obtain that for some constant
$\u^{*}_\zs{\r,\iota}>0$
\begin{equation}\label{sec:Mrk.4-2n}
\EV_\zs{x}\,|\U_\zs{m,\iota}|^{\r}\le \u^{*}_\zs{\r,\iota}\,
\upsilon^{*}_\zs{\r,\iota}(x)\,m^{\r/2}
\le \u^{*}_\zs{\r,\max}\,\upsilon^{*}_\zs{\r,\max}(x)\,n^{\r/2}
\,,
\end{equation}
where $\u^{*}_\zs{\r,\max}=\max_\zs{0\le \iota} \u^{*}_\zs{\r,\iota}$.
Furthermore,
$$
\vert
\overline{W}_\zs{n}(g)|
\vert
\le \frac{r}{n}
\vert \overline{\lambda}(\wt{g})\vert
+\frac{1}{n}\,\sum^{p-1}_\zs{\iota=0}\,\vert \U_\zs{m,\iota}\vert
+\frac{1}{n}
\vert M_\zs{n}\vert
+\frac{1}{n}
\sum^{p-1}_\zs{\iota=r+1}\vert\wt{g}(X_\zs{pm+\iota})\vert\,.
$$
Using the upper bound \eqref{sec:Mrk.4-2n} in this inequality, we obtain the inequality \eqref{sec:Mrk.4-1nE}.
\endproof

We return to the detection problem for Markov processes, assuming
that the sequence $(X_\zs{n})_\zs{n\ge 1}$ is a Markov process, such that
$(X_\zs{n})_\zs{1\le n \le \nu}$ is a homogeneous process
 with the transition (from $x$ to $y$) density $f_\zs{0}(y|x)$  
and $(X_\zs{n})_\zs{n> \nu}$ is 
homogeneous positive  ergodic  with the transition density $f_\zs{1}(y|x)$
and the ergodic (stationary) distribution $\lambda$. 
The densities $f_\zs{0}(y|x)$ and  $f_\zs{1}(y|x)$   are calculated 
with respect to a sigma-finite positive measure $\mu$ on $\Bc$.

In this case, we can represent the process $Z^{k}_\zs{n}$ defined in \eqref{sec:Prbf.7} as
\begin{equation}\label{sec:Mrk.5}
Z^{k}_\zs{n}=\sum^{n}_\zs{j=k+1} g(X_\zs{j},X_\zs{j-1})\,,
\quad
  g(y,x)=\log \frac{f_\zs{1}(y|x)}{f_\zs{0}(y|x)}\,.
\end{equation}
\noindent
Therefore, in this case
\begin{equation}\label{sec:Mrk.6}
\wt{g}(x)=
\int_\zs{\Xc}\,g(y,x)\,f_\zs{1}(y|x)\,\mu(\d y)
\,.
\end{equation}

We now formulate the conditions that are sufficient for the main condition $(\A_2)$ to hold in the case of Markov processes. We write $\EV_\zs{x,0}$  for the expectation with respect to
 the  distribution $\Pb_\zs{x,0}(\cdot)=\Pb_\zs{0}(\cdot|X_\zs{0}=x)$.  \vspace{3mm}

\noindent $(\C_\zs{1})$ {\em Assume that there exists a set $C\in\Bc$ with $\mu(C)<\infty$ such that

\begin{enumerate}
 
 \item[$({\rm C}1.1)$] $f_\zs{*}=\inf_\zs{x,y\in C}\,f_\zs{1}(y|x)>0$.

\item[$({\rm C}1.2)$] There exists  $\Xc\to [1,\infty)$ Lyapunov's function $\V$ such that $\V(x)\ge \wt{g}(x)$ and $\V^{*}=\sup_\zs{x\in C} V(x)<\infty$.

\item[$({\rm C}1.3)$] For some $0<\rho<1$ and $D>0$ and for all  $x\in\Xc$, $\EV_\zs{x,0}\,\V(X_\zs{1})\le (1-\rho) \V(x)+D\Ind{C}(x)$.

\end{enumerate}
}
 
\noindent 
$(\C_\zs{2})$ 
{\em Assume that there exists $\r> 2$
such that 
$$
\check{g}_\zs{\r}=
\sup_\zs{k\ge 1}\,
\EV_\zs{\infty}\,
g^{*}_\zs{\r}(X_\zs{k})\,
<\infty
\quad\mbox{and}\quad
\check{\upsilon}_\zs{\r}=
\sup_\zs{k\ge 1}\,
\EV_\zs{\infty}\,
\upsilon^{*}_\zs{\r}(X_\zs{k})\,
<\infty
\,,
$$
where $g^{*}_\zs{\r}(x)=\sup_\zs{n\ge 1}\EV_\zs{x,0}[g(X_\zs{n},X_\zs{n-1}]^{\r}$ and $\upsilon^{*}_\zs{\r}(x)=\sup_\zs{n\ge 0}\EV_\zs{x,0}[\V(X_\zs{n})]^{\r}$.
}

\begin{theorem} \label{Th.sec:Mrk.1} 
 Conditions $(\C_\zs{1})$ and $(\C_\zs{2})$ imply condition $(\A_\zs{2})$ with $I=\lambda(\wt{g})$.
\end{theorem}

\proof
Note first that in the Markov case 
\begin{equation}\label{sec:Mrk.6-1n}
\wt{Z}_\zs{k+n}^{k}= \frac{1}{n}\,\sum^{n}_\zs{l=1} g(X_\zs{l+k},X_\zs{l+k-1})
- \lambda(\wt{g})\,.
\end{equation}
Therefore, using the fact that the process $(X_\zs{n})_\zs{n\ge \nu+1}$ is homogeneous, we obtain 
$$
\Pb_\zs{\nu} \set{|\wt{Z}_\zs{\nu+n}^\nu | \ge \vae}= \EV_\zs{\nu}\,\Psi_\zs{n}(X_\zs{\nu})= \EV_\zs{\infty}\,\Psi_\zs{n}(X_\zs{\nu})\,,
$$
where  $\Psi_\zs{n}(x)=\Pb_\zs{x,0}(|\wt{W}_\zs{n}| \ge  \vae)$. 
Note now that in view of condition $(\C_\zs{1})$, for any $x\in C$,
\begin{align*}
\Pb_\zs{x,0}(A)=\int_\zs{A}\,f_\zs{1}(y\vert x)\,\mu(\d y)\,\ge\,\int_\zs{A\cap C}\,f_\zs{1}(y\vert x)\,\mu(\d y)
\ge\, f_\zs{*}\mu(A\cap C)=\delta \varsigma(A)\,,
\end{align*}
where $\delta=f_\zs{*} \mu(C)$ and $\varsigma(A)=\mu(A\cap C)/\mu(C)$.
So Theorem~\ref{Th.AMc.1} in Appendix~\ref{B} implies condition $(\B_\zs{1}$), and therefore, Proposition~\ref{Pr.sec:Mrk.1} yields 
$$
\Psi_\zs{n}(x) \le \,W^{*}_\zs{\r}\,
\frac{\left(1+\upsilon^{*}_\zs{\r}(x)+g^{*}_\zs{\r}(x)\right)}{n^{\r/2}\vae^{r}} \quad \text{for any}~ x\in\Xc
\,,
$$
where $W^{*}_\zs{\r}$ is defined in \eqref{sec:Mrk.4}.
Thus, using condition $(\C_\zs{2}$) we obtain that
$$
\sup_\zs{\nu\ge 0}\,
\Pb_\zs{\nu} \set{|\wt{Z}_\zs{\nu+n}^\nu | \ge \vae}\,\le\,
\frac{W^{*}_\zs{\r}\left(1+\check{g}_\zs{\r}+
\check{\upsilon}_\zs{\r}\right)}{n^{\r/2}\vae^{r}}\,.
$$
This implies immediately that for any positive $\vae$ the sum defined in  \eqref{sec:MaRe.3} is bounded as 
$$
\Upsilon^{*}(\varepsilon)
\,
\le\,
\frac{W^{*}_\zs{\r}\left(1+\check{g}_\zs{\r}+
\check{\upsilon}_\zs{\r}\right)}{\vae^{r}}
\,\sum_\zs{n\ge 1}
\frac{1}{n^{\r/2}}\,.
$$
Hence Theorem~\ref{Th.sec:Mrk.1}.
\endproof

Now we obtain sufficient conditions for $(\B'_\zs{1})$ and $(\B'_\zs{2})$. To this end, we denote 
by $f_\zs{p}(y|x)$ the conditional density of $X_\zs{k+p}$ with respect to $X_\zs{k}$. \vspace{3mm}

\noindent $(\C'_\zs{1})$ {\em 
Assume that there exist an integer $p\ge 1$ and a set $C\in\Bc$ with $\mu(C)<\infty$ such that
\begin{enumerate}
 
 \item[$({\rm C^{'}}1.1)$] $f_\zs{*}=\inf_\zs{x,y\in C}\,f_\zs{p}(y|x)>0$ .
 
\item[$({\rm C^{'}}1.2)$] There exists  $\Xc\to [1,\infty)$ Lyapunov's function $\V$
such that $\V(x)\ge \wt{g}(x)$ and $\V^{*}=\sup_\zs{x\in C} V(x)<\infty$.

\item[$({\rm C^{'}}1.3)$] For some $0<\rho<1$ and $D>0$  and for all  $x\in\Xc$, $\EV_\zs{x,1}\,\V(X_\zs{p})\le (1-\rho) \V(x)+D\Ind{C}(x)$.
\end{enumerate}
}
 
\begin{theorem} \label{Th.sec:Mrk.2} 
 Conditions $(\C'_\zs{1})$ and $\C_\zs{2})$ imply condition $(\A_\zs{2})$ with $I=\lambda(\wt{g})$.
\end{theorem}

\proof
First, let $\wt{X}_\zs{l}=X_\zs{k+pl+\iota}$ for some fixed $0\le \iota\le p-1$.
Condition $(\C'_\zs{1})$ implies that 
for any $0\le \iota\le p-1$ the transition probability of the homogeneous Markov process $(\wt{X}_\zs{l})_\zs{l\ge 1}$
 for any $x\in C$
\begin{align*}
\tilde{\Pb}_\zs{x}(A)=\int_\zs{A}\,f_\zs{p}(y\vert x)\,\mu(\d y)\,\ge\,\int_\zs{A\cap C}\,f_\zs{1}(y\vert x)\,\mu(\d y)
\ge\, f_\zs{*}\mu(A\cap C)=\delta \varsigma(A)\,,
\end{align*}
where $\delta=f_\zs{*} \mu(C)$ and $\varsigma(A)=\mu(A\cap C)/\mu(C)$.
So Theorem~\ref{Th.AMc.1} implies condition $(\B'_\zs{1}$)
with the same $R_\zs{\iota}=R$, $\lambda_\zs{\iota}=\lambda$, $V_\zs{\iota}=V$ and $\kappa_\zs{\iota}$
for $0\le \iota\le p-1$.  Hence, in this case $\overline{\lambda}=\lambda$
and, therefore, for any $x\in\Xc$ by Proposition~\ref{Pr.sec:Mrk.1-1}  and
condition $(\C_\zs{2}$) for the process $Z^{k}_\zs{k+n}$ defined in \eqref{sec:Mrk.6-1n}
 we obtain that
$$
\sup_\zs{\nu\ge 0}\,
\Pb_\zs{\nu} \set{|\wt{Z}_\zs{\nu+n}^\nu | \ge \vae}\,\le\,
\frac{W^{*}_\zs{\r}\left(1+\check{g}_\zs{\r}+
\check{\upsilon}_\zs{\r}\right)}{n^{\r/2}\vae^{r}}\,.
$$
This implies immediately that for any $\vae>0$ the sum defined in \eqref{sec:MaRe.3} is bounded as 
$$
\Upsilon^{*}(\varepsilon)
\,
\le\,
\frac{W^{*}_\zs{\r}\left(1+\check{g}_\zs{\r}+
\check{\upsilon}_\zs{\r}\right)}{\vae^{r}}
\,\sum_\zs{n\ge 1}
\frac{1}{n^{\r/2}}\,.
$$
Hence Theorem~\ref{Th.sec:Mrk.2}.
\endproof

\section{Examples}\label{sec:Ex}

We now present several examples that illustrate the general theory developed in Sections~\ref{sec:Bay} and \ref{sec:MaRe}. The main goal is to verify condition $(\A_2)$ in order to be able to apply 
the theorems proved in Sections~\ref{sec:Bay} and \ref{sec:MaRe} and establish asymptotic pointwise and minimax optimality of the SR detection procedure.

\subsection{Example 1: Two dimensional AR process}\label{ssec:LaiCn}
This example motivates the necessity of relaxing conditions \eqref{sec:LaiCn.0-0} and \eqref{sec:LaiCn.0-1} proposed in \cite{LaiIEEE98} in certain interesting problems. 
It shows that both conditions \eqref{sec:LaiCn.0-0} and \eqref{sec:LaiCn.0-1} do not hold, while our uniform complete convergence condition $(\A_\zs{2})$ holds. 

Hereafter the prime in the vector $Y^{\prime}$ denotes the transposition. Consider the two dimensional autoregression (AR) process $X_\zs{k}=(X_\zs{1,k}\,,\,X_\zs{2,k})'$  defined as
\begin{equation}
\label{sec:LaiCn.1}
X_\zs{k}= \left(\Lambda\,\Ind{k\le \nu}+A_\zs{k}\right)\,X_\zs{k-1}+\xi_\zs{k}\,,
\end{equation}
where 
$$
\Lambda=
\left(
\begin{array}{ll}
 \lambda_\zs{1,}\,, &0\\[2mm]
0\,,& \lambda_\zs{2}
\end{array}
\right)
\,,
\quad
A_\zs{k}
=
\left(
\begin{array}{ll}
 \sigma_\zs{1}\,\eta_\zs{1,k}\,, &0\\[2mm]
0\,,& \sigma_\zs{2}\,\eta_\zs{2,k}
\end{array}
\right)
\quad\mbox{and}\quad
\xi_\zs{k}=
\left(
\begin{array}{c}
 \xi_\zs{1,k}\\[2mm]
 \xi_\zs{2,k}
\end{array}
\right)\,.
$$
Here 
the sequences  $(\eta_\zs{1,k})_\zs{k\ge 1}$ and $(\eta_\zs{2,k})_\zs{k\ge 1}$
are iid $\Nc(0,1)$ random variables independent of the sequence $(\xi_\zs{k})_\zs{k\ge 1}$, which is the iid sequence of 
$\Nc(0,Q)$ random vectors with
$$
Q
=
\left(
\begin{array}{cc}
1+\rho^{2} &,\,\rho\\[2mm]
\rho &,\, 1
\end{array}
\right)
$$
and $\rho>0$ is some fixed number which will be specified later.  It is clear that the iid random matrices $(A_\zs{k})_\zs{k\ge 1}$  
in \eqref{sec:LaiCn.1} are such that
$$
\EV[A_\zs{1}\,\otimes\,A_\zs{1}]=
\left(
\begin{array}{clll}
 \sigma^{2}_\zs{1}\,,&0\,,&0\,,&0\\[2mm]
0\,,&0\,,&0\,,&0 \\[2mm]
0\,,&0\,,&0\,,&0 \\[2mm]
0\,,&0\,,&0\,,&\sigma^{2}_\zs{2}
\end{array}
\right)
\,.
$$
As to the coefficients $\sigma_\zs{i}$, we choose them so that this matrix has the modules of its eigenvalues less than one, i.e., 
\begin{equation}
\label{sec:LaiCn.3}
0<\sigma^{2}_\zs{1}<1
\quad\mbox{and}\quad
0<\sigma^{2}_\zs{2}<1\,.
\end{equation}
Under these conditions the process $(X_\zs{k})_\zs{k> \nu}$ has the stationary distribution in $\bbr^{2}$
given by
\begin{equation}
\label{sec:LaiCn.4}
\zeta=
\left(
\begin{array}{l}
 \zeta_\zs{1}\\[2mm]
 \zeta_\zs{2}
\end{array}
\right)
=
\sum^{\infty}_\zs{k=1}\,
\Pi_\zs{k-1}
\,\xi_\zs{k}
\,,
\end{equation}
where $\Pi_\zs{0}=I_\zs{2}$ and $\Pi_\zs{m}=\prod^{m}_\zs{j=1}\,A_\zs{j}$ for $m\ge 1$.
One can deduce directly that this vector, conditioned on $\Gc=\sigma\{\eta_\zs{1,k},\eta_\zs{2,k}\,, k\ge 1\}$, is Gaussian $\Nc(0,\F)$
with
\begin{equation}
\label{sec:LaiCn.4-0}
\F=
\sum^{\infty}_\zs{k=1}\,
\Pi_\zs{k-1}
\,
V
\,
\Pi_\zs{k-1}
=
\left(
\begin{array}{cc}
(1+\rho^{2})\varsigma_\zs{11} &,\,\rho\varsigma_\zs{12}\\[2mm]
\rho \varsigma_\zs{12}&,\, \varsigma_\zs{22}
\end{array}
\right)
\,,
\end{equation}
where
$
\varsigma_\zs{ij}=\sum^{\infty}_\zs{k= 1}\,\sigma^{k-1}_\zs{i}\,\sigma^{k-1}_\zs{j}\,\prod^{k-1}_\zs{l=1}\,\eta_\zs{i,l}\eta_\zs{j,l}
\,.
$

Note now that,  conditioned on $X_\zs{k-1},\ldots,X_\zs{1}$, the random vector $X_\zs{k}$
for $k>\nu$  is Gaussian $\Nc(0,\D_\zs{k-1})$ with $\D_\zs{k-1}=\G(X_\zs{k-1})$, where for $x=(x_\zs{1},x_\zs{2})$
\begin{equation}
\label{sec:LaiCn.4-0n0}
\G(x)
=
\left(
\begin{array}{cc}
1+\rho^{2}+\sigma^{2}_\zs{1}\,x^{2}_\zs{1} &,\,\rho\\[2mm]
\rho &,\, 1+\sigma^{2}_\zs{2}\,x^{2}_\zs{2}
\end{array}
\right)
\,.
\end{equation}
So we can represent the LLR as
\begin{equation} \label{sec:LaiCn.05-0}
Y_\zs{j}=\log \frac{f_\zs{1,j}(X_\zs{j}|\Xb^{j-1})}{f_\zs{0,j}(X_\zs{j}|\Xb^{j-1})}
=\varkappa(X_\zs{j-1})
-\varpi_\zs{j} \,,
\end{equation}
where for any $x=(x_\zs{1},x_\zs{2})'$  
\begin{equation}\label{defkappavarpi}
\varkappa(x)=\frac{1}{2}\,x'\,\Lambda\,\G^{-1}(x)\,\Lambda x'
\quad\mbox{and}\quad
\varpi_\zs{j}=X_\zs{j-1}'\,\Lambda\,D^{-1}_\zs{j-1}\,X_\zs{j}
\,.
\end{equation}

Therefore, by the ergodic theorem, 
$$
\Pb_\zs{0}\set{\lim_\zs{n\to\infty}n^{-1}\sum_\zs{j=1}^{n} Y_\zs{j}=\EV\,\varkappa(\zeta)}=1
\,,
$$
i.e., $I=\EV\,\varkappa(\zeta)$, where the vector $\zeta$ is defined in \eqref{sec:LaiCn.4}. 
Clearly,  condition \eqref{sec:LaiCn.0-1} in this case has the following form: for any $\varepsilon>0$
$$
\lim_\zs{n\to\infty}\,
\sup_\zs{x\in\bbr^{2}}\,
\Pb_\zs{0}\left(
\sum_\zs{j=1}^{n}\, Y_\zs{j}<(I-\varepsilon)n
\,\vert\,X_\zs{0}=x
\right)\,=\,0\,,
$$
where $I=\EV\,\varkappa(\zeta)$.
We now establish that it does not hold by showing that for some  $0<\varepsilon<1$
\begin{equation}
\label{sec:LaiCn.5liminf}
\liminf_\zs{n\to\infty}\,
\sup_\zs{x\in\bbr^{2}}\,
\Pb_\zs{0}\left(
\sum_\zs{j=1}^{n}\, Y_\zs{j}<(I-\varepsilon)n
\,\vert\,X_\zs{0}=x
\right)\,>\,0\,.
\end{equation}
Observe first that for any $n\ge 1$
$$
\sup_\zs{x\in\bbr^{2}}\,
\Pb_\zs{0}\left(
\sum_\zs{j=1}^{n}\, Y_\zs{j}<(I-\varepsilon)n
\,\vert\,X_\zs{0}=x
\right)
\ge 
\lim_\zs{x_\zs{2}\to\infty}\,
\Pb_\zs{0}\left(
\sum_\zs{j=1}^{n}\, Y_\zs{j}<(I-\varepsilon)n
\,\vert\,X_\zs{0}=(0,x_\zs{2})' \right)
$$
and that for any $x_\zs{1}$
$$
\varkappa_\zs{1}(x_\zs{1})=\lim_\zs{\vert x_\zs{2}\vert\to\infty}\,\varkappa(x_\zs{1},x_\zs{2})
=\frac{\lambda^{2}_\zs{1}x^{2}_\zs{1}}{2(1+\rho^{2}+\sigma^{2}_\zs{1}x^{2}_\zs{1})}\,
+
\frac{\lambda^{2}_\zs{2}}{2\sigma^{2}_\zs{2}}
$$
and $\Pb_\zs{0}(\lim_{n\to\infty} n^{-1} \sum_{j=1}^{n}\varkappa_{1}(X_{1,j})=I_{1})=1$,
where 
$$
I_\zs{1}=\EV\,\varkappa_\zs{1}(\zeta_\zs{1})
\le 
\frac{\lambda^{2}_\zs{1}\,\EV\,\zeta^{2}_\zs{1}}{2(1+\rho^{2})}\,
+
\frac{\lambda^{2}_\zs{2}}{2\sigma^{2}_\zs{2}}
=
\frac{\lambda^{2}_\zs{1}\,\EV\,\varsigma_\zs{11}}{2}\,
+
\frac{\lambda^{2}_\zs{2}}{2\sigma^{2}_\zs{2}}
=
\frac{\lambda^{2}_\zs{1}}{2(1-\sigma^{2}_\zs{1})}\,
+
\frac{\lambda^{2}_\zs{2}}{2\sigma^{2}_\zs{2}}
\,.
$$

Let us show now  that there exist $\sigma_\zs{1}>0$ and $\rho>0$ for which $I>I_\zs{1}$. If so, then taking into account Lemma \ref{Le.sec:A_Ex_0} in Appendix~\ref{A}, we obtain
that for come $\varepsilon>0$
$$
\lim_\zs{n\to\infty}
\sup_\zs{x\in\bbr^{2}}\,
\Pb_\zs{0}\left(
\sum_\zs{j=1}^{n}\, Y_\zs{j}<(I-\varepsilon)n
\,\vert\,X_\zs{0}=x
\right)=1 \, .
$$
Hence, \eqref{sec:LaiCn.5liminf} follows.
Indeed, choosing in \eqref{sec:LaiCn.1} the parameter $\sigma^{2}_\zs{1}$ as a function of $\rho$ such that 
$\sigma^{2}\rho^{4}\to 0$ as $\rho\to\infty$, we obtain in view of Lemma  \ref{Le.sec:A_Ex_1} in  Appendix~\ref{A}
that  there exist $\sigma_\zs{1}$ and $\rho>0$ for which $I>I_\zs{1}$. This implies the inequality \eqref{sec:LaiCn.5liminf}, and hence, condition \eqref{sec:LaiCn.0-1} does not hold.  

Note that condition \eqref{sec:LaiCn.0-0} also does not hold. Indeed, this condition  in the case considered has the following form: for any $\varepsilon>0$
$$
\lim_\zs{n\to\infty}\,
\sup_\zs{x\in\bbr^{2}}\,
\Pb_\zs{0}\left(
\sum_\zs{j=1}^{n}\, Y_\zs{j}>(I+\varepsilon)n
\,\vert\,X_\zs{0}=x
\right)\,=\,0\,.
$$
If we put $\rho=0$  then we obtain that
$
\lim_\zs{\vert x_\zs{1}\vert, \vert x_\zs{2}\vert\to\infty}\,\varkappa(x_\zs{1},x_\zs{2})
=\lambda^{2}_\zs{1}/(2\sigma^{2}_\zs{1})\,
+
\lambda^{2}_\zs{2}/(2\sigma^{2}_\zs{2})
:=\kappa^{*}$
and $\kappa^{*}>\varkappa(x_\zs{1},x_\zs{2})$ for any $x_\zs{1}$ and $x_\zs{2}$ from $\bbr$. Therefore, 
$\kappa^{*}>I=\EV \varkappa(\zeta)$. Similarly to the above reasoning we obtain that for some $\varepsilon>0$
\begin{align*}
\lim_\zs{n\to\infty}\,\sup_\zs{x\in\bbr^{2}}\,\Pb_\zs{0}\left(\sum_\zs{j=1}^{n}\, Y_\zs{j}>(I+\varepsilon)n\,\vert\,X_\zs{0}=x\right)
\ge\,
\lim_\zs{n\to\infty}\,
\lim_\zs{\vert x_\zs{1}\vert 	\wedge \vert x_\zs{2}\vert\to\infty}
\Pb_\zs{0}\left(
\sum_\zs{j=1}^{n}\, Y_\zs{j}>(I+\varepsilon)n
\,\vert\,X_\zs{0}=x \right)\,=1\,,
\end{align*}
where $a\wedge b=\min(a,b)$.

On the other hand, our uniform complete convergence condition $(\A_2)$ holds. Indeed, as we will see in  Example 4 below, condition  $(\A_2)$ holds even for a more general vector AR model 
than \eqref{sec:LaiCn.1}.  Thus, the SR procedure is asymptotically minimax.

\subsection{Example 2:  Change in the correlation coefficient of the AR(1) model}\label{ssec:AR1cor}

Consider the change of the correlation coefficient in the first-order AR model
\begin{equation}\label{sec:Ex.1-0-0}
X_\zs{n} =\vartheta_\zs{n}\,X_\zs{n-1}+w_\zs{n}\,,
\end{equation}
where $\vartheta_\zs{n}=\a_\zs{0}\Ind{n \le \nu}+\a_\zs{1}\Ind{n > \nu}$
and  $(w_\zs{n})_\zs{n\ge 1}$ are iid  not necessarily Gaussian random variables 
with $\EV\,w_\zs{1}=0$, $\EV\,w^{2}_\zs{1}=1$ and a known density $\psi(x)$
such that for any $n\ge 1$
\begin{equation}\label{sec:Ex.1-1}
\inf_\zs{-n\le x\le n}\,\psi(x)>0 \,.
\end{equation}
We assume that the parameters $-1<\a_\zs{i}<1$ are known.
In this case the ergodic distributions for $(X_\zs{n})_\zs{n\le \nu}$ and $(X_\zs{n})_\zs{n\ge \nu+1}$ are  given 
by the random variables $\wt{w}_\zs{0}$ and $\wt{w}_\zs{1}$, respectively, which are defined as
\begin{equation}\label{sec:Ex.1-1n}
\wt{w}_\zs{i}=\sum^{\infty}_\zs{j=1}\,(\a_\zs{i})^{j-1}\,w_\zs{j}, \quad i=0,1\,.
\end{equation}
The pre-change and post-change conditional densities are $f_{0}(X_\zs{n}|X_\zs{n-1})= \psi(X_n-\a_\zs{0} X_{n-1})$ for all $1\le n \le \nu$ and $f_{1}(X_\zs{n}|X_\zs{n-1})=\psi(X_n-\a_\zs{1} X_{n-1})$ 
for $n > \nu$, where $X_\zs{0}$ is an initial value independent of the sequence $(w_\zs{n})_\zs{n\ge 1}$. Note that condition
\eqref{sec:Ex.1-1} implies the lower bound (C$1.1$) in condition $(\C_\zs{1})$ for any ``minorization'' set of the form $C=[-n,n]$. 
It is easily seen that 
\begin{equation}\label{sec:Ex.1-2}
g(y,x)=\log\frac{f_\zs{1}(y|x)}{f_\zs{0}(y|x)}=
\log\frac{\psi(y-\a_\zs{1}x)}{\psi(y-\a_\zs{0}x)}
\end{equation}
and 
\begin{equation}\label{sec:Ex.1-3}
\wt{g}(x)=\int_\zs{\bbr}\,\log\frac{\psi(y-\a_\zs{1}x)}{\psi(y-\a_\zs{0}x)}\,
\psi(y-\a_\zs{1}x)\,\d y\,.
\end{equation}
Assume that there exist $q^{*}\ge 1$ and $\iota>0$ such that
\begin{equation}\label{sec:Ex.2-1}
\EV\,|w_\zs{1}|^{\iota}<\infty\,,
\quad
\sup_\zs{y,x\in\bbr}\,
\frac{\vert g(y,x)\vert}{
\left(
1+|y|^{\iota}+|x|^{\iota}
\right)}
\,
 \le\,q^{*}
\quad\mbox{and}\quad
\sup_\zs{x\in\bbr}\,
\frac{\wt{g}(x)}{
\left(
1+|x|^{\iota}
\right)}
\,
 \le\,q^{*}
\,. 
\end{equation}
For example, in the Gaussian case (i.e., $\psi$ is $(0,1)$ Gaussian density), 
$$
g(y,x)=\frac{(y-\a_\zs{0}x)^{2}-(y-\a_\zs{1}x)^{2}}{2}
\quad\mbox{and}\quad
\wt{g}(x)=\frac{(\a_\zs{1}-\a_\zs{0})^{2}\,x^{2}}{2}\,,
$$
 i.e., conditions \eqref{sec:Ex.2-1} are satisfied with $\iota=2$ and
$$
q^{*}=
\max\left\{1,\frac{\vert\a^{2}_\zs{1}-\a^{2}_\zs{0}\vert+(\a_\zs{1}-\a_\zs{0})^{2}+1}{2}\right\} \,.
$$

Define the Lyapunov function as
\begin{equation}\label{sec:Ex.2-2}
\V(x)=q^{*}\left(
1+|x|^{\iota}
\right)\,.
\end{equation}
Obviously, 
$$
\lim_\zs{|x|\to\infty}\,\frac{\EV_\zs{x,0}\,\V(X_\zs{1})}{\V(x)}
=\lim_\zs{|x|\to\infty}\,\frac{1+\EV\,|\a_\zs{1} X+w_\zs{1}|^{\iota}}{1+|x|^{\iota}}
=|\a_\zs{1}|^{\iota}<1\,.
$$
Therefore, for any $|\a_\zs{1}|^{\iota}<\rho<1$ there exist $n\ge 1$ and
$D>0$ such that condition $(\C_\zs{1}$) holds with $C=[-n,n]$. 

Let us check now condition ($\C_\zs{2}$). Assume that there exists  $\r> 2$ for which
\begin{equation}\label{sec:Ex.2-3}
\int_\zs{\bbr}\,|v|^{\r_\zs{1}}\,\psi(v)\,\d v\,<\,\infty,\quad \r_\zs{1}=\iota\r\,.
\end{equation}
This condition implies that 
$\EV\,\vert\wt{w}_\zs{0}\vert^{\r_\zs{1}}<\infty$
and 
$\EV\,\vert\wt{w}_\zs{1}\vert^{\r_\zs{1}}<\infty$.
Moreover, taking into account the ergodicity properties, we obtain that for any $x\in\bbr$ 
\begin{equation}\label{sec:Ex.2-3n}
\lim_\zs{k\to\infty}\,\EV_\zs{x,\infty}\,\vert X_\zs{k}\vert^{\r_\zs{1}}=\EV\,\vert\wt{w}_\zs{0}\vert^{\r_\zs{1}}<\infty
\quad\mbox{and}\quad
\lim_\zs{k\to\infty}\,\EV_\zs{x,0}\,\vert X_\zs{k}\vert^{\r_\zs{1}}=\EV\,\vert\wt{w}_\zs{0}\vert^{\r_\zs{1}}<\infty\,.
\end{equation}
Note also that under the probability $\Pb_\zs{x,0}$ for any $j\ge 1$, $X_\zs{j}=\a^{j}_\zs{1}\,x+\sum^{j}_\zs{l=1}\a^{j-l}_\zs{1} w_\zs{l}$.
Therefore, $\EV_{x,0}\vert X_{j}\vert^{\r_{1}}\le 2^{\r_{1}}(\vert x\vert^{\r_{1}}+\EV_{0,0}\vert X_{j}\vert^{\r_{1}})$,
i.e., using the last convergence in \eqref{sec:Ex.2-3n}, we obtain that for some $C^{*}>0$
$$
M^{*}(x)=\sup_\zs{j\ge 1}\,\EV_\zs{x,0}\,\vert X_\zs{j}\vert^{\r_\zs{1}}\,
\le C^{*} (1+\vert x\vert^{\r_\zs{1}})\,.
$$
Using now the first convergence in \eqref{sec:Ex.2-3n}, we obtain that $\sup_\zs{k\ge 1}\EV\,M^{*}(X_\zs{k})<\infty$.
So the upper bounds in \eqref{sec:Ex.2-1} imply condition ($\C_\zs{2}$). 
By  Theorem~\ref{Th.sec:Mrk.1}, condition $(\A_\zs{2}$)
holds for the model \eqref{sec:Ex.1-0-0} if density $\psi$
of the iid  random variables $(w_\zs{n})_\zs{n\ge 1}$ satisfies conditions  
\eqref{sec:Ex.1-1} and  \eqref{sec:Ex.2-1}. The Kullback--Leibler information number is
$$
I=\lambda(\wt{g})= \int_\zs{\bbr}\,\left(
 \int_\zs{\bbr} \log \frac{\psi(y-\a_\zs{1} x)}{\psi(y-\a_\zs{0} x)}\,\psi(y-\a_\zs{1} x)\d y\right)\, \lambda(\d x)\,,
$$
where $\lambda$ is the  distribution of  $\wt{w}_\zs{1}$ given  in \eqref{sec:Ex.1-1n}.

Hence, by Theorem~\ref{Th.sec:Bay.3} and Theorem~\ref{Th.sec:Cnrsk.2}, the SR procedure is asymptotically minimax with respect to the expected detection delays.
 
 In the particular Gaussian case  where $\psi$ is $\Nc(0,1)$,
 the random variable $\wt{w}$ is $\Nc\left(0,(1-\a_\zs{1})^{-2}\right)$,  and the Kullback--Leibler information number can be calculated explicitly, 
 $I=(\a_\zs{1}-\a_\zs{0})^{2}/2(1-\a_\zs{1}^2)$.


\subsection{Example 3: AR process with ARCH(1) errors}

Consider now the change of the correlation coefficient in the first-order AR 
model with ARCH(1) errors \cite{BoKl01},  assuming that for $n\ge 1$
\begin{equation}\label{sec:Ex.4}
X_\zs{n} =\vartheta_\zs{n}\,X_\zs{n-1}+\left(1+\sigma^{2} X^{2}_{n-1}\right)^{1/2}\,
w_\zs{n}\,,
\end{equation}
where an initial value $X_\zs{0}$ is independent of the sequence $(w_\zs{n})_\zs{n\ge 1}$. The sequence $\vartheta_\zs{n}$ is defined in \eqref{sec:Ex.1-0-0}
with the known parameters $\a_\zs{i}$ such that $\a^{2}_\zs{i}+\sigma^{2}<1$. As in the model  \eqref{sec:Ex.1-0-0}, we assume that 
$(w_\zs{n})_\zs{n\ge 1}$ are iid  not necessarily Gaussian random variables with $\EV\,w_\zs{1}=0$, $\EV\,w^{2}_\zs{1}=1$ and
a known density $\psi(x)$ satisfying condition \eqref{sec:Ex.1-1}. The variance $\sigma^2>0$ is  known.
In just the same way as in the model \eqref{sec:Ex.1-0-0}, we find that the pre-change and post-change conditional densities $f_{0}(X_\zs{n}|X_\zs{n-1})$ and
$f_{1}(X_\zs{n}|X_\zs{n-1})$ are of the form
$$
f_\zs{0}(y\mid x)=\frac{\psi\left(\l_\zs{0}(y,x) \right)}{\sqrt{1+\sigma^{2}x^{2}}} \quad \text{and} \quad f_\zs{1}(y\mid x)=\frac{\psi\left(\l_\zs{1}(y,x) \right)}{\sqrt{1+\sigma^{2}x^{2}}}\,,
$$
where $\l_\zs{0}(y,x)=(y-\a_\zs{0} x)/\sqrt{1+\sigma^{2}x^{2}}$ and $\l_\zs{1}(y,x)=(y-\a_\zs{1} x)/\sqrt{1+\sigma^{2}x^{2}}$.

Obviously, the property \eqref{sec:Ex.1-1} implies the lower bound (C$1.1$) in condition $(\C_\zs{1}$).

The function \eqref{sec:Ex.1-2} is given by $g(y,x)=\log[\psi(\l_\zs{1}(x,y))/\psi(\l_\zs{0}(x,y))]$ and
$$
\wt{g}(x)=\,
\frac{1}{\sqrt{1+\sigma^{2}x^{2}}}
\int_\zs{\bbr}\,
\left(
\log\frac{\psi\left(\l_\zs{1}(y,x)\right)}{\psi\left(\l_\zs{0}(y,x)\right)}
\right)
\,\psi\left(\l_\zs{1}(y,x) \right)\,
\d y
\,.
$$

Assume that there exist $q^{*}\ge 1$ and $\iota>0$ such that
\begin{equation}\label{sec:Ex.4-1}
\sup_\zs{y,x\in\bbr}\,
\frac{\vert g(y,x)\vert}{
\left(
1+|\l_\zs{0}(y,x)|^{\iota}+|\l_\zs{1}(y,x)|^{\iota}
\right)}
\,
 \le\,q^{*}
\quad\mbox{and}\quad
\sup_\zs{x\in\bbr}\,
\wt{g}(x)
\,
 \le\,q^{*}
\,. 
\end{equation}
For example, in the Gaussian case (i.e., $\psi$ is standard Gaussian density), 
$$
g(y,x)=\frac{\l^{2}_\zs{0}(y,x)
-
\l^{2}_\zs{1}(y,x)}{2}
\quad\mbox{and}\quad
\wt{g}(x)=\frac{(\a_\zs{1}-\a_\zs{0})^{2}\,x^{2}}
{2(1+\sigma^{2}x^{2})}\,,
$$
 i.e., conditions \eqref{sec:Ex.4-1} are satisfied with $\iota=2$.

The Lyapunov function is any $\bbr\to (1,+\infty)$ function  which satisfies the drift condition (C$1.3$).
We set $\V(x)=q^{*}(1+|x|^{\delta})$ for $0<\delta<\x_\zs{*}$,
where $\x_\zs{*}=\min(\x_\zs{*,0}\,,\,\x_\zs{*,1})$ and $\x_\zs{*,i}$ is a unique positive root of the equation 
$\check{\kappa}(\x)=1$, where $\check{\kappa}(\x)=\EV\,|\a_\zs{i}+\sigma\, w_\zs{1}|^{\x}$.

It is well known \cite{KlPe04} that if $\EV\,w^{2}_\zs{1}=1$, then $\x_\zs{*}>2$. 
Direct calculations yield 
\begin{equation}\label{sec:Ex.4-2n}
\lim_\zs{|x|\to\infty}\,\frac{\EV_\zs{x,0}\V(X_\zs{1})}{\V(x)}
=\check{\kappa}(\delta)<1\,.
\end{equation}
Therefore, for any $\check{\kappa}(\delta)<\rho<1$ there exist $n\ge 1$ and
$D>0$ for which condition $(\C_\zs{1}$) holds with $C=[-n,n]$.

Next, we verify condition ($\C_\zs{2}$). To this end, note that under the probability $\Pb_0$ we have
$$
\l_\zs{1}(X_\zs{j},X_\zs{j-1})=w_\zs{j}
\quad\mbox{and}\quad
|\l_\zs{0}(X_\zs{j},X_\zs{j-1})|\le |w_\zs{j}|+\frac{| \a_\zs{1}-\a_\zs{0} |}{\sigma} 
\,.
$$
So, for any $\r> 2$ satisfying \eqref{sec:Ex.2-3} with $\iota>0$ from  condition \eqref{sec:Ex.4-1}, we obtain  that, 
for some constant $C^{*}>0$, $g^{*}_\zs{\r}(x)\le C^{*}(1+ \EV\vert w_\zs{1}\vert^{\iota\r})$,
i.e., $\check{g}_\zs{\r}<\infty$. Now we check the last inequality in ($\C_\zs{2}$). Fix $\r>2$ such 
that $\r_\zs{1}=\delta\r<\x_\zs{*}$. Evidently, this is possible for sufficiently small $\delta>0$. Similarly to \eqref{sec:Ex.4-2n} we can obtain that
$$
\lim_\zs{|x|\to\infty}\,\frac{\EV_\zs{x,0}\V_\zs{1}(X_\zs{1})}{\V_\zs{1}(x)}
=\check{\kappa}(\r_\zs{1})<1\,,
$$
where $V_\zs{1}(x)=1+\vert x\vert^{\r_\zs{1}}$\,. Therefore, conditions $(\H_\zs{1})$ and $(\H_\zs{2})$ hold, and 
using Theorem \ref{Th.AMc.1} in Appendix~\ref{B}, we obtain that for some constant $C^{*}>0$, $\sup_\zs{j\ge 1}\EV_\zs{x,0}\vert X_\zs{j}\vert^{\r_\zs{1}}\le C^{*}(1+\vert x\vert^{\r_\zs{1}})$.
Similarly we obtain that $\sup_\zs{j\ge 1}\EV_\zs{\infty}\vert X_\zs{j}\vert^{\r_\zs{1}}<\infty$, i.e., ($\C_\zs{2}$) is satisfied.

Thus, by Theorem~\ref{Th.sec:Mrk.1}, condition $(\A_\zs{2}$) holds for the model \eqref{sec:Ex.4} where  
the iid  random variables $(w_\zs{n})_\zs{n\ge 1}$ have density $\psi(x)$ that satisfies conditions  \eqref{sec:Ex.1-1} 
and \eqref{sec:Ex.2-3} with $\iota>0$ from  condition \eqref{sec:Ex.4-1}.
Note that in this case  there exists
 the stationary distribution $\lambda$ for
$(X_\zs{n})_\zs{n>\nu}$ which in the Gaussian case, 
$w_\zs{n}\sim \Nc(0,1)$, is given by the following random variable
\begin{equation}\label{sec:Ex.5}
\sum^{\infty}_\zs{j=1}\,\prod^{j-1}_\zs{l=1}\upsilon_\zs{l}\,w_\zs{j}\,,
\end{equation}
where $(\upsilon_\zs{l})_\zs{l\ge 1}$ is an iid $\Nc(\a_\zs{1},\sigma^{2})$ sequence independent of $(w_\zs{j})_\zs{j\ge 1}$.
The Kullback--Leibler information is
\begin{equation}\label{sec:Ex.6}
I= 
\int_\zs{\bbr}\,
\left(
\frac{1}{\sqrt{1+\sigma^{2}x^{2}}}
\int_\zs{\bbr}\,
\left(
\log\frac{\psi\left(\l_\zs{1}(y,x)\right)}{\psi\left(\l_\zs{0}(y,x)\right)}
\right)
\,\psi\left(\l_\zs{1}(y,x) \right)\,
\d y
\right)
\,
\lambda(\d x)
=\frac{(\a_\zs{1}-\a_\zs{0})^{2}}{\sqrt{2\pi}}\,\EV\,G(\wt{\upsilon})
\,,
\end{equation}
where
$$
 G(z)=\frac{1}{z}\,\int^{\infty}_\zs{0}\,
\frac{y^{2}\,e^{-\frac{y^{2}}{2z^{2}}}}{1+\sigma^{2}\,y^{2}}\,\d y
\quad\mbox{and}\quad
\wt{\upsilon}=
\left(1+\sum^{\infty}_\zs{j=2}\,\prod^{j-1}_\zs{l=1}\upsilon^{2}_\zs{l}\right)^{1/2}\,.
$$

By Theorem~\ref{Th.sec:Bay.3} and Theorem~\ref{Th.sec:Cnrsk.2}, the SR procedure is asymptotically minimax.

\subsection{Example 4:  Change in the parameters of the multivariate linear difference equation}

Consider the multivariate model in $\bbr^{p}$ given by
\begin{equation}\label{sec:Ex.6-0n}
X_\zs{n} = \left(A_\zs{0,n}\Ind{n\le \nu}+A_\zs{1,n}\Ind{n>\nu}\right)\,X_\zs{n-1}+w_\zs{n}\,,
\end{equation}
where $A_\zs{0,n}$ and $A_\zs{1,n}$ are $p\times p$ random matrixes and $(w_\zs{n})_\zs{n\ge 1}$ is an iid sequence of Gaussian random vectors $\Nc(0, Q_\zs{0})$ in
$\bbr^{p}$ with the positive definite $p\times p$ matrix $Q_\zs{0}$. 
Assume also that $A_\zs{i,n}=A_\zs{i}+B_\zs{n}$ and $(B_\zs{n})_\zs{n\ge 1}$ are iid Gaussian random matrixes $\Nc(0\,,Q_\zs{1})$,
where the $p^{2}\times p^{2}$ matrix $Q_\zs{1}$ is not necessary positive definite. 
Assume, in addition, that $\EV [A_\zs{i,1}\,\otimes\,A_\zs{i,1}]$, $i=0,1$ have the modules less than one.

In this case,  the processes  $(X_\zs{n})_\zs{1\le n\le \nu}$ and $(X_\zs{n})_\zs{n> \nu}$
are ergodic with the ergodic distributions give by the vectors \cite{KlPe04}
$\varsigma_\zs{i}=\sum_\zs{l\ge 1}\prod^{l-1}_\zs{j=1}A_\zs{i,j}\,w_\zs{l}$,
i.e.,  the invariant measures $\lambda_\zs{i}$ on $\bbr^{p}$ are 
defined as $\lambda_\zs{i}(A)=\Pb(\varsigma_\zs{i}\in \Gamma)$ for  any $\Gamma\in\Bc(\bbr^{p})$.
As shown in \cite{FeiginTweedieJTSA85}, there exists a positive definite $p\times p$ matrix $T$ 
and the constant $K_\zs{*}>0$ such that the function $V(x)=c (1+ x'Tx)$ and the set $C=\{x\in\bbr^{p}\,:\, x'Tx\le K_\zs{*}\}$
satisfy condition (C$1.3$) for any $c\ge 1$. The function $g(y,x)$ can be calculated for any $x,y$ from $\bbr^{p}$ as
$$
g(y,x)
=
\frac{
\vert \l_\zs{0}(y,x)\vert^{2}
-
\vert \l_\zs{1}(y,x)\vert^{2}}{2}
=y'G^{-1}(x)(A_\zs{1}-A_\zs{0})x
+
\frac{x'A'_\zs{0}G^{-1}(x)\,A_\zs{0}x-x'A'_\zs{1}\,G^{-1}(x)\,A_\zs{1}\,x}{2}
\,,
$$
where $\l_\zs{i}(y,x)=G^{-1/2}(x)(y-A_\zs{i}x) $ and $G(x)=\E\,B_\zs{1}xx'B'_\zs{1}+Q_\zs{0}$.
From this we obtain that
$$
\wt{g}(x)=
\frac{1}{2}\,\vert G^{-1/2}(x)(A_\zs{1}-A_\zs{0})x\vert
=
\frac{1}{2}\,x'(A_\zs{1}-A_\zs{0})'G^{-1}(x)(A_\zs{1}-A_\zs{0})x\,.
$$

Assume that 
\begin{equation}\label{sec:Ex.6-0n.1}
\sup_\zs{x\in\bbr^{p}}\,\vert G^{-1/2}(x)(A_\zs{1}-A_\zs{0})x\vert<\infty .
\end{equation}
Note that for the model \eqref{sec:LaiCn.1} this condition holds. So under this condition $g^{*}=\sup_\zs{x\in\bbr^{p}}\,\wt{g}(x)<\infty$.
Thus, choosing $V(x)=c^{*}\,(1+(x'Tx)^{\delta})$ with $c^{*}=1+g^{*}$ and any fixed $0<\delta<1$ and using the Jensen inequality
yields condition ($\C_\zs{1}$).  

Let us check now condition ($\C_\zs{2}$). Note that 
under the probability $\Pb_\zs{0}$ we obtain that for any $j\ge 1$ the vector
$\xi_\zs{j}=\l_\zs{1}(X_\zs{j},X_\zs{j-1})$ is $(0,I_\zs{p})$ Gaussian in $\bbr^{p}$. Moreover, in view of condition \eqref{sec:Ex.6-0n.1}
 $$
 \vert \l_\zs{0}(X_\zs{j},X_\zs{j-1})\vert =\vert \xi_\zs{j}+ G^{-1/2}(X_\zs{j-1})(A_\zs{1}-A_\zs{0})X_\zs{j-1}\vert
 \le \vert \xi_\zs{j}\vert+C^{*}
 $$
 for some positive $C^{*}$.  Clearly, $\check{g}_\zs{\r}<\infty$ for any $\r>0$.  We now check the last inequality in ($\C_\zs{2}$).
 First note that, as it is shown in  \cite{FeiginTweedieJTSA85}, under our conditions
 $\EV\vert \varsigma_\zs{i}\vert^{2}<\infty$. Next, observe that under the probability $\Pb_\zs{x,0}$
 $$
 X_\zs{j}=\prod^{l-1}_\zs{j=1}A_\zs{1,j}\,x
 +
 \sum^{j}_\zs{l=1}\,\prod^{l}_\zs{i=l+1}\,
A_\zs{1,i}\,w_\zs{l}\,.
 $$
So, for any $0<q\le 2$, $\EV_\zs{x,0}\vert X_\zs{j}\vert^{q}\le C^{*}(\vert x\vert^{q}+ \EV_\zs{0,0}\vert X_\zs{j}\vert^{q})$.
In view of the ergodicity property we obtain that
$$
\lim_\zs{j\to\infty}\,\EV_\zs{\infty}\vert X_\zs{j}\vert^{q}=\E_\zs{\infty}\vert \varsigma_\zs{0}\vert^{q}<\infty
\quad \text{and}\quad
\lim_\zs{j\to\infty}\,\EV_\zs{0,0}\vert X_\zs{j}\vert^{q}=\EV_\zs{0,0}\vert \varsigma_\zs{1}\vert^{q}<\infty\,,
$$
i.e., $\sup_\zs{j\ge 1}\,\EV_\zs{x,0}\vert X_\zs{j}\vert^{q}\le C^{*}(1+\vert x\vert^{q})$ for some positive $C^{*}$.
So  $\check{\upsilon}_\zs{\r}<\infty$ for any $\r>2$ for which $\delta r\le 2$.

Hence, by Theorem~\ref{Th.sec:Mrk.1}, condition $(\A_2)$ is satisfied 
with $I= \EV\,\wt{g}(\varsigma_\zs{1})$, and by Theorem~\ref{Th.sec:Bay.3} and Theorem~\ref{Th.sec:Cnrsk.2} the SR detection procedure is asymptotically minimax.

\subsection{Example 5:  Change in the correlation coefficients of the AR(p) model}

Let us now generalize the results of Subsection~\ref{ssec:AR1cor} for the problem of detecting the change of the correlation coefficient in the $p$-th order AR process,  assuming that for $n\ge 1$
\begin{equation}\label{sec:Ex.7}
X_\zs{n} =\vartheta_\zs{1,n}\,X_\zs{n-1}+\ldots+\vartheta_\zs{p,n}\,X_\zs{n-p}+w_\zs{n}\,,
\end{equation}
where $\vartheta_\zs{i,n}=\a_\zs{0,i}\Ind{n \le \nu}+\a_\zs{1,i}\Ind{n > \nu}$
and  $(w_\zs{n})_\zs{n\ge 1}$ are iid,  not necessarily Gaussian random variables 
with $\EV\,w_\zs{1}=0$, $\EV\,w^{2}_\zs{1}=1$.  
In the sequel we use the notation  $\a_\zs{i}=(\a_\zs{i,1},\ldots,\a_\zs{i,p})'$.
The process\eqref{sec:Ex.7} is not Markov, but the $p$-dimensional process
\begin{equation}\label{sec:Ex.7-1n}
\check{X}_\zs{n}=
(X_\zs{n},\ldots,X_\zs{n-p+1})'\in\bbr^{p}
\end{equation}
is Markov.  Note that for $n > \nu$
\begin{equation}\label{sec:Ex.8}
\check{X}_\zs{n}=A\check{X}_\zs{n-1}+\check{w}_\zs{n}\,,
\end{equation}
where
$$
A=
\left(
\begin{array}{rc}
 \a_\zs{1,1}&,\ldots,\a_\zs{1,p}\\[2mm]
1&,\ldots,0\\[2mm]
..&\ldots \\[2mm]
0&,\ldots,1,0
\end{array}
\right)
\quad\mbox{and}\quad
\check{w}_\zs{n}=(w_\zs{n},0,\ldots,0)'\in\bbr^{p}\,.
$$
It is clear that
$$
\EV[\check{w}_\zs{n}\,\check{w}'_\zs{n}]=\,B\,=
\left(
\begin{array}{rc}
1&,\ldots,0\\[2mm]
..&\ldots \\[2mm]
0&,\ldots,0
\end{array}
\right)
\,.
$$
Assume that all eigenvalues of the matrix $A$ have the modules less than one.
The ergodic distribution is given by the vector
$\varsigma=\sum_\zs{l\ge 1}A^{l-1}\,\check{w}_\zs{l}\sim \Nc(0,\F)$, where 
\begin{equation}\label{erg_gaus.00}
\F=\sum_\zs{l\ge 0}A^{l}\,B\,(A')^{l}
\,.
\end{equation}
Obviously, condition (C$1.1$) does not hold for the process \eqref{sec:Ex.7-1n}. 

To fulfill this condition
we replace this process by the embedded homogeneous Markov process $Y_\zs{n}=\check{X}_\zs{np+\iota}$ for some $0\le \iota\le p-1$.  
This process can be represented as 
\begin{equation}\label{sec:A.7}
Z_\zs{n}=A^{p}Z_\zs{n-1}+\zeta_\zs{n}\,,
\quad
\zeta_\zs{n}=\sum^{p-1}_\zs{j=0}\,A^{j}\,\check{w}_\zs{np+\iota-j}\,.
\end{equation}
Clearly, $\zeta_\zs{n}$ is Gaussian with the parameters $(0,Q)$, where $Q=\sum^{p-1}_\zs{j=0}\,A^{j}\, B\,(A')^{j}$.
One can check directly that this matrix is positive definite. Define the function $V(x): \bbr^{p}\to\bbr$ as
\begin{equation}\label{sec:A.8}
V(x)= c (1+x'Tx) \,,\quad T=\sum^{\infty}_\zs{l=0}\,(A')^{pl}\,A^{pl}\,,
\end{equation}
where $c\ge 1$ will be specified later. Let $t_\zs{\max}=\max_\zs{\|x\|=1}\,x'Tx$ and $t_\zs{*}=1-1/t_\zs{\max}$. 
Obviously, $t_\zs{\max}>1$, i.e., $0<t_\zs{*}<1$. Now we  set $K=[(1+\EV\,\|\zeta_\zs{1}\|^{2})/\rho]^{1/2}$ with $\rho=(1-t_\zs{*})/2$ and 
$D=1+\|T^{1/2}A^{p}\|^{2}\,K^{2}+\EV \|\zeta_\zs{1}\|^{2}$.
Next we need the minorizing measure in condition $(\H_\zs{1}$) on the Borel $\sigma$-field in $\bbr^{p}$. 
To this end, we define $\check{\nu}(\Gamma) =\mes(\Gamma\cap C)/\mes(C)$ for any Borel set $\Gamma$ in $\bbr^{p}$,
where $\mes(\cdot)$ is the Lebesgue measure in $\bbr^{p}$.

Finally, we show that, for any $0\le \iota<p$, the Markov process \eqref{sec:A.7} satisfies condition (${\rm C}^{'}1.3$). 
Indeed, note that 
\[
\begin{aligned}
\EV(V(Z_\zs{1})|Z_\zs{0}=x)&=c+c \EV\left\|T^{1/2}(A^{p}x+\zeta_\zs{1})\right\|^{2}\\
&= c+c \left(x'(A^{p})'T\,A^{p}x\right) +r\EV\, \left\|\zeta_\zs{1}\right\|^{2}\,.
\end{aligned}
\]
Taking into account that
\[
\frac{x'(A^{p})'TA^{p}x}{x'Tx}=1-\frac{\|x\|^{2}}{x'Tx}\le 1-\frac{1}{t_\zs{max}}=t_\zs{*}\,,
\]
we obtain that, for $\|x\|\ge K$, $\EV(V(Z_\zs{1})\,|\,Z_\zs{0}=x)\le (1-\rho)\,V(x)$. 
Moreover, 
$$
\wt{g}(x)=\frac{1}{2}\,
\left(
\sum^{p}_\zs{j=1}(\a_\zs{1,j}-\a_\zs{0,j})\,x_\zs{j}
\right)^{2}
\le \frac{\sum^{p}_\zs{j=1}(\a_\zs{1,j}-\a_\zs{0,j})^{2}}{2}\,|x|^{2}\,.
$$
Therefore, choosing in \eqref{sec:A.8} $c=1 +[\sum^{p}_\zs{j=1}(\a_\zs{1,j}-\a_\zs{0,j})^{2}/2$, we obtain condition $(\C'_\zs{1})$. 

Condition ($\C_\zs{2}$) can be checked in the same way as in Example 2.

By Theorem~\ref{Th.sec:Mrk.2}, condition  $(\A_2)$ holds with $I=\EV \wt{g}(\varsigma)= \overline{\a}'\F\,\overline{\a}/2$,
where $\overline{\a}=\a_\zs{1}-\a_\zs{0}$ and the matrix $\F$ is defined in \eqref{erg_gaus.00},
and the SR procedure is asymptotically minimax. 

\section*{Acknowledgements}

The work of the first author was  partially supported by 
the Russian Science Foundation (research project No. 14-49-00079, National Research University ``MPEI", Moscow, Russia) and
by the Academic D.I. Mendeleev Fund Program of the  Tomsk State University (research project NU 8.1.55.2015 L). The work of the second author was supported in part by the U.S.\ Air Force Office of Scientific Research under MURI grant FA9550-10-1-0569, by the U.S.\  Defense Advanced Research Projects Agency under grant W911NF-12-1-0034 and by the U.S.\ Army Research Office under grants W911NF-13-1-0073 and W911NF-14-1-0246  at the University of Southern California, Department of Mathematics and at the University of Connecticut, Department of Statistics.

\appendix

\section{Auxiliary results}\label{A}
In this appendix, we present results needed for proofs in Section~\ref{sec:MaRe}.




\subsection{Asymptotic properties of the SR procedure for large threshold values}\label{sec:PrSt}

The following proposition establishes asymptotic properties of the SR procedure for large $h$ regardless of the optimality criteria. While it is being used in the proofs of asymptotic
optimality of the SR procedure under considered criteria, it also interesting independently.

\begin{proposition} \label{Prop: A1}
Let $T(h)$ be the SR procedure defined in \eqref{sec:Prbf.7-00}. 

\noindent {\rm \bf (i)} Assume that there exists a positive and finite number $I$ such that, for all $\varepsilon >0$, the following conditions hold:
\begin{equation}\label{Pmax}
\lim_{M\to\infty} \Pb_0\brc{\frac{1}{M} \max_{1\le n \le M} Z_{n}^0 \ge I (1+\varepsilon)} = 0 
\end{equation}
and, for some $r\ge 1$,
\begin{equation}\label{rcompLeft}
\sum_{n=1}^\infty \, n^{r-1} \, \sup_{\nu \ge 0} \Pb_\nu\brc{\frac{1}{n} Z_{\nu+n}^\nu < I  - \varepsilon} <\infty \, .
\end{equation}
Then
\begin{equation}\label{AAAproxa}
\limsup_{h\to\infty}\, \frac{1}{(\log h)^r} \, \EV_{\nu}[(T(h)-\nu)^+]^r \le \frac{1}{I^r} \quad \text{for all}~ \nu \ge 0,
\end{equation}
\begin{equation}\label{AAAproxb}
\limsup_{h\to\infty}\, \frac{1}{(\log h)^r} \, \EV_{\nu}[(T(h)-\nu)^r|T(h)>\nu] \le \frac{1}{I^r} \quad \text{for all}~ \nu \ge 0,
\end{equation}
and
\begin{equation}\label{AAAprox}
\lim_{h\to\infty}\, \frac{1}{(\log h)^r} \, \sup_{\nu\ge 0}\, \EV_{\nu}[(T(h)-\nu)^+]^r = \frac{1}{I^r} \, .
\end{equation}
\noindent {\rm \bf (ii)} Asymptotic relations \eqref{AAAproxa}, \eqref{AAAproxb} and \eqref{AAAprox} hold if
\begin{equation}\label{rcompA}
\sum_{n=1}^\infty \, n^{r-1} \, \sup_{\nu \ge 0} \Pb_\nu\brc{\abs{\frac{1}{n} Z_{\nu+n}^\nu-I} > \varepsilon} <\infty \quad \text{for all}~\varepsilon >0 .
\end{equation}
\noindent {\rm \bf (iii)} If, in particular, $r=1$, then the uniform complete convergence condition $(\A_\zs{2})$ implies \eqref{AAAproxa}, \eqref{AAAproxb} and \eqref{AAAprox} with $r=1$.
\end{proposition}

\proof
(i) Let $n_0=1+\lfloor\log h /(I-\varepsilon)\rfloor$. We have 
\begin{align}\label{IneqExpTb}
\EV_{\nu}\brcs{(T(h)-\nu)^+}^r  & = \int_0^\infty r t^{r-1} \Pb_\nu\brc{T(h)-\nu > t} \, \d t   \nonumber
\\
& = n_0^r + \sum_{n=0}^{\infty} \int_{n_0+n}^{n_0+n+1} r t^{r-1}  \Pb_\nu (T(h)-\nu > t) \, \d t \nonumber
\\
& = n_0^r + \sum_{n=0}^{\infty} \int_{n_0+n}^{n_0+n+1} r t^{r-1}  \Pb_\nu (T(h)-\nu > n_0+n) \, \d t \nonumber
\\
& = n_0^r + \sum_{n=0}^{\infty} [(n_0+n+1)^r- (n_0+n)^r ] \Pb_\nu (T(h)-\nu > n_0+n) \nonumber
\\
& = n_0^r + \sum_{n=n_0}^{\infty} [(n+1)^r-n^r]  \Pb_\nu (T(h)-\nu > n) \nonumber
\\
 & \le  n_0^{r} +\sum_{n=n_0}^{\infty}   r (n+1)^{r-1}   \Pb_\nu (T(h) > \nu+ n) \nonumber
 \\
 & \le  n_0^{r} +\sum_{n=n_0}^{\infty}   r 2^{r-1} n^{r-1}   \Pb_\nu (T(h) > \nu+ n) \, .
\end{align}
In just the same way as in \eqref{upperprobnu},  we obtain that for all $\nu \ge 0$ and $n\ge n_0$,
\[
\Pb_\nu\brc{T(h) >\nu+n} \le  \Pb_\nu\brc{Z_{\nu+n}^\nu/n < \log h/n} \le \Pb_\nu\brc{Z_{\nu+n}^\nu/n <I-\varepsilon}\, .
\]
Hence 
\[
 \sup_{\nu\ge 0} \EV_{\nu}\brcs{(T(h)-\nu)^+}^r \le  \brc{1+\frac{\log h}{I-\varepsilon}}^r + 
 r2^{r-1} \, \sum_{n=1}^\infty n^{r-1} \sup_{\nu \ge 0} \Pb_\nu\brc{\frac{1}{n} Z_{\nu+n}^\nu < I  - \varepsilon},
\]
where, by condition \eqref{rcompLeft}, the last term on the right-hand side is finite. This immediately implies the 
following upper bounds for the moments of the detection delay (for any $\nu\ge 0$)
\begin{equation}\label{AUppermax}
\limsup_{h\to\infty}\, \frac{1}{(\log h)^r} \, \EV_\nu\brcs{(T(h)-\nu)^+}^r \le\limsup_{h\to\infty}\, \frac{1}{(\log h)^r} \, \sup_\zs{\nu\ge 0}\, \EV_\nu\brcs{(T(h)-\nu)^+}^r \le \frac{1}{I^r} .
\end{equation}
The upper bound \eqref{AAAproxa} follows. To obtain the upper bound \eqref{AAAproxb} for the conditional risk, it suffices to observe that $ \EV_{\nu}[(T(h)-\nu)^r|T(h)>\nu] = \EV_{\nu}[(T(h)-\nu)^r]^+/\Pb_\infty(T(h)>\nu)$
and that $\Pb_\infty(T(h)>\nu) \ge 1-\nu/h\to1$ as $h\to\infty$ for every $\nu\ge 0$. The letter follows easily from the fact that $R_n-n$ is  a zero-mean $\Pb_\infty$-martingale.

Define $M_{\varepsilon,h}=(1-\varepsilon) \log h/(I+d)$, where as before $d=-\log(1-\varrho)$. 
Replacing $\alpha$ in the proof of Theorem~\ref{Th.sec:Bay.1}  by $1/h$, in particular in \eqref{sec:Bay.4}, we obtain that, for any $0<\varepsilon<1$,
\begin{equation} \label{gamma0}
\Pb_0(T(h) \le M_{\varepsilon,h})\le e^{(1+\varepsilon)I  M_{\varepsilon,h}} \, \Pb_{\infty} \brc{T(h) \le M_{\varepsilon,h}}
 +\Pb_0 \brc{\frac{1}{M_{\varepsilon,h}} \max_{1\le n \le M_{\varepsilon,h}} Z_{n}^{0}\ge (1+\varepsilon) I} .
\end{equation}
By Lemma~\ref{Lem:PFASR}, $\PFA(T(h)) \le 1/(h\varrho+1) := \alpha_\zs{*}$,
and as in \eqref{sec:Bay.5}, we have $\Pb_{\infty}\brc{T(h) \le \ell} \le \alpha_\zs{*}\,(1-\varrho)^{-\ell}$. Hence,
$$
e^{(1+\varepsilon)I  M_{\varepsilon,h}} \, \Pb_{\infty} \brc{T(h) \le M_{\varepsilon,h}} 
\le e^{(1+\varepsilon)I  M_{\varepsilon,h}+d M_{\varepsilon,h} + \log \alpha_\zs{*}}
\le h\alpha_\zs{*}\, h^{-\varepsilon^2}\le\, \varrho^{-1} h^{-\varepsilon^2} \, ,
$$
so that the first term in \eqref{gamma0} goes to zero as $h\to\infty$ for any $\varepsilon >0$ and for any $0<\varrho<1$.
By condition \eqref{Pmax}, the second term also goes to zero  as $h\to\infty$, and therefore,
$\lim_{h\to\infty} \Pb_0(T(h) > M_{\varepsilon,h}) = 1$ for any $0<\varepsilon<1$ and any $0<\varrho<1$.
 Finally, Chebyshev's inequality yields
$$
\sup_{\nu \ge 0} \EV_{\nu}[(T(h)-\nu)^+]^r\ge \EV_0 [T(h)]^r \ge M_{\varepsilon,h}^r\,\Pb_0\brc{T(h)>M_\zs{\varepsilon,h}},
$$
so 
$$
\liminf_\zs{h\to\infty}\,\frac{1}{(\log h)^r}
\sup_{\nu \ge 0} \EV_{\nu}[(T(h)-\nu)^+]^r \ge \brc{\frac{1-\varepsilon}{I+d}}^r
$$
for arbitrary $0<\varepsilon <1$ and $0<\varrho<1$, and we obtain the asymptotic lower bound
\begin{equation*}
\lim_\zs{h\to\infty}\, \frac{1}{(\log h)^r} \, \sup_{\nu\ge 0}\, \EV_{\nu}[(T(h)-\nu)^+]^r \ge \frac{1}{I^r} \, ,
\end{equation*}
which along with the upper bound \eqref{AUppermax} completes the proof of \eqref{AAAprox} in (i). 

(ii) The uniform $r-$complete convergence condition \eqref{rcompA} implies both conditions \eqref{Pmax} and \eqref{rcompLeft}, and hence,  \eqref{AAAproxa}, \eqref{AAAproxb} and \eqref{AAAprox} 
hold true under \eqref{rcompA}. 

(iii) Finally, when $r=1$, condition \eqref{rcompA} is nothing but the uniform complete convergence condition $(\A_\zs{2})$, and hence, \eqref{AAAproxa}, \eqref{AAAproxb} and \eqref{AAAprox}
 hold true with $r=1$ under $(\A_\zs{2})$. This completes the proof of all three assertions.
\endproof

\subsection{Auxiliary results for proving asymptotic optimality of the SR procedure}
\label{sec:Absrsk}

The following proposition allow us to compare the classes \eqref{sec:Prbf.3} and \eqref{sec:Prbf.2}.

\begin{proposition} \label{Pr.sec:Absrsk.1} 
For any $0<\beta<1$, $m^{*}\ge |\log (1-\beta)|/[|\log (1-\varrho_\zs{1,\beta})|]- 1$ and $k^{*}> m^{*}$,  the following inclusions hold:
\begin{equation}\label{sec:Absrsk.2}
\Delta(\alpha_\zs{2},\varrho_\zs{2,\beta}) \subseteq \Hc\left(\beta,k^{*},m^{*}\right) \subseteq \Delta(\alpha_\zs{1},\varrho_\zs{1,\beta})\,.
\end{equation}
\end{proposition}
\proof Let $\tau\in \Delta(\alpha_\zs{2},\varrho_\zs{2,\beta})$.
Taking in \eqref{sec:Bay.5} $l=k^{*}$, we obtain
\begin{equation*}
\sup_\zs{1\le k\le k^{*}-m^{*}}\,
\Pb_\zs{\infty}\left(k\le \tau< k+m^{*}\right)\le 
\Pb_\zs{\infty}\left(\tau< k^{*}\right)\le 
\alpha_\zs{2}\,(1-\varrho_\zs{2,\beta})^{-k^{*}}
=\beta\,.
\end{equation*}
Hence, $\tau \in \Hc(\beta,k^{*},m^{*})$, i.e., the first inclusion in \eqref{sec:Absrsk.2}  follows.

Now, if $\tau\in  \Hc\left(\beta,k^{*},m^{*}\right)$, then 
$\Pb_\zs{\infty}\left(\tau < 1+m^{*}\right)\le \beta$. Therefore, 
taking in \eqref{sec:Prbf.1} $\varrho=\varrho_\zs{1,\beta}$, we obtain
\begin{align*}
\sum_\zs{k= 1}^\infty \,\pi_\zs{k}(\varrho_\zs{1,\beta})\,
\Pb_\zs{\infty}\left(\tau < k\right)
&= \sum^{m^{*}+1}_\zs{k=1}\,\pi_\zs{k}(\varrho_\zs{1,\beta})\,
\Pb_\zs{\infty}\left(\tau < k\right)+
\sum^{\infty}_\zs{k=m^{*}+2}\,\pi_\zs{k}(\varrho_\zs{1,\beta})\,
\Pb_\zs{\infty}\left(\tau < k\right)\\
&\le \beta+\sum^{\infty}_\zs{k= m^{*}+2}\,\pi_\zs{k}(\varrho_\zs{1,\beta})
=
 \beta+\left(1-\varrho_\zs{1,\beta}\right)^{m^{*}+1}
=\alpha_\zs{1}\,,
\end{align*}
i.e., $\tau\in \Delta(\alpha_\zs{1},\varrho_\zs{1,\beta})$.
Hence both inclusions in \eqref{sec:Absrsk.2} are proven. \endproof

The following proposition allows us to compare the classes  \eqref{sec:Prbf.4} and \eqref{sec:Prbf.2}.

\begin{proposition} \label{Pr.sec:Cnrsk.1} 
For any $0<\beta<1$, 
$m^{*}\ge |\log (1-\beta)|/[|\log (1-\varrho_\zs{1,\beta})|]- 1$ and $k^{*}\ge m^{*}$, the following inclusions hold:  
\begin{equation}\label{sec:Cnrsk.2}
\Delta(\alpha_\zs{3},\varrho_\zs{2,\beta})
\subseteq \Hc^{*}\left(\beta,k^{*},m^{*}\right) \subseteq
\Delta(\alpha_\zs{1},\varrho_\zs{1,\beta}) \, .
\end{equation}
\end{proposition}

\proof
First we show the left inclusion. Let $\tau\in \Delta(\alpha_\zs{3},\varrho_\zs{2,\beta})$.
Then, taking into account the inequality \eqref{sec:Bay.5}
for any $l=k^{*}$ and using the definition of $\alpha_\zs{3}$, we obtain that
\begin{align*}
\sup_\zs{1\le k\le k^{*}-m^{*}}\,
 \Pb_\zs{\infty}(\tau < k+m^{*} | \tau \ge k)&\le
 \sup_\zs{1\le k\le k^{*}-m^{*}}\,
\frac{\Pb_\zs{\infty}\left(\tau < k+m^{*}\right)}
{\Pb_\zs{\infty}\left( \tau \ge k\right)}\\[2mm]
&\le \frac{\Pb_\zs{\infty}\left(\tau < k^{*}\right)}
{1-\Pb_\zs{\infty}\left( \tau < k^{*}\right)}\\[2mm]
&\le 
\frac{\alpha_\zs{3}\,(1-\varrho_\zs{2,\beta})^{-k^{*}}}
{1-\alpha_\zs{3}\,(1-\varrho_\zs{2,\beta})^{-k^{*}}}
=\beta
\,,
\end{align*}
i.e., $\tau$ belongs to $\Hc(\beta, k^{*},m^{*})$. 

Now we show the right inclusion 
in \eqref{sec:Absrsk.6}. Let $\tau$ be from $\Hc(\beta, k^{*},m^{*})$. Then, using
the definition of the class $\Hc(\beta, k^{*},m^{*})$ in \eqref{sec:Absrsk.2},
we obtain that $\Pb_\zs{\infty}(\tau < m^{*}+1) \le \beta$. Therefore,
similarly to the proof of the right inclusion in
\eqref{sec:Absrsk.2} we obtain that $\tau\in \Delta(\alpha_\zs{1},\varrho_\zs{1,\beta})$, and the proof is complete.
\endproof

\subsection{Auxiliary results for Example~1, Subsection~\ref{ssec:LaiCn}}
\label{sec:A_Ex_1}

Recall that $\kappa(x)$ and $\varpi_j$ are defined in \eqref{defkappavarpi}.

\begin{lemma} \label{Le.sec:A_Ex_0} 
For any $\varepsilon>0$, 
\begin{equation}\label{sec:A_Ex_1.0}
\lim_\zs{n\to\infty}\,\sup_\zs{x\in\bbr^{2}}\,\Pb_\zs{0}\left(\left\vert\sum_\zs{j=1}^{n}\, \varpi_\zs{j}\,\right\vert>\varepsilon n\,\vert\,X_\zs{0}=x\right)=0 \,.
\end{equation}
\end{lemma}

\proof
Indeed, we have 
$$
\EV_\zs{0}\left(\left(\sum^{n}_\zs{j=1}\,\varpi_\zs{j} \right)^{2}\vert X_\zs{0}=x\right)=\sum^{n}_\zs{j=1}\,
\EV_\zs{0}\left(\,\varpi^{2}_\zs{j} \vert X_\zs{0}=x \right)=\EV_\zs{0}\left(\sum^{n}_\zs{j=1}\,\EV_\zs{0} \left(\,\varpi^{2}_\zs{j}
\vert \Fc_\zs{j-1}\right) \vert X_\zs{0}=x
\right)\,,
$$
where $\Fc_\zs{j}=\sigma\{X_\zs{1},\ldots,X_\zs{j}\}$. Using the definition of the sequence $\varpi_\zs{j}$ in \eqref{defkappavarpi}, we obtain that 
$$
\EV_\zs{0} \left(\,\varpi^{2}_\zs{j}\vert \Fc_\zs{j-1}\right) =
X'_\zs{j-1}\Lambda\,G^{-1}(X_\zs{j-1})\Lambda\,X_\zs{j-1} \,\le\,\frac{\sigma^{2}_\zs{1}}{\lambda^{2}_\zs{1}}+\frac{\sigma^{2}_\zs{2}}{\lambda^{2}_\zs{2}}\,.
$$
Thus, $\EV_\zs{1}[ (\sum^{n}_{j=1} \varpi_{j})^{2}\vert X_{0}=X] \le \sigma^{2}_1/\lambda^{2}_1+\sigma^{2}_2/\lambda^{2}_2$,
which implies \eqref{sec:A_Ex_1.0}.
\endproof

\begin{lemma} \label{Le.sec:A_Ex_1} 
Assume that  in \eqref{sec:LaiCn.1} 
 the parameter  $\sigma^{2}_\zs{1}$ 
 is such that
 $
 \lim_\zs{\rho\to\infty}\,
 \sigma^{2}_\zs{1}\,\rho^{4}
 =0\,.$
Then
\begin{equation}\label{sec:A_Ex_1.1}
\lim_\zs{\rho\to\infty}\,
\EV\,\varkappa(\zeta)\,=+\infty\,.
\end{equation}
\end{lemma}

\proof
First note that 
 for any $x\in\bbr^{2}$
the inverse matrix  for \eqref{sec:LaiCn.4-0n0} can be written as
$$
\G^{-1}(x)
=
\frac{1}{\det G(x)}
\left(
\begin{array}{cc}
1+\sigma^{2}_\zs{2}\,x^{2}_\zs{2}
 &,\,-\rho\\[2mm]
-\rho &,\, 1+\rho^{2}+\sigma^{2}_\zs{1}\,x^{2}_\zs{1}
\end{array}
\right)
$$
and
$$
\det \G(x)=(1+\rho^{2}+\sigma^{2}_\zs{1}\,x^{2}_\zs{1})
(1+\sigma^{2}_\zs{2}\,x^{2}_\zs{2})
-\rho^{2}
=t_\zs{0}(x^{2}_\zs{2})+\sigma^{2}_\zs{1}\,x^{2}_\zs{1}\,
t_\zs{1}(x^{2}_\zs{2})
\,,
$$
where $t_\zs{0}(x)=1+\sigma^{2}_\zs{2}\,(1+\rho^{2})\,x$ and $t_\zs{1}(x)=1+\sigma^{2}_\zs{2}\,x$.
This function can be written as
\begin{equation}
\label{sec:LaiCn.5}
\varkappa(x_\zs{1},x_\zs{2})=
\frac{1}{2}\lambda^{2}_\zs{1}\,a(x)
-\rho\,\lambda_\zs{1}\lambda_\zs{2}\,\b(x)
+\frac{1}{2}\lambda^{2}_\zs{2}\,c(x)
\,,
\end{equation}
where $a(x)=x^{2}_\zs{1}\,t_\zs{1}(x^{2}_\zs{2})/\det \G(x)$, $\b(x)=x_\zs{1}\,x_\zs{2}/\det \G(x)$, and
$c(x)=x^{2}_\zs{2}(1+\rho^{2}+\sigma^{2}_\zs{1}\,x^{2}_\zs{1})/\det \G(x)$.
Now to study the function \eqref{sec:LaiCn.5} we represent the coefficient $\a(x)$ as $\a(x)=\a_\zs{0}(x)-\sigma^{2}_\zs{1}\,\a_\zs{1}(x)$,
where $\a_\zs{0}(x)=x^{2}_\zs{1}\,t_\zs{3}(x^{2}_\zs{2})$, $\a_\zs{1}(x)=x^{4}_\zs{1}\,t_\zs{1}(x^{2}_\zs{2})t_\zs{3}(x^{2}_\zs{2})/[t_\zs{0}(x^{2}_\zs{2})+\sigma^{2}_\zs{1}\,x^{2}_\zs{1}\,
t_\zs{1}(x^{2}_\zs{2})]$ and $t_\zs{3}(x)=t_\zs{1}(x)/t_\zs{0}(x)$.
Taking into account that $t_\zs{1}(x)\le t_\zs{0}(x)$ 
and that the random variable $\zeta_\zs{1}$ is $\Gc$-conditionally Gaussian with the parameters $0$ and 
$(1+\rho^{2})\varsigma_\zs{11}$, we obtain 
\begin{equation*}
\EV\,\a_\zs{1}(\zeta)\le \EV\,\zeta^{4}_\zs{1}
=\,3\,(1+\rho^{2})^{2}\EV\,\varsigma^{2}_\zs{11}
\,,
\end{equation*}
where $\varsigma_\zs{11}
=\,\sum^{\infty}_\zs{k= 1}\,\sigma^{2(k-1)}_\zs{1}\,\prod^{k-1}_\zs{l=1}\,\eta^{2}_\zs{1,l}$. We recall that
$\Gc=\sigma\{\eta_\zs{1,k},\eta_\zs{2,k}\,, k\ge 1\}$.
Now, by the Bunyakovsky--Cauchy--Schwarz  inequality, 
$$
\EV\,
\varsigma^{2}_\zs{11}\,\le\,
\frac{1}{1-\sigma^{2}_\zs{1}}\,
\sum_\zs{k\ge 1}\,\sigma^{2(k-1)}_\zs{1}\,\EV\,\prod^{k-1}_\zs{l=1}\,\eta^{4}_\zs{1,l}
=
\frac{1}{1-\sigma^{2}_\zs{1}}\,
\sum_\zs{k\ge 1}\,\left(\sigma^{2}_\zs{1}\,\EV\,\eta^{4}_\zs{1,1}\right)^{k-1}
=
\frac{1}{(1-\sigma^{2}_\zs{1})(1-3\sigma^{2}_\zs{1})}\,
\,.
$$
Therefore, for any $0<\sigma_\zs{0}<1/3$,
\begin{equation*}
\limsup_\zs{\rho\to\infty}\,\frac{1}{(1+\rho^{2})^{2}}\,
\sup_\zs{0<\sigma^{2}_\zs{1}\le \sigma_\zs{0}}
\EV\,\vert\a_\zs{1}(\zeta)\vert
\,<\infty
\,.
\end{equation*}
Now we  calculate the expectation  $\EV\a_\zs{0}(\zeta)$. To this end,  we set
$\check{\zeta}=\zeta_\zs{1}-\check{\kappa}\,\zeta_\zs{2}$ and $\check{\kappa}=\EV(\zeta_\zs{1}\zeta_\zs{2}\vert\Gc)/\EV(\zeta^{2}_\zs{2}\vert\Gc) =\rho(\varsigma_\zs{12}/\varsigma_\zs{22})$.
Conditioned on $\Gc$, the random variable $\check{\zeta}$ is independent of $\zeta_\zs{2}$, and  $\check{\zeta}$  is Gaussian with the parameters 
$\left(0,\EV\,(\check{\zeta}^{2}\vert\Gc)\right)$, where
\begin{equation*}
\EV\,(\check{\zeta}^{2}\vert\Gc)
=\EV\,(\zeta^{2}_\zs{1}\vert\Gc)-\frac{\left(\EV\,(\zeta_\zs{1}\zeta_\zs{2}\vert\Gc)\right)^{2}}{\EV\,(\zeta^{2}_\zs{2}\vert\Gc)}
=(1+\rho^{2})\varsigma_\zs{11}
-
\rho^{2}\frac{\varsigma^{2}_\zs{12}}{\varsigma_\zs{22}}
:=1+\rho^{2}\varsigma_\zs{*}\,.
\end{equation*}
By the  Bunyakovsky--Cauchy--Schwarz  inequality, the random variable $\varsigma_\zs{*}\ge 0$ a.s. Next, 
using the definitions of the random variables $\varsigma_\zs{ij}$ in
\eqref{sec:LaiCn.4-0}, we obtain that $\varsigma_\zs{*}=0$ if and only if for any $k\ge 1$
$$
\sigma^{k-1}_\zs{1}\,\prod^{k-1}_\zs{l=1}\,\eta_\zs{1,l}
=\sigma^{k-1}_\zs{2}\,\prod^{k-1}_\zs{l=1}\,\eta_\zs{2,l}
\,.
$$
So, $\varsigma_\zs{*}>0$ a.s. Thus,
\begin{align*}
\EV(\a_\zs{0}(\zeta)\vert\Gc)&=
\EV(\check{\zeta}^{2}\vert\Gc)\,\EV\left(t_\zs{3}(\zeta^{2}_\zs{2})\vert\Gc\right)
+\check{\kappa}^{2}\EV\left(\zeta^{2}_\zs{2}t_\zs{3}(\zeta^{2}_\zs{2})\vert\Gc\right)\\[2mm]
&\ge (1
+\rho^{2}\,\varsigma_\zs{*})\,
 \EV\left(\frac{1}{1+\sigma^{2}_\zs{2}(1+\rho^{2})\zeta^{2}_\zs{2}}\vert\Gc\right)\,,
\end{align*}
and we obtain that
$$
\liminf_\zs{\rho\to\infty}\,\EV\left(\a_\zs{0}(\zeta)\vert\Gc\right)
\ge \,
\frac{1}{\sigma^{2}_\zs{2}}\,
\varsigma_\zs{*}\,
 \EV\left(\frac{1}{\zeta^{2}_\zs{2}}\vert\Gc\right)
=+\infty
\quad\mbox{a.s.,}
$$
i.e.,
$
\lim_\zs{\rho\to\infty}\,
\EV\,\a_\zs{0}(\zeta)\,=+\infty.
$
Setting now 
 $\Gc_\zs{1}=\sigma\{\xi_\zs{2,l}\,,\,\eta_\zs{2,l}\,,\,l\ge 1\}$, we obtain 
$$
\EV\,\left(\check{\kappa} \vert \Gc_\zs{1}\right)=
\rho
\frac{\EV\,\left(\varsigma_\zs{12}\vert \Gc_\zs{1}\right)}{\varsigma_\zs{22}}
=\rho
\sum^{\infty}_\zs{k= 1}\,\sigma^{k-1}_\zs{1}\,\sigma^{k-1}_\zs{2}\,
\frac{1}
{\varsigma_\zs{22}}
\EV\,\left(\prod^{k-1}_\zs{l=1}\,\eta_\zs{1,l}\eta_\zs{2,l}\vert \Gc_\zs{1}\right)
=\,
\frac{\rho}{\varsigma_\zs{22}}\,.
$$
Moreover,  $\b(x)=\b_\zs{0}(x)-\sigma^{2}_\zs{1}\,\b_\zs{1}(x)$,
where $\b_\zs{0}(x)=x_\zs{1}\,x_\zs{2}/t_\zs{0}(x^{2}_\zs{2})$ and $\b_\zs{1}(x)=x^{3}_\zs{1}\,x_\zs{2}t_\zs{3}(x^{2}_\zs{2})/[t_\zs{0}(x^{2}_\zs{2})+\sigma^{2}_\zs{1}\,x^{2}_\zs{1}
t_\zs{1}(x^{2}_\zs{2})]$.
Therefore, taking into account that
$\EV\,\left(\check{\zeta}\vert\Gc\right)=0$, we obtain 
$$
\EV\,\left(\b_\zs{0}(\zeta)\vert\Gc\right)
=\EV\left(\frac{\zeta_\zs{1}\,\zeta_\zs{2}}{t_\zs{0}(\zeta^{2}_\zs{2})}\vert\Gc\right)
=
\EV\,\left(\check{\kappa}\frac{\zeta^{2}_\zs{2}}{t_\zs{0}(\zeta^{2}_\zs{2})}\vert\Gc\right)\,.
$$
 This implies immediately that
$$
\EV\,\b_\zs{0}(\zeta)=
\EV\,\check{\kappa}\frac{\zeta^{2}_\zs{2}}{t_\zs{0}(\zeta^{2}_\zs{2})}
=
\EV\,\left(\frac{\zeta^{2}_\zs{2}}{t_\zs{0}(\zeta^{2}_\zs{2})}\EV\,\left(\check{\kappa} \vert \Gc_\zs{1}\right)\right)
=\,
\rho\,
\EV\,\frac{\zeta^{2}_\zs{2}}{t_\zs{0}(\zeta^{2}_\zs{2})\varsigma_\zs{22}}
\,.
$$
It is easy to check that
$$
\lim_\zs{\rho\to\infty}\,\rho\,\EV\,\b_\zs{0}(\zeta)=
\frac{1}{\sigma^{2}_\zs{2}}\,
\EV\,\frac{1}{\varsigma_\zs{22}}
\le\, 
\frac{1}{\sigma^{2}_\zs{2}}
\,.
$$
Taking into account here that $t_\zs{3}(x)\le 1$
we obtain that $\vert\,\b_\zs{1}(x)\,\vert\,\le\,\vert x_\zs{1}\vert^{3}/\sigma_\zs{2}$, i.e.,
$$
\EV\,\vert\b_\zs{1}(\zeta)\vert\le\,\frac{1}{\sigma_\zs{2}}\,\EV\vert\zeta_\zs{1}\vert^{3}
= \frac{2\sqrt{2}}{\sigma_\zs{2}\sqrt{\pi}}\,(1+\rho^{2})^{3/2}\,
\EV\,
\varsigma^{3/2}_\zs{11}
\le
\frac{2\sqrt{2}}{\sigma_\zs{2}\sqrt{\pi}}\,(1+\rho^{2})^{3/2}\,
\EV\,
\varsigma^{2}_\zs{11} 
\,.
$$
Therefore, for any $0<\sigma_\zs{0}<1/3$
\begin{equation*}
\limsup_\zs{\rho\to\infty}\,\frac{1}{(1+\rho^{2})^{3/2}}\,
\sup_\zs{0<\sigma^{2}_\zs{1}\le\sigma_\zs{0}}
\EV\,\vert\b_\zs{1}(\zeta)\vert\,
<\infty
\,.
\end{equation*}
Clearly, 
$
c(x)\le \sigma^{-2}_\zs{2}.
$
Thus, using the condition of this lemma we obtain
\eqref{sec:A_Ex_1.1}.
\endproof

\section{Auxiliary non-asymptotic bounds for the concentration inequalities}\label{B}

\subsection{Correlation inequality}

We now give the important correlation inequality proved in \cite{GaPe13}.

\begin{proposition}\label{Pr.sec:Bi.1}
 Let $(\Omega,\Fc,(\Fc_\zs{j})_\zs{1\le j\le n},\Pb)$ be a filtered probability space and
$(u_\zs{j}, \Fc_\zs{j})_\zs{1\le j\le n}$ be a sequence of random variables such that, for some $\r\ge 2$, $\max_{1\le j\le n}\,\EV\,|u_{j}|^{\r}<\infty$.
Define
\begin{equation}\label{subsec:AMc.1}
b_\zs{j,n}(\r)=
\left(
\EV\,
\left(
|u_\zs{j}|\,\sum^{n}_\zs{k=j}
|\EV\,(u_\zs{k}|\Fc_\zs{j})|
\right)^{\r/2}
\right)^{2/\r}\,.
\end{equation}
 Then
$$
\EV\, \Big | \sum^{n}_\zs{j=1}\,u_\zs{j} \Big |^{\r} \le\,(2\r)^{\r/2} \left(\sum^{n}_\zs{j=1}\,b_\zs{j,n}(\r)\right)^{\r/2}\,.
$$
\end{proposition}

\subsection{Geometric ergodicity for homogeneous Markov processes}\label{subsec:AMc}
We follow the Meyn--Tweedie approach \cite{MeTw93}. 
We recall some definitions from \cite{MeTw93} and \cite{GaPe14}
for a  homogeneous Markov process $(X_\zs{n})_\zs{n\ge 0}$ defined on a measurable state space
$(\Xc, \Bc(\Xc))$. Denote by $P(x,\cdot)\,, x\in\Xc\,,$ the transition probability 
of this process, i.e., for any $A\in\Bc(\Xc), x\in\Xc$,
\begin{equation}\label{subsec:AMc.1-1}
P(x,A)\,=\,\Pb_\zs{x}(X_\zs{1}\in A)
=
\Pb(X_\zs{1}\in A|X_\zs{0}\,=\,x)\,.
\end{equation}
The $n-$step transition probability is $P^{n}(x,A)\,=\,\Pb_\zs{x}(X_\zs{n}\in A)$.

We recall that a measure $\lambda$ on $\Bc(\Xc))$ is called 
{\it invariant} (or {\em stationary} or {\em ergodic}) for this process if,
for any $A\in\Bc(\Xc)$, 
\begin{equation}\label{subsec:AMc.1-2}
\lambda(A)\,=\,\int_{\Xc}\,P(x,A)\lambda(\d x)\,.
\end{equation}

\noindent If there exists an invariant positive measure 
$\lambda$ with $\lambda(\Xc)\,=\,1$ then the process is
called {\it positive}.

Assume that the process $(X_\zs{n})_\zs{n\ge 0}$ satisfies the following {\em minorization} condition: \vspace{2mm}

\noindent $(\D_\zs{1})$
{\em
There exist  $\delta>0$, a set $C\in\Bc(\Xc)$ and a probability measure $\varsigma$ on $\Bc(\Xc)$ with
$\varsigma(C)=1$, such that for any $A\in \Bc(\Xc)$, for which $\varsigma(A)>0$, $\inf_\zs{x\in C}P(x,A)>\delta\,\varsigma(A)$.
}

\noindent 
Obviously, this condition implies that
$
\eta\,=\,\inf_\zs{x \in C}\,
P(x,C)-\delta>0.
$

Now we impose the {\em drift} condition. \vspace{2mm}

\noindent $(\D_\zs{2})$
 {\em
There exist
a $\Xc\to [1,\infty)$ function $\V$, constants $0<\rho<1$,  $D\ge 1$ and
 a set $C$ from $\Bc(\Xc)$ such that $\V^{*}=\sup_\zs{x\in C}|\V(x)|<\infty$ and, for all $x\in\Xc$},
\begin{equation}\label{subsec:AMc.2}
\EV_\zs{x}\left(\V(X_\zs{1})\right)\,
\le\,(1-\rho)\V(x)\,+\,D\Ind{C}(x)\,.
\end{equation}
In this case, we call $\V$ the {\em Lyapunov function.}

In this paper, we use the following theorem from \cite{GaPe14}.

\begin{theorem} \label{Th.AMc.1}
Let $(X_\zs{n})_\zs{n\ge 0}$ be a homogeneous  Markov process
satisfying conditions $(\D_\zs{1})$ and $(\D_\zs{2})$ with the same set 
$C\in\Bc(\Xc)$. Then $(X_\zs{n})_\zs{n\ge 0}$ is a positive geometric ergodic process, i.e.,
\begin{equation}\label{sec:AMc.3}
\sup_\zs{n\ge 0}\,
e^{\kappa^{*} n}\,
\sup_\zs{x\in\Xc}\,
\sup_\zs{0\le g\le \V}
\frac{1}{\V(x)}
|
\EV_\zs{x}\,g(X_\zs{n})
-
\lambda(\wt{g})
|
\le R^{*}
\end{equation}
for some positive constants $\kappa^{*}$ and  $R^{*}$ which are given in {\rm \cite{GaPe14}}.
\end{theorem}


\end{document}